%% file: main.tex
\documentclass{article}

\usepackage{arxiv}

\usepackage[utf8]{inputenc} 
\usepackage[T1]{fontenc}    
\usepackage{hyperref}       
\usepackage{url}            
\usepackage{booktabs}       
\usepackage{amsfonts}       
\usepackage{nicefrac}       
\usepackage{microtype}      
\usepackage{doi}
\usepackage{style}
 
\title{The lowest-order Neural Approximated Virtual Element Method on polygonal elements}

\author{ \href{https://orcid.org/0000-0001-8642-4258}{\includegraphics[scale=0.06]{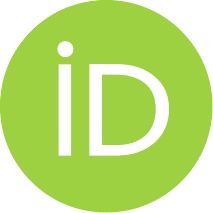}\hspace{1mm}Stefano~Berrone} \\
	Dipartimento di Scienze Matematiche\\
	``G. L. Lagrange''\\
	Politecnico di Torino, TO, 10129 \\
	\texttt{stefano.berrone@polito.it} \\
	\And
	\href{https://orcid.org/0000-0002-2837-2792}{\includegraphics[scale=0.06]{orcid.pdf}\hspace{1mm}Moreno~Pintore} \\
	Inria, \\
	Laboratoire Jacques-Louis Lions, \\
	Sorbonne Universit\'e, \\
	4 place Jussieu, 75005 Paris, France \\
	\texttt{moreno.pintore@inria.fr} \\
	\And
	\href{https://orcid.org/0000-0002-8540-3639}{\includegraphics[scale=0.06]{orcid.pdf}\hspace{1mm}Gioana~Teora} \\
	Dipartimento di Scienze Matematiche\\
	``G. L. Lagrange''\\
	Politecnico di Torino, TO, 10129 \\
	\texttt{gioana.teora@polito.it} \\
}

\renewcommand{\shorttitle}{NAVEM}

\hypersetup{
pdftitle={The lowest-order Neural Approximated Virtual Element Method on polygonal elements},
pdfsubject={q-bio.NC, q-bio.QM},
pdfauthor={Stefano~Berrone, Moreno~Pintore, Gioana~Teora},
pdfkeywords={Virtual Element Method, Neural Network, Basis Functions, Polygonal Meshes},
}

\begin{document}
\maketitle

\begin{abstract}
The lowest-order Neural Approximated Virtual Element Method on polygonal elements is proposed here. This method employs a neural network to locally approximate the Virtual Element basis functions, thereby eliminating issues concerning stabilization and projection operators, which are the key components of the standard Virtual Element Method. We propose different training strategies for the neural network training, each correlated by the theoretical justification and with a different level of accuracy. Several numerical experiments are proposed to validate our procedure on general polygonal meshes and demonstrate the advantages of the proposed method across different problem formulations, particularly in cases where the heavy usage of projection and stabilization terms may represent challenges for the standard version of the method.
Particular attention is reserved to triangular meshes with hanging nodes which assume a central role in many virtual element applications.
\end{abstract}


\input{sections}

\section*{Acknowledgements}
The author S.B. kindly acknowledges partial financial support provided by PRIN project ``Advanced polyhedral discretisations of heterogeneous PDEs for multiphysics problems'' (No. 20204LN5N5\_003), by PNRR M4C2 project of CN00000013 National Centre for HPC, Big Data and Quantum Computing (HPC) (CUP: E13C22000990001) and the funding by the European Union through project Next Generation EU, M4C2, PRIN 2022 PNRR project P2022BH5CB\_001 ``Polyhedral Galerkin methods for engineering applications to improve disaster risk forecast and management: stabilization-free operator-preserving methods and optimal stabilization methods.''. The author G.T. kindly acknowledges the financial support provided by the MIUR programme ``Programma Operativo Nazionale Ricerca e Innovazione 2014 - 2020'' 
~ (CUP: E11B21006490005) and by INdAM - GNCS Project CUP\_E53C23001670001.

\bibliographystyle{IEEEtran}
\bibliography{biblio.bib}

\end{document}

%% file: sections.tex
\section{Introduction}
The Virtual Element Method (VEM in short), introduced in \cite{LBe13} for the Laplace problem and then extended to the general second-order elliptic problem in \cite{LBe16}, can be considered as a generalization of the Finite Element Method (FEM) which introduces in the local space suitable non-polynomial functions, as well as the standard polynomials. The introduction of these non-polynomial functions, which are not required in a closed form, allows to work with polytopal elements in a very simple way preserving the polynomial accuracy.
These Virtual Element functions are solutions to local PDE problems inside each element of the tessellation that are actually never solved explicitly, neither exactly nor approximately. Since these functions are not explicitly known inside the elements, the discrete bilinear form used in the VEM discretization of the problem is just an approximation of the continuous counterpart that exploits some computable polynomial projections of VE functions to access their point-wise evaluation. Indeed, the core idea of the standard VEM method is to define suitable local spaces and degrees of freedom that allow to exactly compute the entries of the stiffness bilinear form when at least one of the two entries is a polynomial. The remaining entries, which account for the non-polynomial part, are replaced by a stabilization term to produce results that are of the right order of magnitude and satisfy stability properties. However, a unique prescription of this stabilization term is not provided from the virtual element theory and its selection is mainly guided by numerical experiments, becoming highly problem-dependent \cite{Russo2024}. Furthermore, the presence of the stabilization term can limit the accuracy of the method in case of strongly anisotropic problems due to its intrinsic isotropic nature \cite{MarconTeora2023, Teora2023}. Finally, the need to introduce some polynomial projectors to access the point-wise evaluation of VE functions can represent a limitation in the post-processing phase and may induce many issues also in complex non-linear problems \cite{Cangiani2019, Adak2019}.

Recently, various efforts have been made to address these limitations. In \cite{MarconTeora2023, Marcon2024}, the first stabilization-free methods have been proposed in which the discrete bilinear forms only involve polynomial projections on enhanced polynomial spaces, whose polynomial degrees strongly depend on the geometry of the underlying polygons. In \cite{Credali2024}, a reduced basis method is proposed to cheaply reconstruct approximations of VE basis functions which could be exploited to properly design stabilization terms or for post-processing of the solution. In \cite{TrezziZerbinati2024}, a lightning virtual element method is developed which actually computes the VE basis functions by solving a PDE problem on each element with the iterative \textit{Laplace Solver} proposed in \cite{Gopal2019}. Lastly, in \cite{PintoreTeora2024}, the Neural Approximated Virtual Element Method (NAVEM in short), which employs the neural network to approximate the VE basis functions, has been briefly presented and tested for the case of quadrilateral elements. In the context of the latter two methods, it is clear that the usage of ``virtual'' term just refers to the underlying local space, since the local construction of VE basis functions allows to get rid of any stabilization terms or polynomial projectors which represent the main features of the Virtual Element Method.

In the last few years, thanks to the availability of easily customizable machine learning libraries like Tensorflow \cite{tensorflow2015-whitepaper}, Pytorch \cite{pytorch} or JAX \cite{jax2018github}, numerous novel machine learning enhanced numerical methods have been proposed. In the context of this new research field, known as Scientific Machine Learning (SciML) \cite{Cuomo2022}, we present the NAVEM method on general polygonal meshes. Initially introduced by the authors for solving the Laplace problem on quadrilateral meshes in \cite{PintoreTeora2024}, this method is now extended to polygons with more than four vertices, allowing also for the presence of hanging nodes. Inspired by \cite{MarconTeora2023,Gopal2019} and by the recent success of other SciML techniques, this method leverages the neural network to approximate VE basis functions as a linear combination of harmonic functions, segregating the main computational effort needed to compute such approximations to the offline stage. Indeed, in the online phase, since the need for computing the local projection matrices and defining a stability operator is circumvented, NAVEM acts like a FEM method on polygonal meshes.
Specifically, we propose an enhanced version of the original NAVEM method, modifying the local approximation spaces to include new harmonic functions in addition to the harmonic polynomials to better capture singularities near the vertices of the polygon, and refining the neural network architecture to reduce oscillations between interpolation points. We also explore several training strategies, each offering varying levels of accuracy and aimed at minimizing distinct loss functions, with theoretical justifications provided for each approach. Numerical experiments validate the viability of our procedure on different polygonal meshes and show the advantages of using this new procedure, especially when addressing highly non-linear problems.

The paper is structured as follows. Section \ref{sec:vem} briefly introduces the VEM formulation, which is essential for developing an appropriate architecture and training strategy for the neural network. The Neural Approximated method is presented in Section \ref{sec:navem}, while Section \ref{sec:neural_network} details the network architecture and training strategy. Finally, Section \ref{sec:numerical_results} presents various numerical experiments on polygonal meshes to demonstrate the method's performance, including its application to anisotropic and nonlinear problems, which can pose challenges for the standard procedure.

\section{The Model Problem and The Virtual Element Method}\label{sec:vem}

Let us now introduce some notations used throughout the paper. 
Given $k\in\N$, we denote by $\|\cdot \|_{\sob{k}{\genericset}}$ the norm in the Sobolev space $\sob{k}{\genericset}$ on some open subset $\genericset \subset \R^2$. Furthermore, we use the symbol $(\cdot,\cdot)_{\genericset}$ to denote both the scalar product in $\leb{2}{\genericset}$ and in $\leb{2}{\genericset} \times \leb{2}{\genericset}$. We recall that given two vector functions $\vv = \begin{bmatrix} v_1, v_2
\end{bmatrix}^T$ and $\ \uu = \begin{bmatrix} u_1, u_2
\end{bmatrix}^T$, the scalar product in $\leb{2}{\genericset} \times \leb{2}{\genericset}$ is defined as
\begin{linenomath}
\begin{equation*}
\left(\vv,\uu\right)_{\genericset} = \int_{\genericset} (v_1 u_1 + v_2 u_2), \quad \|\vv \|_{\leb{2}{\genericset}} = \sqrt{\left(\vv,\vv\right)_{\genericset}}.
\end{equation*}
\end{linenomath}

Let us consider an open, bounded, convex polygonal domain $\Omega \subset \R^2$ with boundary $\Gamma$. Given $f \in \leb{2}{\Omega}$, we consider the following Poisson problem:
\begin{equation}
\begin{cases}
- \Delta u = f & \text{in } \Omega,\\
u = 0 & \text{on } \Gamma.
\end{cases}
\label{eq:laplaceproblem}
\end{equation}
The variational formulation of problem \eqref{eq:laplaceproblem} reads as: \textit{Find $u\in \V = \sob[0]{1}{\Omega}$ such that:}
\begin{equation}
\dbilin{u}{v} = \scal[\Omega]{f}{v} \quad \forall v \in \V,
\label{eq:laplace_variational_problem}
\end{equation}
where the bilinear form $\dbilin{}{} : \V \times \V \to \R$ is given by:
\begin{equation}
\dbilin{u}{v} = \scal[\Omega]{\nabla u}{\nabla v}\quad \forall u,v \in \V.
\label{eq:cont_bilinear_form}
\end{equation}

\subsection{The Virtual Element Space}

Let $\Th$ be a decomposition of $\Omega$ into polygons $E$ and let $\Eh$ be the set of edges of the elements in $\Th$. Furthermore, we denote by $\Nv[E]$ the number of vertices (and of edges), by $\Eh[E]$ the set of edges and by $h_E$ the diameter of the element $E \in \Th$. As usual, $h$ denotes the maximum diameter of the polygons in $\Th$. 
We assume that the following mesh assumptions hold true \cite{LBe16}.
\begin{assum}[Mesh assumptions]\label{ass:mesh}
There exists a positive constant $\rho$, independent of
$E$ and $h$, such that
\begin{itemize}
\item each polygon $E \in \Th$ is star-shaped with respect to a ball of radius $\geq \rho h_E$;
\item for each edge $e \in \Eh[E]$, it holds: $\vert e \vert \geq \rho h_E$.
\end{itemize}
\end{assum}

Given a polygon $E$, for each integer $k \geq 0$, we denote by $\Poly{k}{E}$ the set of two-dimensional polynomials of degree up to $k$ defined on $E$, of dimension $n_k = \dim \Poly{k}{E} = \frac{(k+1)(k+2)}{2}$. Furthermore, we introduce 
the set
\begin{linenomath}
\begin{equation*}
\Bk{1}{\partial E} = \left\{v \in C^0\left(\partial E\right) : v_{|e}\in \Poly{1}{e} \forall e \in \Eh[E]\right\},
\end{equation*}
\end{linenomath}
whose dimension is $\dim \Bk{1}{\partial E} = \Nv[E]$.
For all $E \in \Th$, we define the lowest-order local virtual element space \cite{LBe13} as the set
\begin{linenomath}
\begin{align}
 \Vh[E]{1} = \Big\{ v \in \sob{1}{E}:\quad &(i)\ \Delta v = 0, \quad (ii) \ v_{| \partial E} \in \Bk{1}{\partial E}\Big\}, 
 \label{eq:prop}
\end{align}
\end{linenomath}
with dimension $\Ndof[E] = \dim \Vh[E]{1} = \Nv[E]$, and we consider the value of $v_h \in \Vh[E]{1}$ at the vertices of $E$ as local degrees of freedom. 

The key property of the Virtual Element Method is that, thanks to this definition of the degrees of freedom, we are able to exactly (up to machine precision) compute the local projection $\proj{E,\nabla}{1} v_h$ of each VE function $v_h \in \Vh[E]{1}$, where the \textit{computable} local polynomial projector $\proj{E,\nabla}{1} : \sob{1}{E} \to \Poly{1}{E}$ is defined such that, for each $E\in \Th$,
\begin{linenomath}
\begin{equation*}
        \scal[E]{\nabla v_h - \nabla \proj{E, \nabla}{1} v_h}{\nabla p} = 0,\quad \forall p \in \Poly{1}{E} \text{ and } \int_{\partial E} \proj{E, \nabla}{1} v = \int_{\partial E} v.
    \end{equation*}
    \end{linenomath}
Finally, the global Virtual Element space is obtained by gluing together the local spaces as
\begin{linenomath}
\begin{equation*}
\Vh{1} = \Big\{ v \in \V \cap \con{0}{\overline{\Omega}}:\ v_{h|E} \in  \Vh[E]{1} \quad \forall E \in \Th \Big\}.
\end{equation*}
\end{linenomath}
\subsection{The Virtual Element Discretization}

Initially, we can observe that the continuous bilinear form \eqref{eq:cont_bilinear_form} can be split according to the tessellation $\Th$ as
\begin{linenomath}
\begin{equation*}
\dbilin{u}{v} = \sum_{E \in \Th} \dbilin[E]{u}{v},\quad \dbilin[E]{u}{v} = \scal[E]{\nabla u}{ \nabla v}\quad \forall u,v \in \V. 
\end{equation*}
\end{linenomath}
Then, we note that, in general, we are not able to compute the quantity
\begin{linenomath}
\begin{equation*}
\dbilin[E]{u_h}{v_h} = \scal[E]{\nabla u_h}{ \nabla v_h}\quad \forall u_h,\ v_h \in \Vh[E]{1}, 
\end{equation*}
\end{linenomath}
since we do not know the virtual element functions in a closed-form in the interior of each element $E\in \Th$.
To overcome this issue, the main idea of the Virtual Element Method is to substitute the continuous bilinear form with a \textit{computable} discrete counterpart $\dbilinh[E]{}{} : \Vh[E]{1} \times  \Vh[E]{1} \to \R$ which satisfies the two following properties \cite{LBe13}:
\begin{itemize}
\item \textit{Consistency}: For all $p \in \Poly{1}{E}$ and for all $v_h \in \Vh[E]{1}$
\begin{linenomath}
\begin{equation*}
\dbilinh[E]{p}{v_h} = \dbilin[E]{p}{v_h}.
\end{equation*}
\end{linenomath}
\item \textit{Stability}: There exist two positive constants $\alpha_{\ast},\ \alpha^{\ast}$ independent of $h$ such that
\begin{equation}
    \alpha_{\ast} \dbilin[E]{v}{v} \leq \dbilinh[E]{v}{v} \leq \alpha^{\ast} \dbilin[E]{v}{v},\quad \forall v \in \Vh[E]{1}:\ \proj{E,\nabla}{1}v = 0.
    \label{eq:stab_property}
\end{equation}
\end{itemize}
To build a discrete bilinear form which satisfies the consistency and stability properties, the local continuous bilinear form is first split as 
\begin{equation}
\dbilin[E]{u_h}{v_h} = \dbilin[E]{\proj{E,\nabla}{1}u_h}{\proj{E,\nabla}{1}v_h} + \dbilin[E]{(I-\proj{E,\nabla}{1})u_h}{(I-\proj{E,\nabla}{1})v_h},
\label{eq:cont_vem_bilinear_form}
\end{equation}
where the equality is due to the orthogonality of $\proj{E,\nabla}{1}$ with respect to the scalar product induced by $\dbilin[E]{}{}$. The first term in the right-hand side of \eqref{eq:cont_vem_bilinear_form} is computable thanks to the definition of the local degrees of freedom, whereas the second one could be approximated by any \textit{computable} symmetric positive definite bilinear form $\stab[E]{}{}$ that satisfies the stability property \eqref{eq:stab_property}.

Finally, it can be shown that the local discrete bilinear form
\begin{linenomath}
\begin{equation*}
\dbilinh[E]{u_h}{v_h} = \dbilin[E]{\proj{E,\nabla}{1}u_h}{\proj{E,\nabla}{1}v_h} + \stab[E]{(I-\proj{E,\nabla}{1})u_h}{(I-\proj{E,\nabla}{1})v_h}
\end{equation*}
\end{linenomath}
is computable and satisfies the consistency and the stability property \cite{LBe13}.

Now, let us define $\dof_i^E$, for each $i=1,\dots,\Ndof[E]$ and each $E\in \Th$, as the operator that associates with each sufficiently smooth function $\varphi$ its $i$-th local degree of freedom $\dof_i^E(\varphi)$. A standard choice for the stabilization term for the two-dimensional case is given by the \textit{dofi-dofi} stabilization term 
\begin{equation}
    \stab[E]{u_h}{v_h} = \sum_{i=1}^{\Ndof[E]} \dof_i^E(u_h) \dof_i^E(v_h).
    \label{eq:dofidofi}
\end{equation}
We observe that, when dealing with more general elliptic equations, this stabilization is usually pre-multiplied by a constant $C_{s}$, which accounts for the magnitude of the diffusion coefficients.
Other stabilization methods have been proposed in the literature, which may take integral forms \cite{LBe17} or be a variant of the dofi-dofi stabilization, such as the $D$-recipe version introduced in \cite{BEIRAODAVEIGA20171110}. The $D$-recipe form aims to prevent the stabilization from becoming too small in magnitude with respect to the consistency term when high-order methods are considered. 

Finally, the virtual element discretization of problem \eqref{eq:laplace_variational_problem} reads as: \textit{Find $u_h \in \Vh{1}$ such that:} 
\begin{linenomath}
\begin{equation}
\sum_{E \in \Th} \dbilinh[E]{u_h}{v_h} = \sum_{E \in \Th} \rightlinh[E]{v_h} \quad \forall v_h \in \Vh{1},
\label{eq:vem_problem_discretization}
\end{equation} 
\end{linenomath}
where $\rightlinh[E]{} : \Vh[E]{1} \to \R$ is the discrete version of the local forcing term given by:
\begin{linenomath}
\begin{equation*}
\rightlinh[E]{v_h} = \scal[E]{f}{\proj{\partial E, 0}{0} v_h},\quad \text{with}\quad \proj{\partial E, 0}{0} v_h = \int_{\partial E} v_h.
\end{equation*}
\end{linenomath}

\section{The Neural Approximated Virtual Element Method}\label{sec:navem}

Let us introduce the set of the VEM Lagrangian basis functions $\{\varphi_{i}\}_{i=1}^{\Ndof}$ corresponding to the aforementioned degrees of freedom, each of them associated with a different internal vertex $v_i$ of the tessellation $\Th$. We denote by $\mathbb{S}_i = \supp{\varphi_i} = \bigcup_{j=1}^{N_{v_i}} E_j$ the support of $\varphi_i$, i.e. the union of the $N_{v_i}$ elements $E_j \in \Th$ adjacent to the vertex $v_i$. Furthermore, given an element $E\in \Th$, for the sake of brevity, we denote by $\{\varphi_{j,E}\}_{j=1}^{\Ndof[E]}$ the set of the restrictions to $E$ of the Lagrangian basis functions related to the vertices of $E$. Clearly, the local and the global virtual element spaces can be written as 
\begin{linenomath}
\begin{equation*}
    \Vh[E]{1} = \myspan \{ \varphi_{j,E}: \ j = 1,\dots,\Ndof[E]\}
\end{equation*}
\end{linenomath}
and
\begin{linenomath}
\begin{equation*}
    \Vh{1} = \myspan \{ \varphi_{i}: \ i = 1,\dots,\Ndof\}.
\end{equation*}
\end{linenomath}

Let us denote by $\approxspace[j,E]$ a set of harmonic functions in which we approximate the VE functions, which will be characterised in Section \ref{sec:approximation_space}. Our goal is to approximate both the VE basis functions $\varphi_{j,E}$ and their gradients $\nabla \varphi_{j,E}$ with a neural network-based approximation $(\varphi^\NN_{j,E}, \nqq_{j,E})$, $\forall j = 1,\dots,\Ndof[E]$ and $\forall E \in \Th$. More specifically, the neural network aims to learn the following highly non-linear map:
\begin{equation}
    (v_j, E) \mapsto (\nvarphi_{j,E},  \nqq_{j,E})  \in \approxspace[j,E] \times \nabla \approxspace[j,E], \text{ for each vertex $v_j$ of $E$ \text{ and } $\forall E \in \Th$,}
    \label{eq:nn_map}
\end{equation} 
finding the best linear combinations of some suitable harmonic functions in $\approxspace[j,E]$ and of their gradients in $\nabla \approxspace[j,E]$ which minimize the distance between the pair $(\varphi^\NN_{j,E}, \nqq_{j,E})$ and the target $(\varphi_{j,E}, \nabla \varphi_{j,E})$ at the boundary of the element $E$, where the virtual element functions are well-known.

\subsection{The Local Neural Approximated Virtual Element Space}\label{sec:theory}

Here and in the sequel $C$ will denote a generic positive constant, with different meanings in different occurrences.

Given the approximations $\{\nvarphi_{j,E}\}_{j=1}^{\Ndof[E]}$, we define the local lowest-order NAVEM space as the set
\begin{linenomath}
\begin{equation*}
\nVh[E]{1} = \myspan\{\nvarphi_{j,E},\ j=1,\dots,\Ndof[E]\}.
\end{equation*}
\end{linenomath}

Firstly, we note that the functions $\nvarphi_{j,E}$ should belong to the VEM space $\Vh[E]{1}$  to represent a good approximation of $\varphi_{j,E}$, and, in particular, they should locally satisfy Properties $(i)$ and $(ii)$ defined in \eqref{eq:prop}. In this regard, we note that Property $(i)$ is trivially satisfied by the functions $\nvarphi_{j,E}$ by construction, since the functions in $\approxspace[j,E]$ are harmonic. Instead, Property $(ii)$ is, in general, not satisfied by functions belonging to $\approxspace[j,E]$. Nevertheless, we overcome this issue by training the neural network to learn functions $\nvarphi_{j,E}$ mimicking the VE Lagrangian basis functions $\varphi_{j,E}$ at the boundary $\partial E$ of the element $E$, where all the virtual functions are known in a closed form. In particular, our goal is to minimize the distance between the traces of the functions $\nvarphi_{j,E}$ and $\varphi_{j,E}$ on $\partial E$, i.e.
\begin{equation}
\epsilon_{j,E} = \norm[{\sob{1/2}{\partial E}}]{ \nvarphi_{j,E} - \varphi_{j,E}},
\label{eq:approximation_error}  
\end{equation}
for all $E$ and $j=1,\dots,\Ndof[E]$, to learn the non-linear relationship \eqref{eq:nn_map}.

Thanks to the harmonicity of both the virtual element functions and of the NAVEM basis functions, we can exploit the same steps performed in \cite{TrezziZerbinati2024} to state the following proposition.

\begin{proposition}
For all $E \in \Th$ and for all $j=1,\dots,\Ndof[E]$, it holds
\begin{equation}
\norm[\sob{1}{E}]{\varphi_{j,E} - \nvarphi_{j,E}} \leq C_1 \epsilon_{j,E},\qquad \norm[\leb{\infty}{\partial E}]{\varphi_{j,E} - \nvarphi_{j,E}} \leq C_2 \epsilon_{j,E}, 
\end{equation}
where $C_1$ and $C_2$ depends on $E$ and $\partial E$. 
\end{proposition}

This proposition states that the NAVEM functions $\nvarphi_{j,E}$ could be a good approximation for the related VEM Lagrangian basis functions on the entire element $E$ in the $H^1$-norm, i.e. 
\begin{equation}
\nvarphi_{j,E} \approx \varphi_{j,E} \text{ on } E, \forall j =1,\dots,\Ndof[E], \ \forall E \in \Th. 
\end{equation}
Concerning the vector of functions $\nqq_{j,E}$ in \eqref{eq:nn_map}, we observe that we are able to compute it exactly as $\nqq_{j,E} = \nabla \nvarphi_{j,E}$. However, as we will describe in Section \ref{sec:real_architecture}, sometimes better results are obtained approximating $\nabla \nvarphi_{j,E}$ independently from $\nvarphi_{j,E}$.

\subsection{The NAVEM discretization and the Online Phase}

At this point, we observe that, since the approximation of each virtual element basis function $\varphi_i$, with $i=1,\dots,\Ndof$, is computed locally, the corresponding global approximate function $\nvarphi_i$ is element-wise defined as
\begin{linenomath}
\begin{equation*}
    \nvarphi_i = 
\begin{cases}
    \nvarphi_{j,E} & \text{if $v_i$ is the $j$-the vertex of $E$ and $E \in \mathbb{S}_i$},  \\
    0 & \text{otherwise},
\end{cases}
\end{equation*}
\end{linenomath}
while the global lowest-order neural approximate virtual element space reads as
\begin{equation}
\nVh{1} = \myspan\{    \nvarphi_i:\ i =1,\dots,\Ndof\}.
\end{equation}
Thus, the NAVEM basis functions are not continuous functions across elements and they may have jumps at element interfaces. Nonetheless, we highlight that the degrees of freedom are not decoupled in our framework. Since NAVEM functions are no longer continuous across elements, we need to consider a broken version of the continuous bilinear form $\dbilin{}{}$. 
Therefore, the NAVEM discretization of problem \eqref{eq:cont_bilinear_form} reads as: \textit{Find $u_h^{\NN} \in \nVh{1}$ such that:}
\begin{equation}
\dbilinhn{u_h^{\NN}}{v_h^{\NN}} = \sum_{E \in \Th} \dbilin[E]{u_h^{\NN}}{v_h^{\NN}} = \sum_{E \in \Th} \scal[E]{f}{v_h^{\NN}} \quad \forall v_h^{\NN} \in \nVh{1}.
\label{eq:navem_var_discretization}
\end{equation} 

Firstly, we note that also in the standard Virtual Element Method we must consider a broken version of the global bilinear form due to the local definition of both the projection and stability operators.

Secondly, we observe that the lack of continuity of the functions in $\nVh{1} \nsubseteq \V$ introduces a kind of consistency error in the approximation of the solution $u \in \V$ \cite{Boffi2013}. Indeed, using integration by parts, we obtain 
\begin{linenomath}
\begin{align*}
    \dbilinhn{u}{v_h^{\NN}} &= \sum_{E \in \Th} \scal[E]{-\Delta u}{v_h^{\NN}} + \sum_{E \in \Th} \scal[\partial E]{\nabla u \cdot \nn}{v_h^{\NN}}\\
    &= \sum_{E \in \Th} \scal[E]{f}{v_h^{\NN}} + \sum_{E \in \Th} \scal[\partial E]{\nabla u \cdot \nn}{v_h^{\NN}} \\
    &= \sum_{E \in \Th} \scal[E]{f}{v_h^{\NN}} + \sum_{e \in \Eh} \scal[e]{\nabla u}{\llbracket v_h^{\NN} \rrbracket_e},
\end{align*}
\end{linenomath}
where $\llbracket v_h^{\NN} \rrbracket_e = v^{\NN}_{h|E_1} \nn_{E_1} + v^{\NN}_{h|E_2} \nn_{E_2}$ with $E_1,\ E_2$ being the elements sharing the edge $e$.  This last term measures the extent to which the continuous solution $u$ fails to satisfy the NAVEM formulation \eqref{eq:navem_var_discretization} \cite{DeDios2014}. 

Furthermore, we observe that the bilinear form $\dbilinhn{}{}$ is still symmetric positive definite with respect to the broken norm and that $\dbilinhn{}{}$ has the trivial kernel, i.e. the constant functions, and reduces to $\dbilin{}{}$ on $\V$. Thus, using the Strang's Lemma \cite{brennerscott} and the same steps developed in \cite{TrezziZerbinati2024}, we can deduce the following error bound:
\begin{equation}\label{eq:navem_error_bound}
    \norm[\NN]{u - u^{\NN}_h} \leq C \left(h + h^{-2} \varepsilon \right)
\end{equation}
where $\norm[\NN]{} = \sqrt{\dbilinhn{}{}}$ and $\varepsilon = \displaystyle \max_{j =1,\dots,\Ndof} \max_{E \in \Th} \epsilon_{j,E}$.

Finally, we highlight that, in the assembling phase, our method fully reduces to a standard finite element method since we limit to 
\begin{linenomath}
\begin{equation*}
    \mathrm{ENCODE} \longrightarrow \mathrm{PREDICT} \longrightarrow \mathrm{COMPUTE}
\end{equation*}
\end{linenomath}
that is, for each element $E \in \Th$ and for each of $j = 1,\dots,\Ndof[E]$, we
\begin{enumerate}
\item encode the information $(v_j,E)$ to generate the input of the neural network.
\item predict the coefficients of the corresponding NAVEM basis function with respect to functions contained in the related approximation space $\approxspace[j, E]$.
\item compute the integrals involved in the discretization of the problem.
\end{enumerate}

We note the $\mathrm{ENCODE}$ and $\mathrm{PREDICT}$ phases correspond to the online phase of our neural network.

\section{The Neural Network} \label{sec:neural_network}
In this section, we focus on the role of the involved neural networks. In particular, we describe the encoding of the input data $(v_j,E)$ in Section \ref{sec:input_encoding}, the approximation space $\approxspace[j,E]$ in Section \ref{sec:approximation_space} and the architectures of the neural networks and the related training procedures in Sections \ref{sec:architecture} and \ref{sec:real_architecture}.

Given the encoding of the pair $(v_j,E)$, which represents the input for our neural network \eqref{eq:nn_map}, the corresponding output is represented by the set of coefficients which express $\nvarphi_{j,E}$ and $\nabla \nvarphi_{j,E}$ with respect to the basis functions of $\approxspace[j,E]$ and $\nabla \approxspace[j,E]$. Since a neural network assumes that its input and output have constant dimensions, changing these cardinalities implies using a different neural network. In particular, we subdivide the polygons into different classes such that, in each class, all the polygons can be encoded into vectors of the same size and the corresponding basis functions can be approximated by exploiting the same number of harmonic functions.

\subsection{Input Encoding and Data Prediction}\label{sec:input_encoding}

As mentioned before, in order to predict the VEM basis functions, the first step is the encoding of the pair $(v_j, E)$ in a vector $\xx_0^{\rm{coef}}$ of a given dimension $N_0^{\rm{coeff}}$. The encoding of the input is performed in three consecutive steps: \textit{Polygon Classification}, \textit{Variability Reduction} and \textit{Input Reduction}. The first one is a mandatory step, whereas the other two steps are useful to enhance the performance of the neural network but they could be omitted.

The Polygon Classification step consists in subdividing the polygons into different classes. Pairs $(v_j,E)$ related to polygons $E$ belonging to different classes are encoded into vectors $\xx_0^{\rm{coef}}$ of different dimensions. The general rule for the classification is that, if two polygons have different numbers of vertices $\Nv$, then they belong to different classes. This trivial classification is dictated by our encoding of the pair $(v_j, E)$ into the vector $\xx_0^{\rm{coef}}$ whose dimension depends only on the number of the vertices $\Nv$ of the polygon. 

Since the input dimension is fixed for a given neural network, we need to train a different neural network for each class of polygons and thus for each value of $\Nv \geq 4$. We observe that, for $\Nv = 3$, the lowest-order virtual element method coincides with the finite element method, eliminating the need for a neural network to access point-wise evaluation of the virtual basis functions.

In our framework, the only exception to the general rule of classification is represented by the case of triangles with hanging nodes, where we devise a different strategy in order to improve the method accuracy given their importance in real-life applications \cite{Benedetto2016, Fassino2023}. We refer to Section \ref{sec:triangle_with_hanging} for the encoding of pairs $(v_j, E)$ which are related to triangles with hanging nodes. However, we highlight that, as in the virtual element framework, a triangle with one hanging node could be classified as a quadrilateral, a triangle with two hanging nodes as well as a quadrilateral with one hanging node can be classified as a pentagon and so on. This differentiation for the triangles with hanging nodes is only made to achieve very good accuracy with a very simple neural network architecture.

The Variability Reduction step is performed to reduce the variability of the elements in the datasets and enhance the neural networks' accuracy. For this purpose, we exploit the affine isomorphism defined in \cite{Teora2024}, mapping each element $E$ in the dataset in a new polygon $\srescale{E} = F_E^{-1}(E)$ which is centred at the axes origin and has unit diameter and unit anisotropic ratio. In particular, the anisotropic ratio of an element $E$ is here defined as the ratio between the maximum and the minimum eigenvalues of the inertia tensor of $E$. We recall that this map, in the absence of aligned or quasi-aligned edges, tends to uniform the elements within the same class in terms of their main geometric features, reducing the variability of the elements seen by the network.

The main role of the Input Reduction step is to shrink the dimension $N_0^{\rm{coeff}}$ of $\xx_0^{\rm{coef}}$. For this purpose, assuming a polygon classification based on the number of vertices $\Nv$, for each vertex $v_j$ of $\srescale{E}$, we consider a second affine isomorphism $G_{j,\srescale{E}}^{-1}$ that scales and rotates the element $\srescale{E}$ into a polygon $\widetilde E_j$ such that the vertex $v_j$ is mapped into the point $(1,0)$. Since the coordinates of $v_j$ are fixed, they can be excluded from the vector $\xx_0^{\rm{coef}}$. Thus, denoting by $(\widetilde x_1^r, \widetilde x_2^r)$ the coordinates of $G_{j,\srescale{E}}^{-1}(v_{r})$, we define $\xx_0^{\rm{coef}} = \begin{bmatrix} \widetilde x_1^{j+1} & \widetilde x_2^{j+1} & \dots & \widetilde x_1^{j+\Nv-1} & \widetilde x_2^{j+\Nv-1} \end{bmatrix} \in \R^{2(\Nv-1)}$, where all indeces $j$ are intended up to module $\Nv$. We observe that this map acts as a compression since the pair $(v_j, E)$ is jointly encoded into a vector of size $N_0^{\rm{coeff}} = 2(\Nv-1)$, whereas a naive encoding would require a vector of size $2\Nv+1$, i.e.  $\Nv$ $x_1$-coordinates and $\Nv$ $x_2$-coordinates plus the information about the index $j$. 

\begin{figure}[!ht]
\centering
\begin{subfigure}{0.48\textwidth}
    \includegraphics[width=\textwidth, height=3.5cm]{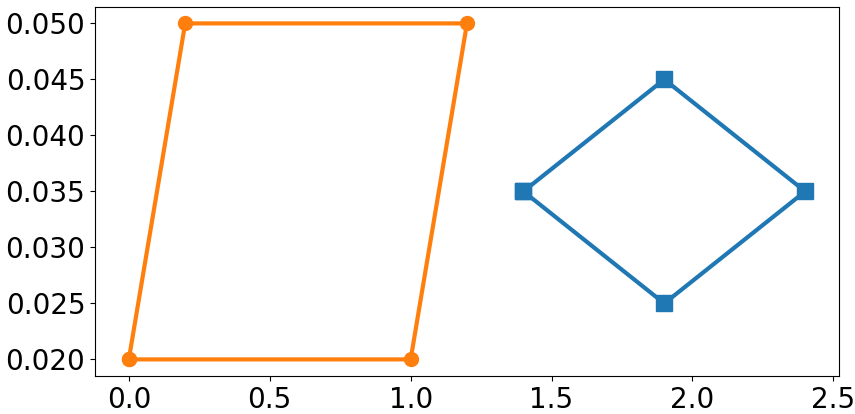}
    \caption{}
    \label{fig:original_elements}
\end{subfigure}
\begin{subfigure}{0.24\textwidth}
    \includegraphics[width=\textwidth, height=3.5cm]{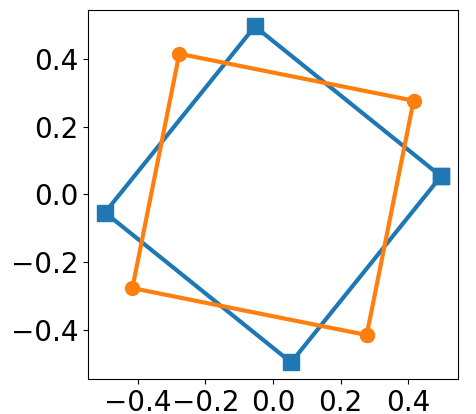}
        \caption{}
        \label{fig:variability_reduction}
\end{subfigure}
\begin{subfigure}{0.24\textwidth}
    \includegraphics[width=\textwidth, height=3.5cm]{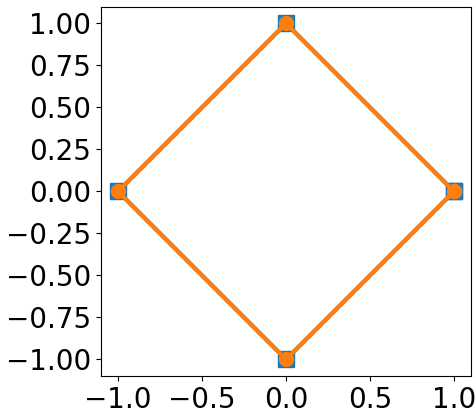}
        \caption{}
        \label{fig:input_reduction}
\end{subfigure}\quad
\caption{Input Encoding. Left: Original elements. Center: Variability Reduction. Right: Input Reduction.}
\label{fig:input_encoding}
\end{figure}

We note that this type of Input Reduction also performs a variability reduction, even though the diameters of these elements are no longer exactly $1$, but still scale as $1$. Let us explain it with a very simple example. Consider two distinct parallelograms as in Figure \ref{fig:original_elements}. The inertial mapping proposed in \cite{Teora2023} maps these two parallelograms in the same square defined up to a rotation as noted in Figure \ref{fig:variability_reduction}. Finally, the Input Reduction step fixes the rotation as highlighted in Figure \ref{fig:input_reduction} transforming the two original elements into the same element.

Actually, since we are interested in performing point-wise evaluations of the basis functions and since the input of each neural network must be a vector of fixed dimension, the final input vector $\xx_0$ for the neural network will be the concatenation of the evaluation point $\xx_0^{\rm{fun}}\in E$ and of the vector $\xx_0^{\rm{coeff}}$ which encodes the pair $(v_j, E)$. 

\subsubsection{Special case: triangles with hanging nodes}\label{sec:triangle_with_hanging}

The ability to handle meshes with hanging nodes, especially triangular meshes characterized by a copious number of hanging nodes, is very important in many contexts such as Discrete Fracture Networks \cite{Benedetto2016} or adaptive strategies \cite{Fassino2023}.
We highlight that, for such elements, the procedure described in Section \ref{sec:input_encoding} can be used but it may become very expensive or inaccurate when the number of hanging nodes grows. We thus decide to treat this case separately because of its importance in applications.

Let us consider an element $E$, which has the shape of a triangle and it is characterized by one or more hanging nodes. In order to encode the pair $(v_j, E)$, let us now consider an element $E'$ which is obtained from $E$ by removing all the hanging nodes with the only exception of $v_{j-1}$, $v_j$ and $v_{j+1}$ if these vertices are hanging nodes. As usual, $v_{j-1}=v_{\Nv-1}$ if $j=0$ and $v_{j+1}=v_0$ if $j=\Nv-1$. 
We observe that in this way we obtain en element $E'$ with at most $6$ vertices: 3 vertices which define the shape of the triangle and at most 3 hanging nodes. Moreover, let us denote by $\varphi_j'$ the VE basis function associated with the vertex $v_j$ but defined on $E'$. Since $\varphi_j$ and $\varphi_j'$ are the solution of the same Laplace problem, we can state that the removed hanging nodes do not contribute to define the shape of $\varphi_j$ and they can thus be neglected.

In the NAVEM framework, the elimination of the hanging nodes that do not influence the function $\varphi_j$ is very important to provide only useful information as the input of the neural network, limiting the input dimension and the number of possible configurations. Indeed, we observe that there exist only 6 different configurations, up to a reflection in the role of $v_{j-1}$ and $v_{j+1}$,  that is
\begin{enumerate}
    \item $v_{j-1}$ and $v_j$ are vertices of the physical triangle and $v_{j+1}$ is an hanging node;
    \item $v_{j-1}$ and $v_{j+1}$ are vertices of the physical triangle and $v_{j}$ is an hanging node;
    \item $v_{j-1}$ is a vertex of the physical triangle and $v_{j}$ and $v_{j+1}$ are hanging nodes;
    \item $v_{j-1}$, $v_j$ and $v_{j+1}$ are hanging nodes;
    \item $v_{j}$ is a vertex of the physical triangle and $v_{j-1}$ and $v_{j+1}$ are hanging nodes;
    \item $v_{j-1}$, $v_j$ and $v_{j+1}$ are vertices of the physical triangle.
\end{enumerate}
Such configurations are summarized in Table \ref{tab:triangles_configurations}. Furthermore, in order to perform variability and input reduction, we map the triangle $E'$ into the equilateral triangle of vertices $\left\{(-1,0), \left(0.5, -\frac{\sqrt{3}}{2}\right), \left(0.5, \frac{\sqrt{3}}{2}\right)\right\}$. For configurations 1, 2, 3 and 4 the triangle $E'$ is rotated such that all the hanging nodes are on the vertical edge, whereas in configuration 5 the vertex $v_j$ is in $\left(0.5, -\frac{\sqrt{3}}{2}\right)$. 

\begin{table}[!ht]
\centering
\begin{tabular}{c|c|c|c|}
\textbf{Configuration} & \textbf{$v_{j-1}$} & \textbf{$v_j$} & \textbf{$v_{j+1}$} \\ \hline
1                      & V                  & V              & H                  \\ \hline
2                      & V                  & H              & V                  \\ \hline
3                      & V                  & H              & H                  \\ \hline
4                      & H                  & H              & H                  \\ \hline
5                      & H                  & V              & H                  \\ \hline
6                      & V                  & V              & V                  \\ \hline
\end{tabular}
\caption{Existing configurations for triangles with hanging nodes. The letters V and H denote an actual vertex of the underlying physical triangle and an hanging nodes, respectively.}
\label{tab:triangles_configurations}
\end{table}
\begin{figure}[!ht]
\centering
\begin{subfigure}{0.25\textwidth}
    \includegraphics[width=\textwidth]{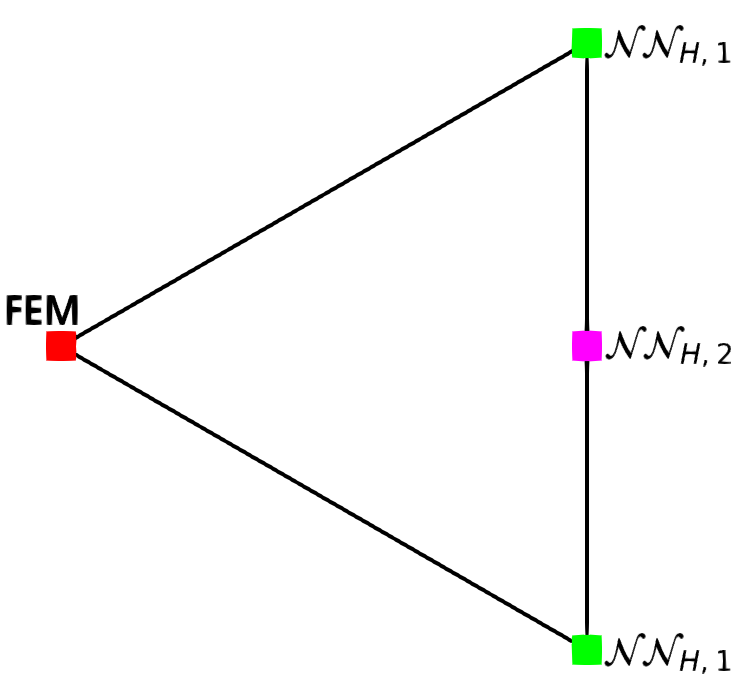}
    \caption{}
    \label{fig:conf1_triangle}
\end{subfigure}
\begin{subfigure}{0.25\textwidth}
    \includegraphics[width=\textwidth]{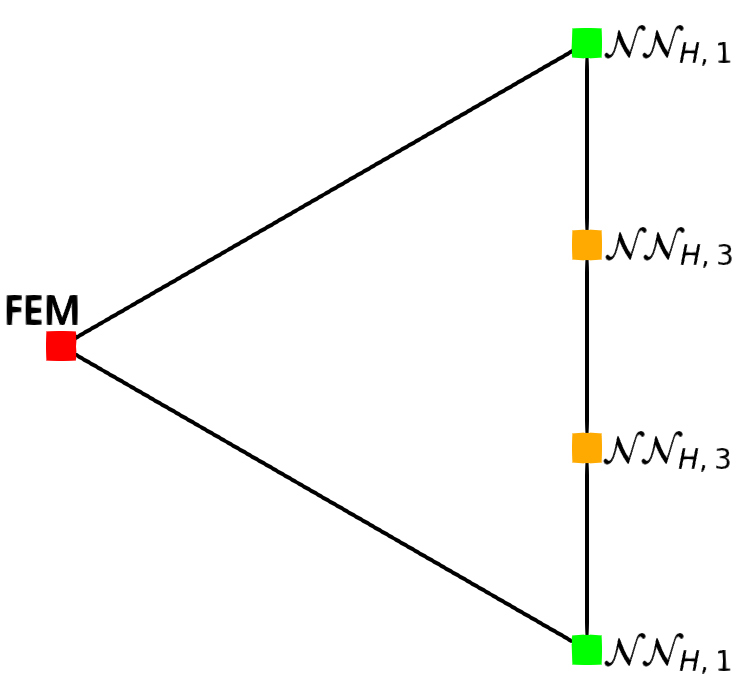}
        \caption{}
        \label{fig:conf2_triangle}
\end{subfigure}\quad
\begin{subfigure}{0.25\textwidth}
    \includegraphics[width=\textwidth]{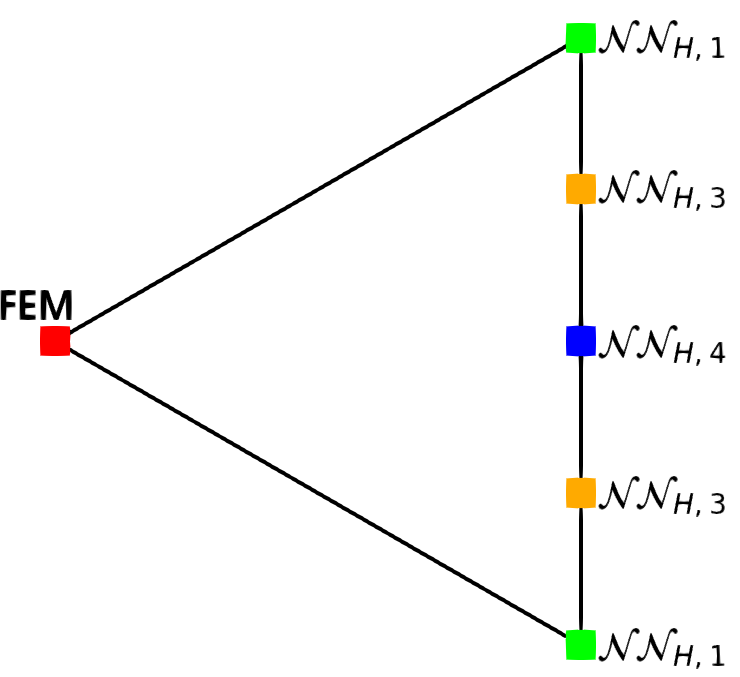}
        \caption{}
        \label{fig:conf3_triangle}
\end{subfigure}\\
\begin{subfigure}{0.25\textwidth}
    \includegraphics[width=\textwidth]{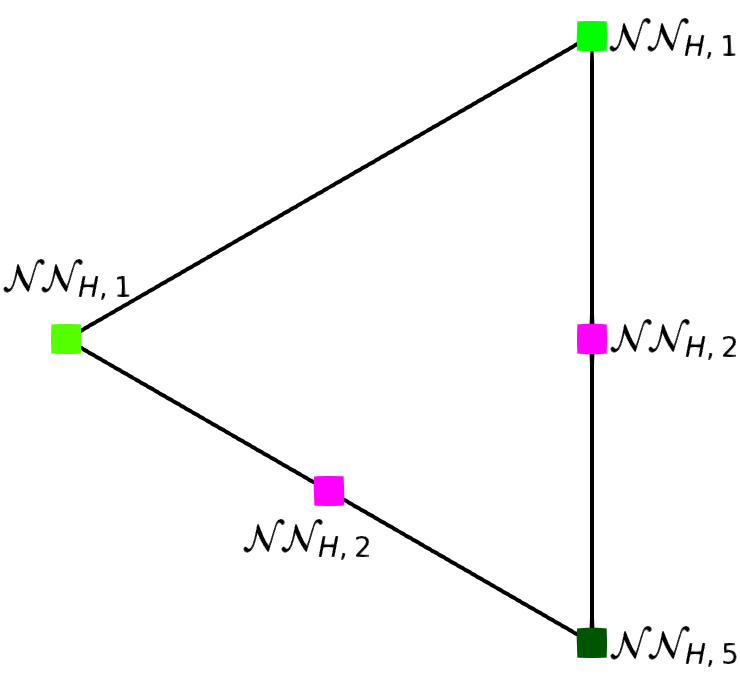}
        \caption{}
        \label{fig:conf4_triangle}
\end{subfigure}\quad
\begin{subfigure}{0.25\textwidth}
    \includegraphics[width=\textwidth]{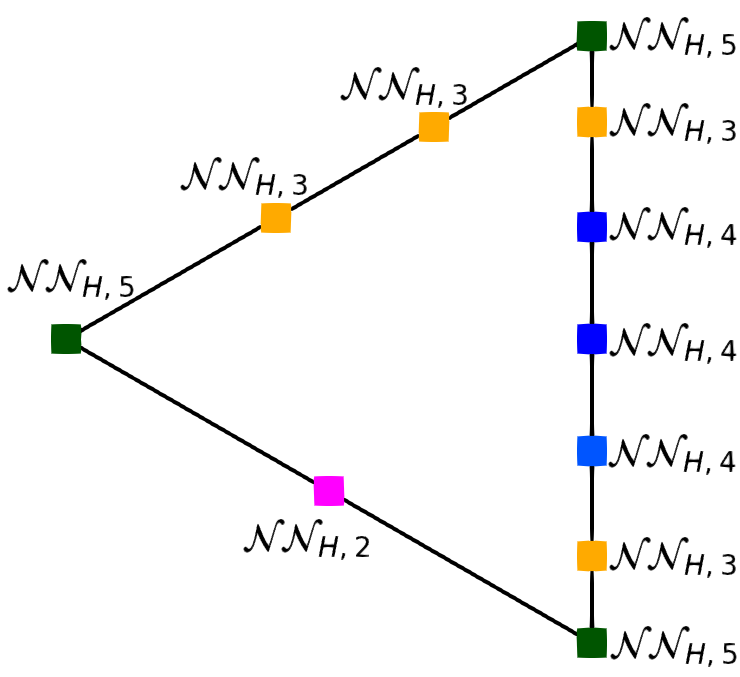}
        \caption{}
        \label{fig:general_case_conf}
\end{subfigure}
\caption{Examples of triangles with hanging nodes. The colors and labels are associated with the different configurations. The required rotations and reflections are already taken into account.}
\label{fig:triangle_configurations}
\end{figure}

To be as accurate as possible, for each configuration $i$, we train a neural network $\NN_{H,i}$, $i=1,\dots,5$, whereas for the 6-th configuration, without loss of generality, we use the known finite element basis functions, i.e. we set
\begin{linenomath}
\begin{gather*}
\nvarphi_{0,E}(x_1,x_2) = \frac{1}{3}\left(1 - 2 x_1\right),\quad \nvarphi_{1,E}(x_1,x_2) = \frac{1}{3}(x_1 - \sqrt{3}x_2 + 1), \\ \nvarphi_{2,E}(x_1,x_2) = \frac{1}{3}(x_1 + \sqrt{3}x_2 + 1),
\end{gather*}
\end{linenomath}
for the vertex $(-1,0)$, $\left(0.5, -\frac{\sqrt{3}}{2}\right)$ and $\left(0.5, \frac{\sqrt{3}}{2}\right)$ of the equilateral triangle, respectively. 

Furthermore, since the underlying physical triangle is always the same, we can avoid including the three vertices which define the shape of the triangle in the input for the network, leaving as the only inputs the curvilinear coordinates of the hanging nodes. We further note that the information about the index of the basis function is not included in the input since it is implicitly considered in the neural network configuration. Thus, the dimension of $\xx_0^{\rm{coeff}}$ of the neural network $\NN_{H,i}$ is $N_0^{\rm{coeff}}=1$ when $i=1,2$, $N_0^{\rm{coeff}}=2$ when $i=3,5$ and $N_0^{\rm{coeff}}=3$ when $i=4$. We observe that these values are very small compared to $2(\Nv[E] - 1)$, which represents the input dimension when using the general encoding procedure shown in the previous section. It is important to highlight that $\Nv[E]$ counts all the hanging nodes in $E$, which could be significantly much more than the ones in $E'$.

The possible configurations of elements $E'$, up to rotation or reflection, and a more complex configuration are shown in Figure \ref{fig:triangle_configurations}, where we label each vertex with the different neural network that should be used to predict $\nvarphi_{j,E'}$.

\subsection{The Approximation Spaces \texorpdfstring{$\approxspace[j,E]$}{H{NN}{j,E}}}\label{sec:approximation_space}

Let us introduce a reference squared region $\refregion = [-\refedge, \refedge]^2 \subset \R^2$, centred at the axes origin and with an edge length of $2\refedge$. We define the harmonic polynomial space $\HPoly{\refell}{\refregion}$, consisting of all the harmonic polynomials of degree up to $\refell \geq 0$, as the span of the following scaled harmonic polynomial basis, i.e.
\begin{equation}
      \Big\{ 1,\ \Re\left(\left(\frac{z}{\refedge}\right)^\ell\right), \Im\left(\left(\frac{z}{\refedge}\right)^\ell\right),\ \ell = 1,\dots, \refell \Big\}
      \label{eq:scaled_harmonic_polynomial}
\end{equation}
where, for simplicity, we use the complex notation $z = x_1 + i x_2$ for each point $\xx = \begin{bmatrix}
    x_1 & x_2
\end{bmatrix}^T \in \R^2$.
We observe that the dimension of $\HPoly{\refell}{\refregion}$ is $\dim \HPoly{\refell}{\refregion} = 2 \refell + 1$. Furthermore, the harmonic scaled polynomials \eqref{eq:scaled_harmonic_polynomial} and their gradients could be easily retrieved thanks to the recursive strategy presented in \cite{Perot2021}.
We then construct an orthonormal polynomial basis $\{\tilde{p}_{\beta}\}_{\beta=1}^{2\refell + 1}$ for $\HPoly{\refell}{\refregion}$ by orthogonalizing the scaled polynomial basis \eqref{eq:scaled_harmonic_polynomial} using the modified Gram-Schmidt algorithm applied twice to the Vandermonde matrix, whose columns contain the evaluations of the scaled polynomials at points forming a lattice built over $\refregion$.

\begin{figure}[!ht]
\centering
\begin{subfigure}{0.32\textwidth}
    \includegraphics[width=\textwidth]{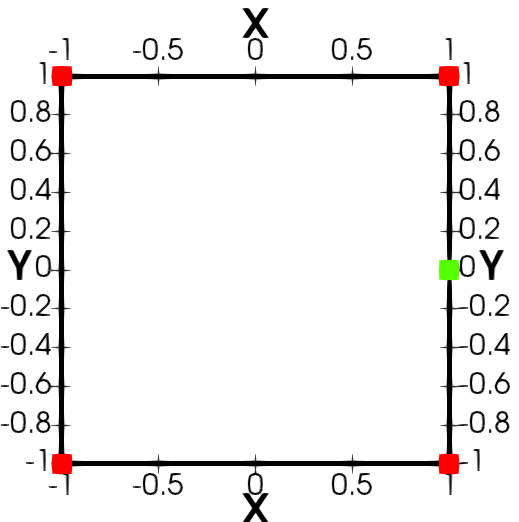}
    \caption{}
    \label{fig:hf_domain}
\end{subfigure}\quad
\begin{subfigure}{0.45\textwidth}
    \includegraphics[width=\textwidth]{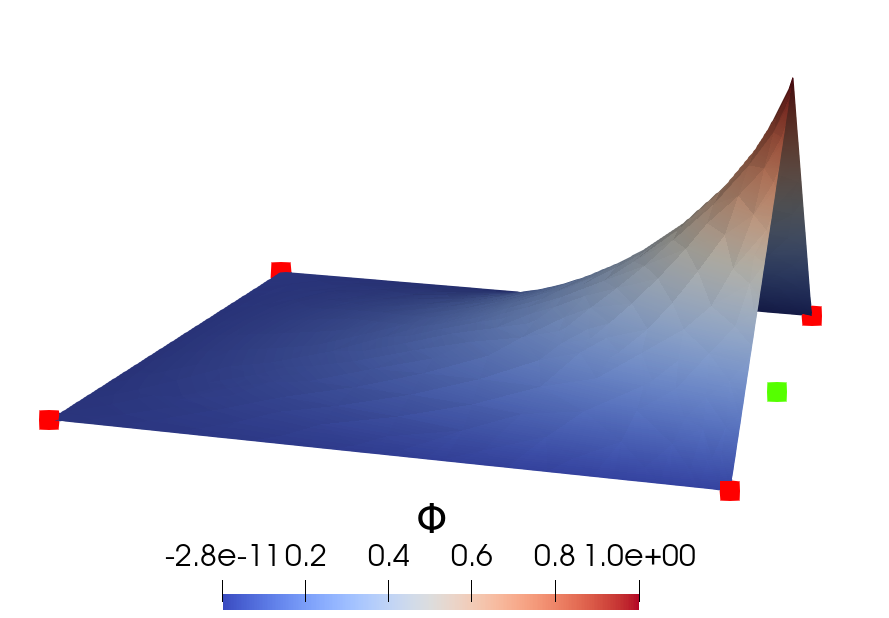}
    \caption{}
    \label{fig:hf_interp}
\end{subfigure}
\caption{Left: Domain $\Omega_\Phi$ of the function $\Phi$. Right: Shape of the function $\Phi$.}
\label{fig:hf}
\end{figure}

Next, we introduce a suitable harmonic function $\Phi$ which represents a least squares approximation of the solution $ \tilde{\Phi}$ to the following Laplace problem
\begin{equation*}
\begin{cases}
        \Delta \tilde{\Phi} = 0 & \text{in } \Omega_{\Phi} = (-1,1)^2,\\
        \tilde{\Phi} = 1 + x_2 & \text{on } \Gamma_{\Phi, 1} = \{x_1 = 1 \text{ and } -1 \leq x_2 \leq 0\},\\
        \tilde{\Phi} = 1 - x_2 & \text{on } \Gamma_{\Phi, 2} = \{x_1 = 1 \text{ and } 0 \leq x_2 \leq 1\},\\
        \tilde{\Phi} = 0 & \text{on } \partial \Omega_{\Phi} \setminus \{\Gamma_{\Phi, 1} \cup \Gamma_{\Phi, 2}\},
\end{cases}
\end{equation*}
which is computed exploiting a simplified version of the method presented in \cite{Gopal2019}.
More precisely, we determine the set of coefficients $\Big\{\{c^{\text{1}}_{\alpha}\}_{\alpha=1}^{N^{\text{1}}},\ \{c^{\text{2}}_{\beta}\}_{\beta=1}^{N^{\text{2}}}\Big\}$ of the following linear combination of harmonic functions
\begin{equation}
    \Phi(z) = \sum_{\alpha =1}^{N^{\text{1}}} c^{\text{1}}_{\alpha} \Re\left(\frac{d_{\alpha}}{z - z_{\alpha}}\right) + \sum_{\beta =0}^{N^{\text{2}}} c^{\text{2}}_{\beta} \Re\left(\left(\frac{z}{2}\right)^{\beta}\right),
    \label{eq:hanging_function}
\end{equation}
which minimizes the distance between $\Phi$ and $\tilde{\Phi}$ at the boundary of the domain $\Omega_{\Phi}$, shown in Figure \ref{fig:hf_domain}.
In Equation \eqref{eq:hanging_function}, the points $z_{\alpha} = 1 + 2 \exp\left(-4 (\sqrt{N_1} - \sqrt{\alpha})\right)$, $\alpha=1,\dots,N^{\text{1}}$ represent $N^{\text{1}}$ poles exponentially distributed along the unit exterior angle bisector $\bm{b}_{\text{green}}$ at the green vertex $z_{\text{green}} = 1 + i0$ of the domain $\Omega_{\Phi}$, whereas $d_{\alpha} = | z_{\text{green}} - z_{\alpha}| = 2 \exp\left(-4 (\sqrt{N_1} - \sqrt{\alpha}\right)$, for each $\alpha=1,\dots,N^{\text{1}}$. Figure \ref{fig:hf_interp} illustrates the shape of the function $\Phi$ obtained choosing $N^{\text{1}} = 50$ and $N^{\text{2}} = 25$.

\begin{figure}[!ht]
\centering
\includegraphics[width=0.4\textwidth]{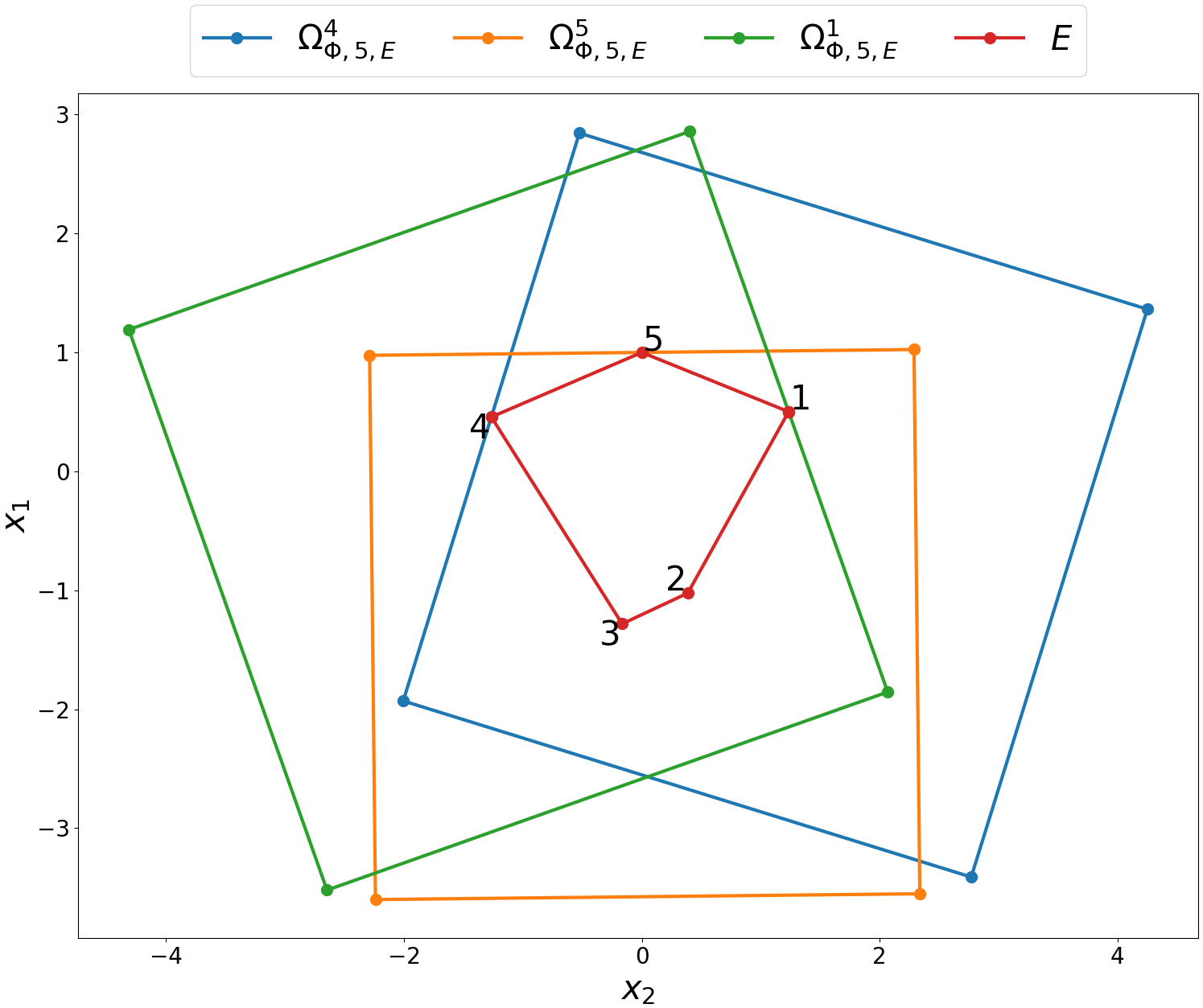}
\caption{Domain $\Omega_{\Phi,j,E}^{i}$ of the auxiliary function $\Phi^{i}_{j,E}$, for each $i=j-1,j,j+1$.}
\label{fig:phi_function}
\end{figure}

Finally, for each polygon $E$ and for each function $\varphi_{j,E}$, with $j=1,\dots,\Ndof[E]$, we introduce the following approximation space:
\begin{equation}
    \approxspace[j,E] = \myspan \Big\{ \{\tilde{p}_{\beta}\}_{\beta=1}^{2\refell + 1},\ \Phi^{j-1}_{j,E}, \Phi^{j}_{j,E}, \Phi^{j+1}_{j,E} \Big\},
    \label{eq:approxspace}
\end{equation}
where the three auxiliary functions $\Phi^{j-1}_{j,E},\ \Phi^{j}_{j,E}$ and $\Phi^{j+1}_{j,E}$ are suitable mappings of the function $\Phi$ on new domains $\Omega_{\Phi,j,E}^{j-1},\ \Omega_{\Phi,j,E}^{j}$ and $\Omega_{\Phi,j,E}^{j+1}$, respectively. These new domains are obtained through an affine isomorphism which maps $\Omega_{\Phi}$ into three corresponding squared regions $\Omega_{\Phi,j,E}^i$, for each $i=j-1,j,j+1$, defined such that $E \subset \Omega_{\Phi,j,E}^{i}$ and $z_{\text{green}}$ is mapped into the $i$-th vertex of $E$ by aligning the exterior angle bisector $\bb_{\text{green}}$ to the exterior angle bisector at the $i$-th vertex of $E$. The whole procedure is outlined in Figure \ref{fig:phi_function}. 
Although the introduction of these functions $\Phi^{j}_{j,E}$, for $i=j-1,j,j+1$, could be tricky, it helps us to capture the singularities of the function $\varphi_{j,E}$ near the vertices of the polygon $E$ conducting the same task of the functions $\Re\left(\frac{d_{\alpha}}{z - z_{\alpha}}\right)$, with $\alpha=1,\dots,N^1$ in \cite{Gopal2019}. In the latter, these functions are introduced directly in the approximation space for all the vertices of $E$ to approximate the single function $\varphi_{j,E}$ causing the number of the coefficients to predict increasing dramatically. 

We recall that, for each class of polygons, the number of outputs and thus the cardinality of the approximation space must be fixed. Nonetheless, the approximation space used to predict a basis function may differ from the approximation space used to predict another basis function, even if these basis functions are related to elements belonging to the same class. In particular, all the approximation spaces related to the same class of polygons share the same set of harmonic polynomials $\{\tilde{p}_{\beta}\}_{\beta=1}^{2\refell + 1}$, while the three auxiliary functions $\Phi^{j-1}_{j,E},\ \Phi^{j}_{j,E}$ and $\Phi^{j+1}_{j,E}$ depend on the single pair $(v_j,E)$. 

Since we fix the harmonic polynomial basis for a given neural network, $\ell^{\NN}$ and $\refedge$ are constant across polygons belonging to the same class. In particular, 
$\refedge$ is chosen in such a way these polynomials are well-scaled for each polygon $E$ predicted with the same neural network.

We further highlight that this kind of construction allows us to orthogonalize the harmonic polynomials and compute the auxiliary function $\Phi$ just once. 

\begin{remark}
We observe that, even if we are formally using rational functions as well as polynomials, we compute the integrals in Equation \eqref{eq:navem_var_discretization} employing the standard Gauss quadrature formula used in the lowest-order virtual element framework, which is exact for integrating polynomial functions of degree up to $2$. We further note that we are committing a second variational crime using this formula, but we decide to use it since we find out that its employment does not limit the convergence of the method. Other quadrature formulas could be employed to compute the integrals in Equation \eqref{eq:navem_var_discretization}, which are much more suitable for both rational functions and polynomials \cite{GAUTSCHI2001}. However, the discussion about the errors introduced by quadrature formulas is beyond the scope of this work.
\end{remark}

\subsection{The Neural Network Architecture and the Training Phase}\label{sec:architecture}

As more deeply discussed in Section \ref{sec:input_encoding}, the final input vector $\xx_0 \in \R^{2\Nv}$ can be seen as the concatenation of two vectors $\xx_0^{\rm{fun}}\in\R^2$ and $\xx_0^{\rm{coef}}\in\R^{2(\Nv-1)}$ for the general case. The former encodes the information about the point where the approximate function is evaluated, whereas the latter is associated with the polygon and the index of the VE basis function that we are approximating. Coherently with this structure, we approximate the value of a target function $\varphi$ in  $\xx_0^{\rm{fun}}$ with a linear combination of harmonic functions belonging to an approximation space $\approxspace(\xx_0^{\rm{coef}})$ whose coefficients depend only on $\xx_0^{\rm{coef}}$. We underline that we use the notation $\approxspace(\xx_0^{\rm{coef}})$ to highlight that the basis functions for $\approxspace(\xx_0^{\rm{coef}})$ actually depend on the underlying polygon (see Section \ref{sec:approximation_space}).

The architecture of the involved neural networks for a single input $\xx_0 = \begin{bmatrix}
\xx_0^{\rm{fun}} \\ \xx_0^{\rm{coeff}}
\end{bmatrix}$ is:
\begin{equation} 
  \begin{aligned}
    & \xx_0^{\varphi, \rm{coeff}} = \xx_0^{\rm{coef}}, \\
 &\xx_\ell^{\varphi,\rm{coeff}} = \rho(A^{\varphi}_\ell \xx_{\ell-1}^{\varphi, \rm{coeff}} + b^{\varphi}_\ell), \hspace{2cm} \ell = 1,...,L-1, \\
 &{\mathbf c^{\varphi}}(\xx_0^{\rm{coef}}) = A^{\varphi}_{L} \xx_{L-1}^{\varphi,\rm{coeff}}  + b^{\varphi}_L,\\[0.2cm]
 & \varphi^{\NN}(\xx_0^{\rm{fun}};  \xx_0^{\rm{coef}})  = \mathbf{V}(\xx_0^{\rm{fun}}; \xx_0^{\rm{coef}})\ {\mathbf c^{\varphi}}(\xx_0^{\rm{coef}}),\\
   & \qq^{\NN}_1(\xx_0^{\rm{fun}};  \xx_0^{\rm{coef}})  = \mathbf{V}_{x_1}(\xx_0^{\rm{fun}}; \xx_0^{\rm{coef}})\ {\mathbf c^{\varphi}}(\xx_0^{\rm{coef}}),\\
  & \qq^{\NN}_2(\xx_0^{\rm{fun}};  \xx_0^{\rm{coef}})  = \mathbf{V}_{x_2}(\xx_0^{\rm{fun}}; \xx_0^{\rm{coef}})\ {\mathbf c^{\varphi}}(\xx_0^{\rm{coef}}),
  \end{aligned}
  \label{eq:nn_formula}
\end{equation}
where
\begin{itemize}
    \item $ \varphi^{\NN}(\xx_0^{\rm{fun}};  \xx_0^{\rm{coef}})$ approximates the value $\varphi(\xx_0^{\rm{fun}})$ and 
    \begin{linenomath}
    \begin{equation*}
        \qq^{\NN}(\xx_0^{\rm{fun}};  \xx_0^{\rm{coef}}) = \begin{bmatrix}
    \qq^{\NN}_1(\xx_0^{\rm{fun}};  \xx_0^{\rm{coef}}) \\
    \qq^{\NN}_2(\xx_0^{\rm{fun}};  \xx_0^{\rm{coef}})
\end{bmatrix} = \nabla \varphi^{\NN}(\xx_0^{\rm{fun}};  \xx_0^{\rm{coef}});
    \end{equation*}
    \end{linenomath}
    \item ${\mathbf c^{\varphi}}(\cdot)$ is the trainable Multi-Layer Perceptron (MLP) (or fully-connected feed-forward neural network) \cite{popescu2009multilayer, delashmit2005recent};
    \item $L$ is the number of layers;
    \item the matrices and vectors $A^{\varphi}_\ell\in\R^{N_\ell\times N_{\ell-1}}$ and $b^{\varphi}_\ell\in\R^{N_\ell}$ store the network weights, for each $\ell=1,...,L$. Several choices are available in the literature to initialize such weights. Here, we will employ the Glorot normal initialization \cite{glorot2010understanding};
    \item $\rho:\R\rightarrow\R$ is a nonlinear activation function acting on its input vector component-wise, i.e. $\rho(\mathbf y)=\left[\rho(y_1),...,\rho(y_{n_y})\right]$ for any vector $\mathbf y\in\R^{n_y}$) \cite{sharma2020activation}. Standard activation functions are, for example, $\rho(x) = {\rm{ReLU}}(x) = \max(0,x)$, $\rho(x)=\tanh(x)$ and $\rho(x)=1/(1+e^{-x})$. In our numerical experiments, we always use the hyperbolic tangent $\rho(x)=\tanh(x)$ as the activation function.
    
    \item the \textit{Vandermonde} vectors $\mathbf{V}(\xx_0^{\rm{fun}}; \xx_0^{\rm{coef}}) \in \R^{\dim \approxspace(\xx_0^{\rm{coef}})}$, $\mathbf{V}_{x_1}(\xx_0^{\rm{fun}}; \xx_0^{\rm{coef}}) \in \R^{\dim \approxspace(\xx_0^{\rm{coef}})}$ and $\mathbf{V}_{x_2}(\xx_0^{\rm{fun}}; \xx_0^{\rm{coef}}) \in \R^{\dim \approxspace(\xx_0^{\rm{coef}})}$ are assembled such that their $i$-th element contains the evaluation of the $i$-th basis function of $\approxspace(\xx_0^{\rm{coef}})$ in $\xx_0^{\rm{fun}}$, its $x_1$-derivative and $x_2$-derivative, respectively. 
\end{itemize}

Finally, we observe that the Vandermonde vectors (or matrices when multiple evaluations are performed at the same time) do not contain trainable weights. The architecture \eqref{eq:nn_formula} can thus be seen as a standard MLP with a final multiplication for a particular matrix whose entries are non-trainable and explicitly depend on the input vector.


The weights characterizing this architecture are optimized in the following way. For each $\Nv \geq 4$, let us consider a training set $\mathcal{T}^{\mathrm{train}}$ of polygons with $\Nv$ vertices which satisfy the mesh assumptions \ref{ass:mesh}, we train the related neural network to minimize the following $L^2$-loss function
\begin{linenomath}
\begin{align}
    \nonumber \mathcal{L}_0 &= \frac{1}{\#  \mathcal{T}^{\mathrm{train}} \Nv}\sum_{E\in \mathcal{T}^{\mathrm{train}}} \sum_{j=1}^{\Nv} \norm[\leb{2}{\partial \tilde{E}_j}]{\varphi_{j,\tilde{E}_j} - \nvarphi_{j,\tilde{E}_j}}^2 \\
    &\approx  \frac{1}{\#  \mathcal{T}^{\mathrm{train}} \Nv} \sum_{E\in \mathcal{T}^{\mathrm{train}}} \sum_{j=1}^{\Nv} \sum_{q=1}^{N^q}\omega_q \left(\varphi_{j,\tilde{E}_j}(\tilde{\xx}_q) - \nvarphi_{j,\tilde{E}_j}(\tilde{\xx}_q) \right)^2
    \label{eq:l2_loss}
\end{align}
\end{linenomath}
where $\{(\omega_q, \tilde{\xx}_q)\}_{q=1}^{N^q}$ is an appropriate quadrature formula on $\partial \tilde{E}_j$ and $ \nvarphi_{j,\tilde{E}_j}(\tilde{\xx}_q) = \nvarphi(\tilde{\xx}_q; \xx_0^{\rm{coeff}})$ is the neural network output with $\xx_0^{\rm{coeff}}$ representing the encoding of $(v_j,E)$.

\subsection{An effective variant for the Training strategy}\label{sec:real_architecture}

We must observe that the inexact quadrature rule used to estimate the loss function \eqref{eq:l2_loss} and the underlying nonlinear optimization process may lead to small oscillations in the function $\nvarphi_{j,\srescale{E}}$, which may result in a poor approximation $\nqq_{j,E}$ for the gradient of the VEM basis functions $\nabla \varphi_{j,E}$.  Thus, we decide to employ a second neural network with output $\nqq_{j,E}$ to approximate the gradient of $\varphi_{j,E}$, even if in this way we introduce a new consistency error because $\nqq_{j,E}$ is, in general, different from $\nabla\nvarphi_{j,E}$.

The whole architecture of the involved neural networks for a single input $\xx_0 = \begin{bmatrix}
\xx_0^{\rm{fun}} \\ \xx_0^{\rm{coeff}}
\end{bmatrix}$ is:
\begin{linenomath}
\begin{equation*} 
  \begin{aligned}
    & \xx_0^{\varphi, \rm{coeff}} = \xx_0^{\rm{coef}}, \\
 &\xx_\ell^{\varphi,\rm{coeff}} = \rho(A^{\varphi}_\ell \xx_{\ell-1}^{\varphi, \rm{coeff}} + b^{\varphi}_\ell), \hspace{2cm} \ell = 1,...,L-1, \\
 &{\mathbf c^{\varphi}}(\xx_0^{\rm{coef}}) = A^{\varphi}_{L} \xx_{L-1}^{\varphi,\rm{coeff}}  + b^{\varphi}_L,\\[0.2cm]
 & \varphi^{\NN}(\xx_0^{\rm{fun}};  \xx_0^{\rm{coef}})  = \mathbf{V}(\xx_0^{\rm{fun}}; \xx_0^{\rm{coef}})\ {\mathbf c^{\varphi}}(\xx_0^{\rm{coef}}),\\
 & \\
     & \xx_0^{\qq, \rm{coeff}} = \xx_0^{\rm{coef}}, \\
 &\xx_\ell^{\qq,\rm{coeff}} = \rho(A^{\qq}_\ell \xx_{\ell-1}^{\qq, \rm{coeff}} + b^{\qq}_\ell), \hspace{2cm} \ell = 1,...,L-1, \\
 &{\mathbf c^{\qq}}(\xx_0^{\rm{coef}}) = A^{\qq}_{L} \xx_{L-1}^{\qq,\rm{coeff}}  + b^{\qq}_L,\\[0.2cm]
   & \qq^{\NN}_1(\xx_0^{\rm{fun}};  \xx_0^{\rm{coef}})  = \mathbf{V}_{x_1}(\xx_0^{\rm{fun}}; \xx_0^{\rm{coef}})\ {\mathbf c^{\qq}}(\xx_0^{\rm{coef}}),\\
  & \qq^{\NN}_2(\xx_0^{\rm{fun}};  \xx_0^{\rm{coef}})  = \mathbf{V}_{x_2}(\xx_0^{\rm{fun}}; \xx_0^{\rm{coef}})\ {\mathbf c^{\qq}}(\xx_0^{\rm{coef}}).
  \end{aligned}
\end{equation*}
\end{linenomath}
In this case, we initialize the weights of the first neural network $\mathbf c^{\varphi}$ as before and optimize them by minimizing the $L^2$-loss function \eqref{eq:l2_loss}, whereas the weights of the second neural network $\mathbf c^{\qq}$ are initialized with the final ones of the first network and fine-tuned to minimize the following $H^1$-loss function 
\begin{equation}
    {\cal L}_1 = \frac{1}{\#  \mathcal{T}^{\mathrm{train}} \Nv} \sum_{E\in {\mathcal{T}^{\mathrm{train}}}} \sum_{j=1}^{\Nv} \left( \hloss[j,\tilde{E}_j] \right)^2,
    \label{eq:h1_loss}
\end{equation}
where
\begin{equation}
\hloss[j,\tilde{E}_j] = \norm[{\leb{2}{\partial \tilde{E}_j}}]{ \nqq_{j,\tilde{E}_j}\cdot\tt  - \nabla \varphi_{j,\tilde{E}_j} \cdot \tt}.
\label{eq:local_loss_h1}  
\end{equation}
As in Section \ref{sec:architecture}, ${\cal L}_1$ is computed employing suitable quadrature formulas.

\begin{remark}
    We remark that, to avoid loss in the accuracy of the approximation for $\nabla \varphi_{j,E}$, a different training procedure is proposed in \cite{PintoreTeora2024}, where a single neural network is introduced which minimizes a suitable combination of $\widehat {\cal L}_0$ and $\widehat {\cal L}_1$ at the same time. Such alternative formulation can be used without additional technical complexities, but we do not focus deeply on it since we observed that it is sub-optimal in the presence of very small edges. 
\end{remark}


Note that the vector of functions $\nqq_{j,E}$ can be seen as the gradient of a function $\tilde{\varphi}^{\NN}_{j,E}$, which possibly differs from $\nvarphi_{j,E}$, defined through the minimization of the loss \eqref{eq:l2_loss}. A theoretical justification for the definition of the loss function \eqref{eq:h1_loss} is offered by the following Proposition, which shows that if two harmonic functions share the same tangential derivatives at the boundary $\partial E$, then they have the same gradient inside $E$.



\begin{lemma}
Let us assume that $u$ and $v$ are two harmonic functions with the same tangential derivatives on $\partial E$. Then,
\begin{equation}
\nabla u = \nabla v \quad \text{in } E.
\end{equation}
\end{lemma}
\begin{proof}
Since $u$ and $v$ are harmonic and share the same tangential derivatives, the function $w = u - v$ is such that 
\begin{linenomath}
\begin{equation*}
\Delta w = 0 \quad \text{in } E,\quad \nabla_{\tt} w = \bm{0}\quad \text{on } \partial E,
\end{equation*}
\end{linenomath}
where $\nabla_{\tt} w = \nabla w - (\nabla w \cdot \nn) \nn$ is the tangential derivative of $w$.
This implies that $w$ is constant on the boundary $\partial E$ and, since $w$ is harmonic, we can state that $w$ is constant everywhere on $E$. Thus, 
\begin{linenomath}
\begin{equation*}
\nabla w =  0 \quad \text{in } E \quad \Rightarrow \quad \nabla u = \nabla v \quad \text{in } E.
\end{equation*}
\end{linenomath}
\end{proof}

 We now prove that minimizing a loss function of the form \eqref{eq:h1_loss} ensures a good approximation of the gradient also inside the polygons. Note that this property is crucial since the integral forms involved in the PDEs are evaluated through quadrature rules with nodes inside the elements.

\begin{proposition}\label{prop:loss_l1}
It holds
\begin{equation}
\norm[{\leb{2}{E}}]{  \nqq_{j,E} -  \nabla \varphi_{j,E}}  \leq C \hloss[j,E],
\end{equation}
where $C$ depends on the polygon $E$.
\end{proposition}
\begin{proof}
We assume that $h_E = 1$, since the loss function \eqref{eq:local_loss_h1} is computed on mapped elements which scale as $1$. 
Let us set $\nqq_{j,E} = \nabla \tilde{\varphi}_{j,E}^{\NN}$ and $\tilde{\Psi}_{j,E} = \varphi_{j,E} - \tilde{\varphi}^{\NN}_{j,E}$. Since $\tilde{\varphi}_{j,E}^{\NN}$ is defined up to a constant, we define it in such a way $\tilde{\Psi}_{j,E}$ has zero mean value, so that the second Poincaré inequality holds \cite{CanutoNochetto2024}:
\begin{linenomath}
\begin{equation*}
	\norm[{\leb{2}{E}}]{\tilde{\Psi}_{j,E}} \leq C \norm[{\leb{2}{E}}]{\nabla \tilde{\Psi}_{j,E}}.
\end{equation*}
\end{linenomath}
From the trace theorem, we have
\begin{linenomath}
\begin{align*}
\norm[{\leb{2}{\partial E}}]{\tilde{\Psi}_{j,E}} \leq C \norm[\leb{2}{E}]{ \nabla \tilde{\Psi}_{j,E}}.
\end{align*}
\end{linenomath}
Furthermore, since $\tilde{\Psi}_{j,E}$ is harmonic, we obtain
\begin{linenomath}
\begin{align*}
\norm[{\leb{2}{E}}]{\nabla \tilde{\Psi}_{j,E}}^2 &= \int_E \nabla \tilde{\Psi}_{j,E}\cdot \nabla \tilde{\Psi}_{j,E} = \int_{\partial E} \tilde{\Psi}_{j,E} \nabla\tilde{\Psi}_{j,E} \cdot \nn \\
&\leq \norm[{\leb{2}{\partial E}}]{\tilde{\Psi}_{j,E}} \norm[{\leb{2}{\partial E}}]{\nabla \tilde{\Psi}_{j,E} \cdot \nn}  \\
&\leq C \norm[{\leb{2}{E}}]{\nabla \tilde{\Psi}_{j,E}} \norm[{\leb{2}{\partial E}}]{\nabla \tilde{\Psi}_{j,E} \cdot \nn}.
\end{align*}
\end{linenomath}
Thus,
\begin{linenomath}
\begin{align*}
\norm[{\leb{2}{E}}]{\nabla \tilde{\Psi}_{j,E}}  \leq C  \norm[{\leb{2}{\partial E}}]{\nabla \tilde{\Psi}_{j,E} \cdot \nn} .
\end{align*}
\end{linenomath}
Now, we recall that the $\leb{2}{\partial E}$-norm of the tangential derivatives is equivalent to the $\leb{2}{\partial E}$-norm of the normal derivative for harmonic functions, as a consequence of the Rellich's Identity \cite{Ammari2004}. Thus, we can conclude that
\begin{linenomath}
\begin{align*}
\norm[{\leb{2}{E}}]{  \nqq_{j,E} -  \nabla \varphi_{j,E}} &= \norm[\leb{2}{E}]{\nabla \tilde{\Psi}_{j,E}}\\
&\leq C  \norm[{\leb{2}{\partial E}}]{\tilde{\Psi}_{j,E} \cdot \nn}\\
&\leq C  \norm[{\leb{2}{\partial E}}]{\tilde{\Psi}_{j,E} \cdot \tt} = C  \norm[{\leb{2}{\partial E}}]{\nqq_{j,E} \cdot \tt -  \nabla \varphi_{j,E} \cdot \tt}\\
&= C \hloss[j,E].
\end{align*}
\end{linenomath}
\end{proof}

\section{Numerical Results}\label{sec:numerical_results}

In this section, we perform three numerical experiments to validate our procedure on different families of polygonal meshes and for different kinds of partial differential equations. In particular, after evaluating the performance of our method on a general advection-diffusion-reaction problem, we investigate its behaviour in solving anisotropic problems and non-linear problems, where the use of a stabilizing form or polynomial projectors may limit the performance of the method.

Denoting by $u$ the exact solution of the underlying problem, for each family of meshes, we test the performance of the NAVEM by looking at the behaviour of the following errors
\begin{equation}
    \mathrm{err}_0^\NN = \sqrt{\sum_{E\in\Th} \norm[0,E]{u - u^{\NN}_h}^2},\quad \mathrm{err}_1^\NN = \sqrt{\sum_{E\in\Th} \norm[0,E]{\nabla u - \nabla u^{\NN}_h}^2},
    \label{eq:nn_errors}
\end{equation}
as the family mesh parameter $h$ decreases. 

We further compare the performance of our method with the standard VE method which is available in the literature for the corresponding problem. Since we consider problems with variable coefficients, in the following numerical experiments, we will adopt the virtual element discretization introduced in \cite{LBe16} for the linear case, which is based on the definition of an \textit{enhanced} space. This alternative formulation allows to compute the $L^2$-projections of the virtual element functions on polynomial spaces of higher polynomial degrees in order to avoid loss of accuracy in the presence of variable coefficients. Thus, in the following, for each $k \geq 0$ and $E \in \Th $, we will denote by $\proj{E,0}{k}: \Vh[E]{1} \to \Poly{k}{E}$ the $L^2$-polynomial projection of virtual element functions. Furthermore, without loss of generality, we use the same symbol also to denote the $L^2$ polynomial projection of vector-valued functions. In the following, we further employ the standard \textit{dofi-dofi} stabilization term \eqref{eq:dofidofi} as a stabilizing form for the virtual element method since we deal with the lowest-order discretization.
Moreover, since we can not access to the point-wise evaluation of virtual functions, we define the VEM errors as usual, that is
\begin{equation}
    \mathrm{err}_0^{\rm{VEM}} = \sqrt{\sum_{E\in\Th} \norm[0,E]{u - \proj{E,0}{1} u_h^{\rm{VEM}}}^2},\quad \mathrm{err}_1^{\rm{VEM}} = \sqrt{\sum_{E\in\Th} \norm[0,E]{\nabla u - \proj{E,0}{0} \nabla u_h^{\rm{VEM}}}^2}.
    \label{eq:vem_errors}
\end{equation}

To train the neural networks, we use a combination of the ADAM optimizer \cite{kingma2014adam} with the BFGS optimizer \cite{wright1999numerical} to optimize the weights of the neural networks.
Moreover, to avoid problems related to overfitting, it is often advisable to add a regularization term to the loss function. We thus adopt a standard regularization technique penalizing the $L^2$-norm of the trainable coefficients and we set the regularization coefficients to $10^{-8}$.

Finally, for simplicity, we choose $\ell^{\NN} = 20$ and $\refedge = 3$ for all the neural networks. A better fine-tuning strategy of the parameter $\ell^{\NN}$ could be performed to further improve the efficiency of the method, reducing the number of function evaluations, but this is beyond the scope of the manuscript.

\subsection{Meshes and Training sets}

\begin{figure}[!ht]
\centering
\begin{subfigure}{0.31\textwidth}
    \includegraphics[width=\textwidth]{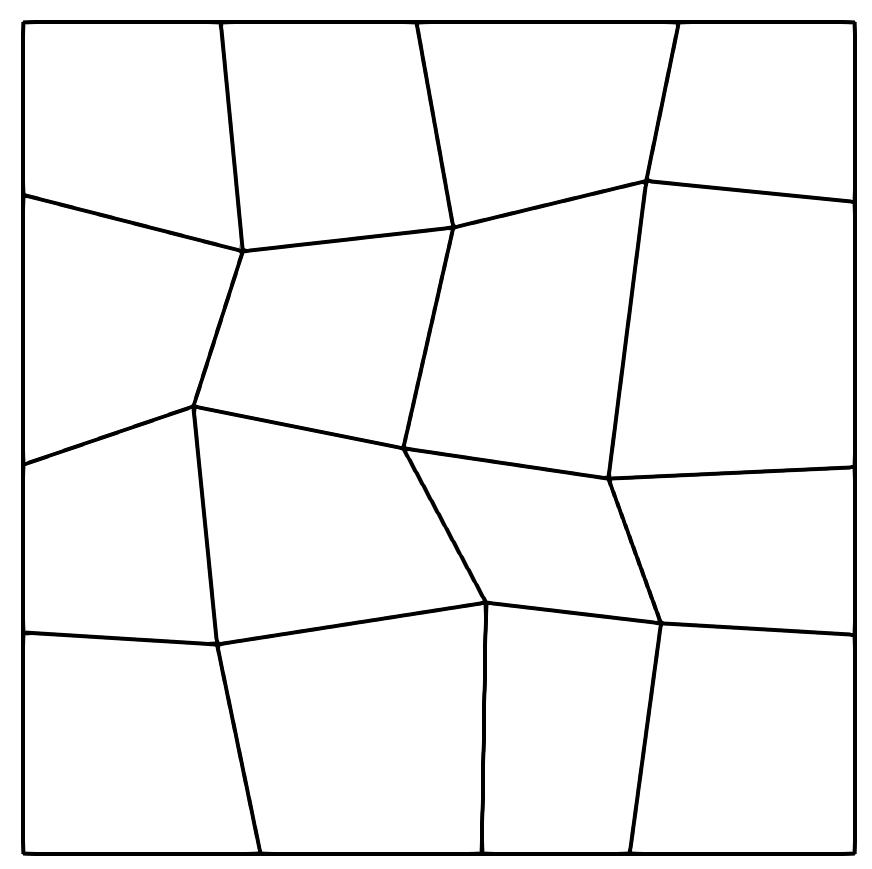}
    \caption{}
    \label{fig:rand_dist_mesh}
\end{subfigure}\quad
\begin{subfigure}{0.31\textwidth}
    \includegraphics[width=\textwidth]{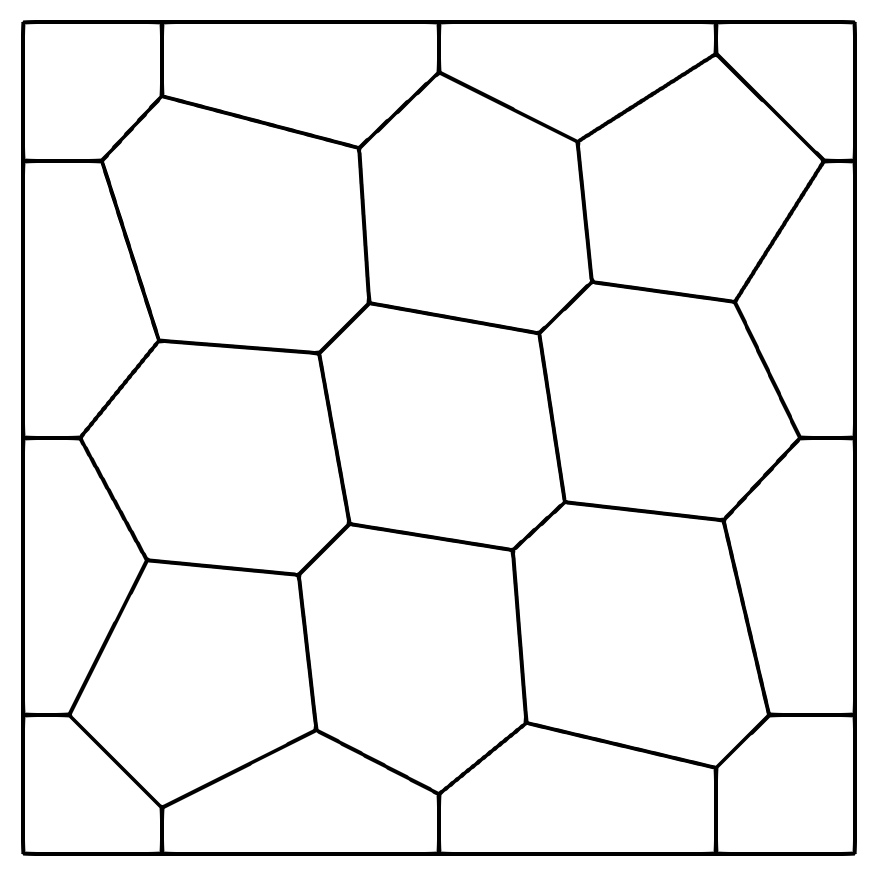}
    \caption{}
    \label{fig:voronoi_mesh}
\end{subfigure}\quad
\begin{subfigure}{0.31\textwidth}
    \includegraphics[width=\textwidth]{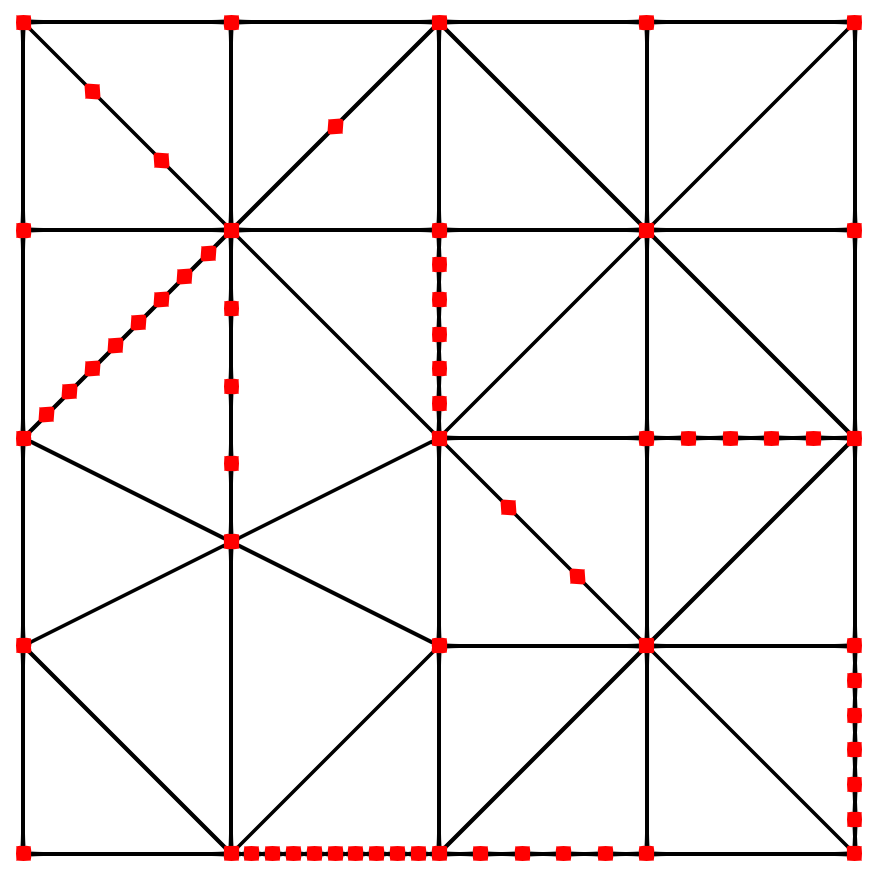}
    \caption{}
    \label{fig:triangle_with_hang_mesh}
\end{subfigure}
\caption{Left: First mesh of RDQM family. Center: First mesh of VM family. Right: First mesh of HTM family. The red dots denote the hanging nodes.}
\label{fig:meshes}
\end{figure}

In this section, we describe the three families of meshes used in the numerical experiments and the training sets that we use to train the related neural networks. Each family of meshes is made up of four meshes with decreasing mesh parameters $h$. The first mesh for each family is shown in Figure \ref{fig:meshes}. 

\begin{table}[]
\centering
\resizebox{\hsize}{!}{
\begin{tabular}{ccccccccccccccc}
                                                                & \multirow{2}{*}{\textbf{\# \{E\}}} & \multicolumn{2}{c}{\textbf{$N_v^E$}}     & \multirow{2}{*}{\textbf{}} & \multicolumn{2}{c|}{\textbf{Area}}               & \multicolumn{2}{c|}{\textbf{Diameter}}           & \multicolumn{2}{c|}{\textbf{Anisotropic Ratio}}  & \multicolumn{2}{c|}{\textbf{Edge ratio}}         & \textbf{$\mathcal{L}^0$}  & \textbf{$\mathcal{L}^1$}                       \\
                                                                &                                    & \textbf{min}       & \textbf{max}        &                            & \textbf{min} & \multicolumn{1}{c|}{\textbf{max}} & \textbf{min} & \multicolumn{1}{c|}{\textbf{max}} & \textbf{min} & \multicolumn{1}{c|}{\textbf{max}} & \textbf{min} & \multicolumn{1}{c|}{\textbf{max}} & \textbf{avg}              & \textbf{avg}                                   \\
\multirow{2}{*}{\textbf{$\mathcal{T}^{\mathrm{test}}_{RDQM}$}}  & \multirow{2}{*}{736}               & \multirow{2}{*}{4} & \multirow{2}{*}{4}  & $E$                        & 1.20e-03     & \multicolumn{1}{c|}{9.27e-02}     & 5.32e-02     & \multicolumn{1}{c|}{4.46e-01}     & 1.02e+00     & \multicolumn{1}{c|}{4.10e+00}     & 1.05e+00     & \multicolumn{1}{c|}{2.43e+00}     & \multirow{2}{*}{2.00e-03} & \multicolumn{1}{c|}{\multirow{2}{*}{4.18e-03}} \\
                                                                &                                    &                    &                     & $\srescale{E}$             & 4.69e-01     & \multicolumn{1}{c|}{5.00e-01}     & 1.00e+00     & \multicolumn{1}{c|}{1.00e+00}     & 1.00e+00     & \multicolumn{1}{c|}{1.00e+00}     & 1.01e+00     & \multicolumn{1}{c|}{2.30e+00}     &                           & \multicolumn{1}{c|}{}                          \\ \hline
\multirow{2}{*}{\textbf{$\mathcal{T}^{\mathrm{test}}_{VM}$}}    & \multirow{2}{*}{1054}              & \multirow{2}{*}{5} & \multirow{2}{*}{7}  & $E$                        & 3.09e-04     & \multicolumn{1}{c|}{7.62e-02}     & 2.84e-02     & \multicolumn{1}{c|}{3.68e-01}     & 1.00e+00     & \multicolumn{1}{c|}{6.09e+00}     & 1.00e+00     & \multicolumn{1}{c|}{9.04e+00}     & \multirow{2}{*}{6.91e-03} & \multicolumn{1}{c|}{\multirow{2}{*}{1.96e-02}} \\
                                                                &                                    &                    &                     & $\srescale{E}$             & 4.97e-01     & \multicolumn{1}{c|}{7.10e-01}     & 1.00e+00     & \multicolumn{1}{c|}{1.00e+00}     & 1.00e+00     & \multicolumn{1}{c|}{1.00e+00}     & 1.00e+00     & \multicolumn{1}{c|}{4.53e+00}     &                           & \multicolumn{1}{c|}{}                          \\ \hline
\multirow{2}{*}{\textbf{$\mathcal{T}^{\mathrm{test}}_{HTM}$}}   & \multirow{2}{*}{542}               & \multirow{2}{*}{3} & \multirow{2}{*}{20} & $E$                        & 1.90e-03     & \multicolumn{1}{c|}{4.69e-02}     & 7.11e-02     & \multicolumn{1}{c|}{3.75e-01}     & 1.00e+00     & \multicolumn{1}{c|}{6.84e+00}     & 1.00e+00     & \multicolumn{1}{c|}{1.79e+01}     & \multirow{2}{*}{1.52e-03} & \multicolumn{1}{c|}{\multirow{2}{*}{4.69e-03}} \\
                                                                &                                    &                    &                     & $\srescale{E}$             & 1.30e+00     & \multicolumn{1}{c|}{1.30e+00}     & 1.73e+00     & \multicolumn{1}{c|}{1.73e+00}     & 1.00e+00     & \multicolumn{1}{c|}{1.00e+00}     & 1.00e+00     & \multicolumn{1}{c|}{1.00e+01}     &                           & \multicolumn{1}{c|}{}                          \\ \hline
                                                                &                                    &                    &                     &                            &              &                                   &              &                                   &              &                                   &              &                                   &                           &                                                \\ \hline
\multirow{2}{*}{\textbf{$\mathcal{T}^{\mathrm{train}}_{RDQM}$}} & \multirow{2}{*}{1000}              & \multirow{2}{*}{4} & \multirow{2}{*}{4}  & $E$                        & 4.03e-02     & \multicolumn{1}{c|}{9.02e-01}     & 1.00e+00     & \multicolumn{1}{c|}{1.41e+00}     & 1.01e+00     & \multicolumn{1}{c|}{9.37e+02}     & 1.04e+00     & \multicolumn{1}{c|}{5.13e+01}     & \multirow{2}{*}{4.62e-03} & \multicolumn{1}{c|}{\multirow{2}{*}{1.00e-02}} \\
                                                                &                                    &                    &                     & $\srescale{E}$             & 4.34e-01     & \multicolumn{1}{c|}{5.00e-01}     & 1.00e+00     & \multicolumn{1}{c|}{1.00e+00}     & 1.00e+00     & \multicolumn{1}{c|}{1.00e+00}     & 1.01e+00     & \multicolumn{1}{c|}{4.24e+01}     &                           & \multicolumn{1}{c|}{}                          \\ \hline
\multirow{2}{*}{\textbf{$\mathcal{T}^{\mathrm{train}}_{VM}$}}   & \multirow{2}{*}{18701}             & \multirow{2}{*}{4} & \multirow{2}{*}{7}  & $E$                        & 1.97e-05     & \multicolumn{1}{c|}{9.02e-01}     & 7.14e-03     & \multicolumn{1}{c|}{1.41e+00}     & 1.00e+00     & \multicolumn{1}{c|}{9.37e+02}     & 1.00e+00     & \multicolumn{1}{c|}{5.13e+01}     & \multirow{2}{*}{2.44e-03} & \multicolumn{1}{c|}{\multirow{2}{*}{5.79e-03}} \\
                                                                &                                    &                    &                     & $\srescale{E}$             & 4.34e-01     & \multicolumn{1}{c|}{7.14e-01}     & 1.00e+00     & \multicolumn{1}{c|}{1.00e+00}     & 1.00e+00     & \multicolumn{1}{c|}{1.00e+00}     & 1.00e+00     & \multicolumn{1}{c|}{4.24e+01}     &                           & \multicolumn{1}{c|}{}                          \\ \hline
\multirow{2}{*}{\textbf{$\mathcal{T}^{\mathrm{train}}_{HTM}$}}  & \multirow{2}{*}{7670}              & \multirow{2}{*}{4} & \multirow{2}{*}{6}  & $E$                        & 1.30e+00     & \multicolumn{1}{c|}{1.30e+00}     & 1.73e+00     & \multicolumn{1}{c|}{1.73e+00}     & 1.00e+00     & \multicolumn{1}{c|}{1.00e+00}     & 2.02e+00     & \multicolumn{1}{c|}{5.61e+02}     & \multirow{2}{*}{4.19e-03} & \multicolumn{1}{c|}{\multirow{2}{*}{1.65e-02}} \\
                                                                &                                    &                    &                     & $\srescale{E}$             & 1.30e+00     & \multicolumn{1}{c|}{1.30e+00}     & 1.73e+00     & \multicolumn{1}{c|}{1.73e+00}     & 1.00e+00     & \multicolumn{1}{c|}{1.00e+00}     & 2.02e+00     & \multicolumn{1}{c|}{5.61e+02}     &                           & \multicolumn{1}{c|}{}                          \\ \hline
\end{tabular}}
\caption{Geometric properties of training and test meshes used throughout the experiments.}
\label{tab:geometric_data}
\end{table}

\begin{figure}[!ht]
\centering
\begin{subfigure}{0.24\textwidth}
    \includegraphics[width=\textwidth]{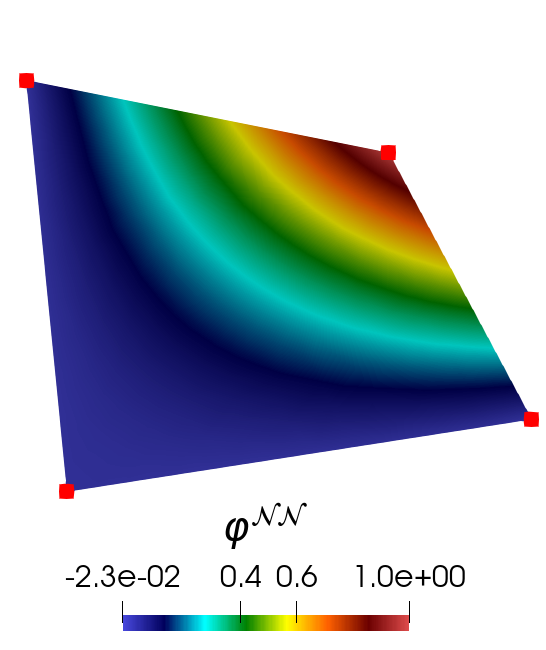}
    \caption{}
\end{subfigure}
\begin{subfigure}{0.24\textwidth}
    \includegraphics[width=\textwidth]{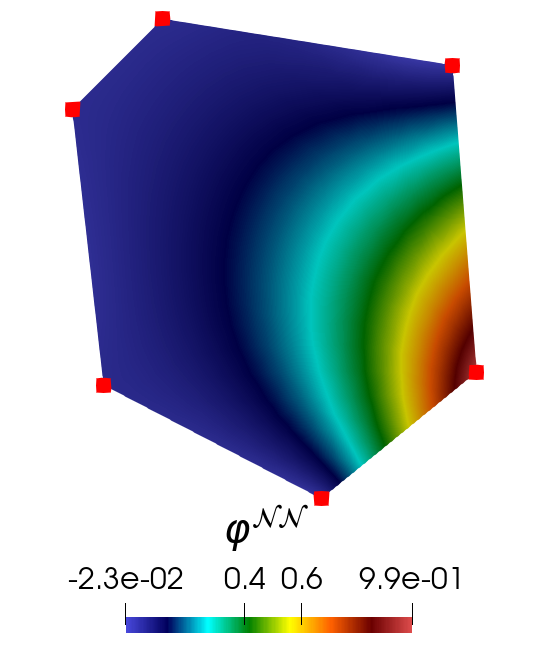}
    \caption{}
\end{subfigure}
\begin{subfigure}{0.24\textwidth}
    \includegraphics[width=\textwidth]{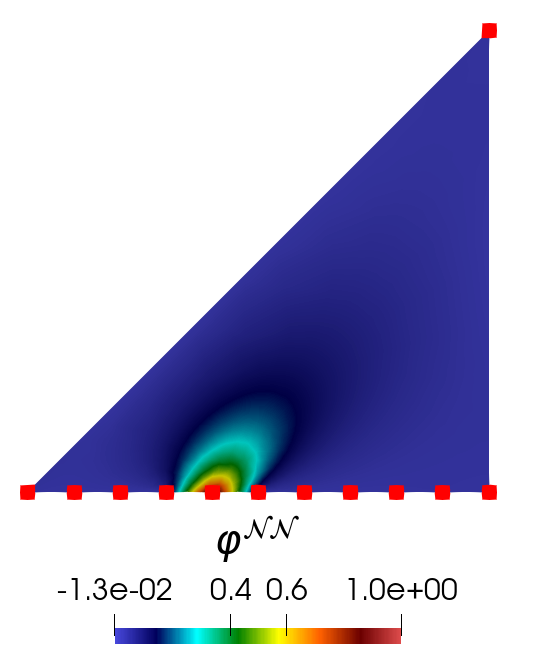}
    \caption{}
\end{subfigure}
\begin{subfigure}{0.24\textwidth}
    \includegraphics[width=\textwidth]{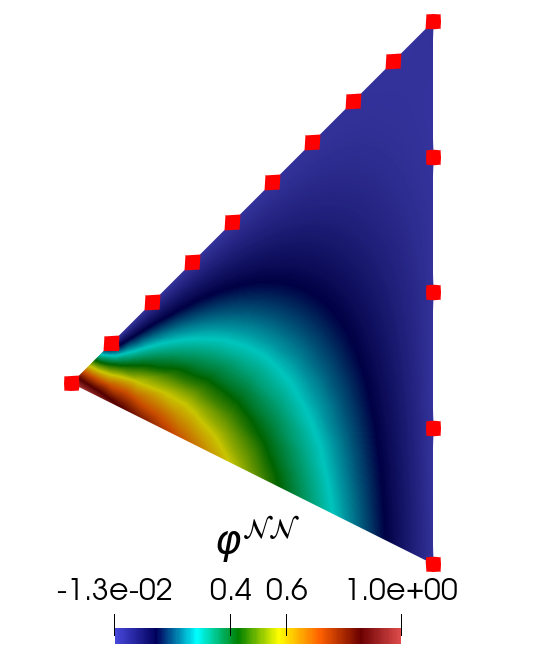}
    \caption{}
\end{subfigure}
\caption{Contour plots of some basis functions predicted by the neural network on test elements. The red dots denote the vertices of the polygons.}
\label{fig:basis_functions}
\end{figure}

\paragraph{Random Distorted Quadrilateral Meshes (RDQM)} The first family of meshes is obtained starting from a family of Cartesian meshes made up of identical squares. The vertices of such Cartesian meshes are then randomly perturbed to generate meshes with random quadrilaterals, as the one represented in Figure  \ref{fig:rand_dist_mesh}. For this kind of mesh, a single neural network $\NN_4$ is needed since all the polygons have the same number of vertices. More precisely, two neural networks are trained to minimize the loss functions \eqref{eq:l2_loss} and \eqref{eq:h1_loss}, respectively, over a training set of 1000 convex quadrilaterals which are randomly generated through the Python Library \href{https://pypi.org/project/polygenerator/}{polygenerator}.  The chosen neural network architectures comprise 4 hidden layers with 40 neurons in each layer. 

\paragraph{Voronoi meshes (VM)} The second family of meshes is a set of Voronoi meshes. The coarsest mesh is shown in Figure \ref{fig:voronoi_mesh}. The elements included in these meshes are convex quadrilaterals, pentagons, hexagons and heptagons. For the quadrilaterals, we use the neural networks $\NN_4$ trained for the family RDQM, whereas for the other polygons, we train three networks $\NN_{V,i}$, with $i=5,6,7$, all with $5$ layers and $50$ neurons in each layer. The training sets are obtained starting from different refinements of Voronoi meshes and then splitting the elements into different training sets according to the number of their vertices. 

\paragraph{Triangular Meshes with Hanging nodes (HTM)} The last family is a set of four triangular meshes with hanging nodes. These meshes are generated starting from standard triangular meshes and randomly selecting a subset of edges to which we add a random number between $1$ and $10$ of equispaced hanging nodes. The coarsest mesh is shown in Figure \ref{fig:triangle_with_hang_mesh}, where the red dots represent the hanging nodes. As described in Section \ref{sec:triangle_with_hanging}, we train 5 neural networks to approximate the basis functions related to this kind of mesh. Thus, we consider a set of identical neural networks comprising $5$ hidden layers with $50$ neurons in each layer. 
For each neural network, we construct the corresponding training dataset by adding $1$, $2$ or $3$ hanging nodes (depending on the configuration) on the suitable edges of the reference equilateral triangle, as described in Section \ref{sec:triangle_with_hanging}. We choose $10^{-2}$ as the minimum distance allowed between two consecutive vertices.

Table \ref{tab:geometric_data} summarizes the main geometric properties of the elements belonging to the training and test meshes and reports the square root of the related losses \eqref{eq:l2_loss} and \eqref{eq:h1_loss}. In particular, we report statistics about the area, the diameter, the anisotropic ratio and the edge ratio, which is given by the ratio between the maximum and the minimum lengths of the edges of the elements. These geometric properties are reported for both the original $E$ and the mapped elements $\srescale{E}$. We recall that the mapped elements for the $RDQM$ and $VM$ are the elements obtained through the inertial mapping introduced in \cite{Teora2023}, while all the physical triangles related to the $HTM$ family are mapped to the same reference equilateral triangle. 

In particular, concerning the test and the training meshes related to the family $RDQM$ and $VM$, which are grouped in $\mathcal{T}^{\mathrm{test}}_{RDQM}$-$\mathcal{T}^{\mathrm{train}}_{RDQM}$ and $\mathcal{T}^{\mathrm{test}}_{VM}$-$\mathcal{T}^{\mathrm{train}}_{VM}$, respectively, we can note that thanks to the inertial mapping we are able to strongly reduce the variability of the elements seen by the neural network obtaining elements which approximately the same area, with unit diameter and unit anisotropic ratio. Concerning the test and the training sets for the $HTM$ family which are grouped in $\mathcal{T}^{\mathrm{test}}_{HTM}$-$\mathcal{T}^{\mathrm{train}}_{HTM}$, we can see that the values for the area, diameter and the anisotropic ratio are equal for all the mapped elements since the last corresponds to the same physical triangle. The only variability regards the value for the edge ratio, which takes into account the presence of the hanging nodes. Furthermore, we recall that the training set for the HTM family is made up of different copies of the same reference triangle, thus actually the elements $E \in \mathcal{T}^{\mathrm{train}}_{HTM}$ are equal to the related mapped element $\srescale{E}$. 
Finally, we can note that the losses related to the training set and the test set are very similar to each other and this means that the test elements are well represented by the chosen training sets. In Figure \ref{fig:basis_functions}, we report the contour plots of some basis functions predicted on test elements, which behave as expected.

\begin{figure}[!ht]
\centering
\begin{subfigure}{0.32\textwidth}
    \includegraphics[width=\textwidth]{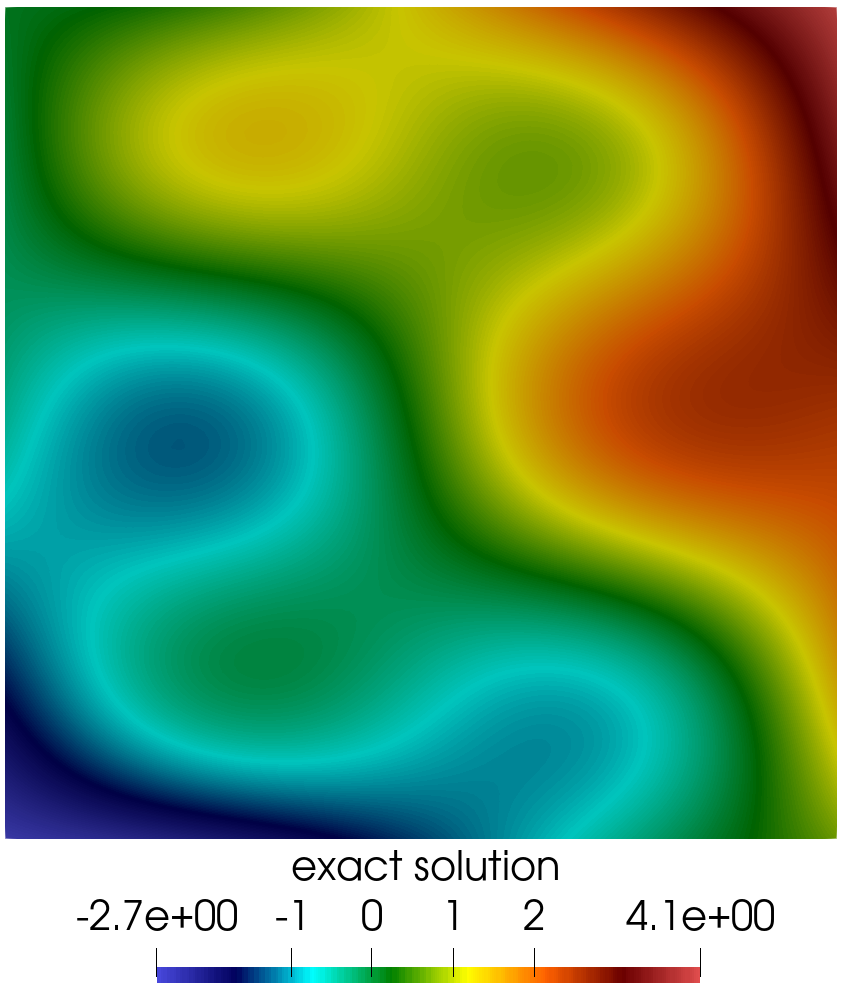}
    \caption{}
    \label{fig:test1_solutions}
\end{subfigure}\quad
\begin{subfigure}{0.32\textwidth}
    \includegraphics[width=\textwidth]{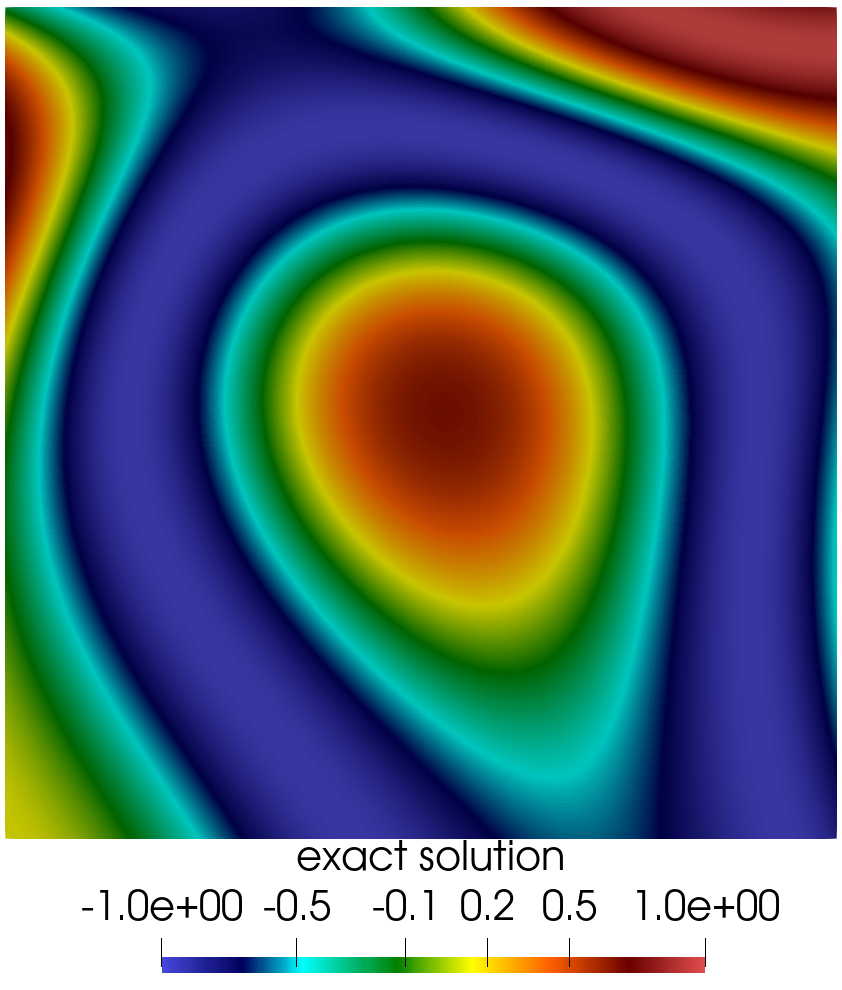}
    \caption{}
    \label{fig:test2_solutions}
\end{subfigure}
\begin{subfigure}{0.32\textwidth}
    \includegraphics[width=\textwidth]{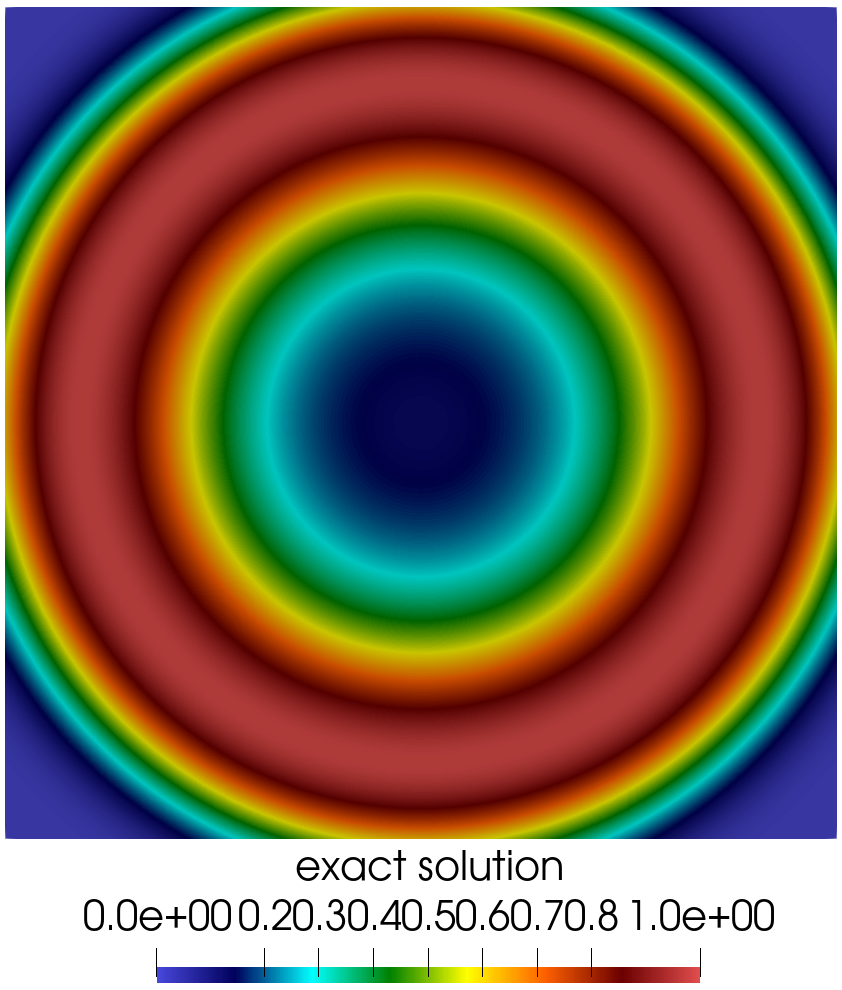}
    \caption{}
    \label{fig:test3_solutions}
\end{subfigure}
\caption{Contour plots of the exact solutions \eqref{eq:test1_exact_solution}, \eqref{eq:test2_exact_solution} and \eqref{eq:test3_exact_solution} from right to left. Test 1 (Left), Test 2 (Center) and Test 3 (Right).}
\label{fig:exact_solutions}
\end{figure}

\subsection{Test problem 1: Advection-Diffusion-Reaction problem} \label{sec:conv_rates}

\begin{figure}[t]
\centering
\begin{subfigure}{0.32\textwidth}
    \includegraphics[width=\textwidth, clip]{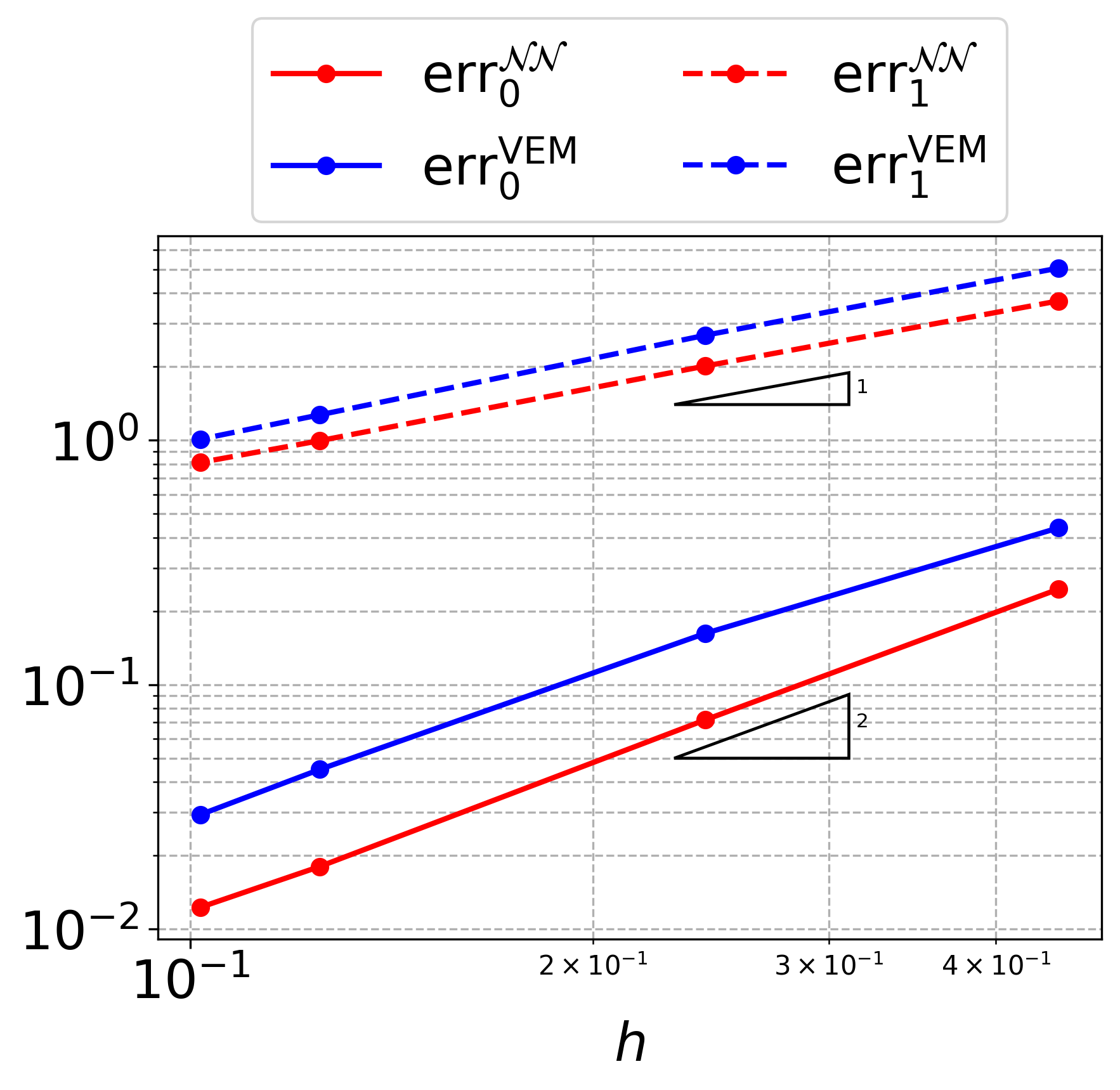}
    \caption{}
    \label{}
\end{subfigure}
\begin{subfigure}{0.32\textwidth}
    \includegraphics[width=\textwidth, clip]{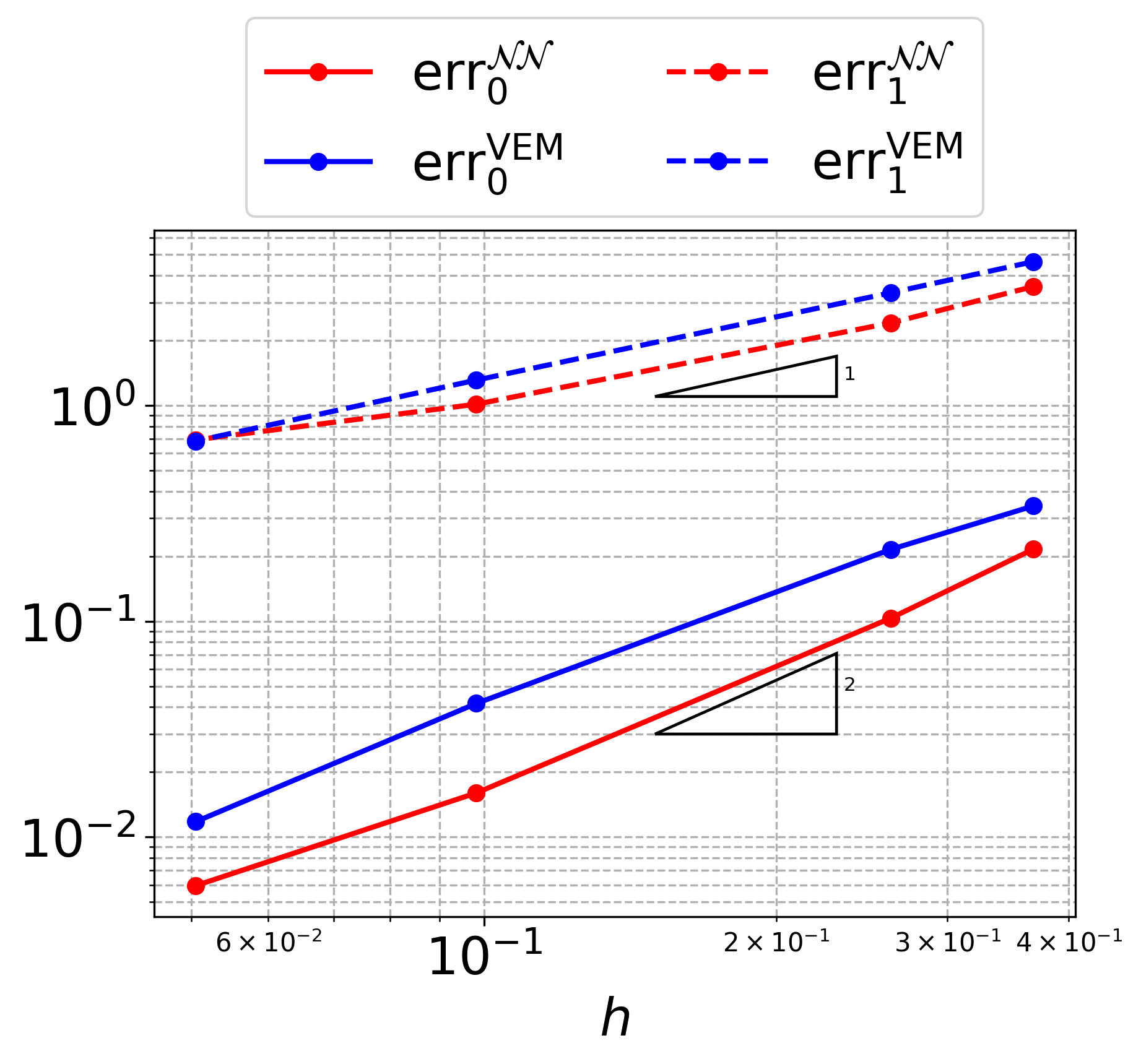}
    \caption{}
    \label{}
\end{subfigure}
\begin{subfigure}{0.32\textwidth}
    \includegraphics[width=\textwidth, clip]{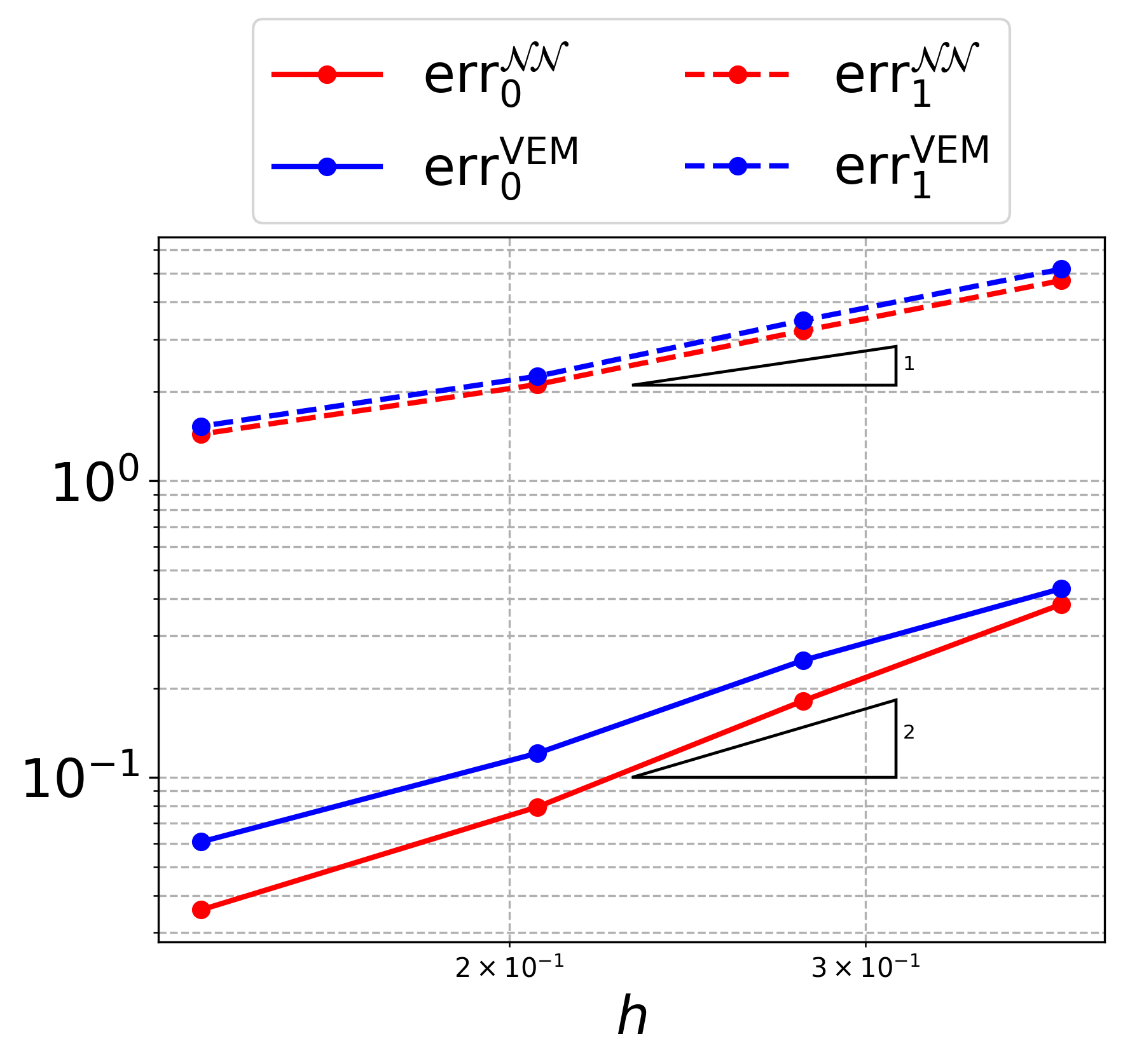}
    \caption{}
    \label{}
\end{subfigure}
\caption{Test 1: NAVEM and VEM errors w.r.t $h$. RDQM (Left), VM (Center) and HTM (Right).}
\label{fig:test1_errors}
\end{figure}

In this first experiment, we test the NAVEM method on a simple diffusion-advection-reaction problem. In particular, we consider the following boundary value problem on $\Omega = (0,1)^2$
\begin{equation}
    \begin{cases}
       - \nabla \cdot \left( \DD(\xx) \nabla u\right) + \bbeta(\xx) \cdot \nabla u + \gamma(\xx) u = f & \text{in } \Omega,\\
        u = g_D & \text{on } \Gamma,
    \end{cases}
    \label{eq:diffreacadv}
\end{equation} 
where 
\begin{linenomath}
\begin{equation*}
    \DD(\xx) = \begin{bmatrix}
        1 + x_2^2 & -x_1 x_2\\
        -x_1 x_2 & 1 + x_1^2
    \end{bmatrix}, \quad\quad \bbeta(\xx) = \begin{bmatrix}
     x_1 \\
     -x_2
    \end{bmatrix}, \quad\quad \gamma(\xx) = x_1 x_2,
\end{equation*}
\end{linenomath}
while the Dirichlet boundary condition $g_D$ and the forcing term $f$ are chosen such that the exact solution is
\begin{align}
u(\xx) &= 3 \left((x_1 - 0.2) + \frac{x_2 - 0.3}{2}\right)^2 + 2 \left( \frac{x_1 - 0.7}{2} + (x_2 - 0.8) \right)^3 +\sin(2 \pi x_1) \sin(3 \pi x_2),
     \label{eq:test1_exact_solution}
\end{align}
which is shown in Figure \ref{fig:test1_solutions}. We observe that this is the same test performed in \cite{PintoreTeora2024} which is now extended to new polygonal meshes. 



We solve problem \eqref{eq:diffreacadv} using both the NAVEM and VEM methods, and we plot the corresponding errors with respect to $h$ in Figure \ref{fig:test1_errors}. Since the desired solution is regular enough, the VEM $L^2$-error $\mathrm{err}_0^{\rm{VEM}}$
and $H^1$-error $\mathrm{err}_1^{\rm{VEM}}$ decrease with expected rates of 
$O(h^2)$ and $O(h^1)$, respectively. Empirical observations indicate that the NAVEM solution converges at the same rate as the VEM in both the 
$L^2$ and $H^1$ norms. We further note that these numerical results are coherent with the expected ones, since the values of the loss functions are small enough on the elements of the meshes. Indeed, we recall that the elements of our meshes are well-represented by the polygons in the training sets, as discussed previously. Additionally, the absence of the projection and stabilization operator in the NAVEM is manifested as a downward shift in the convergence curves.

\subsection{Test problem 2: Anisotropic problem}

\begin{figure}[!ht]
\centering
\begin{subfigure}{0.32\textwidth}
    \includegraphics[width=\textwidth]{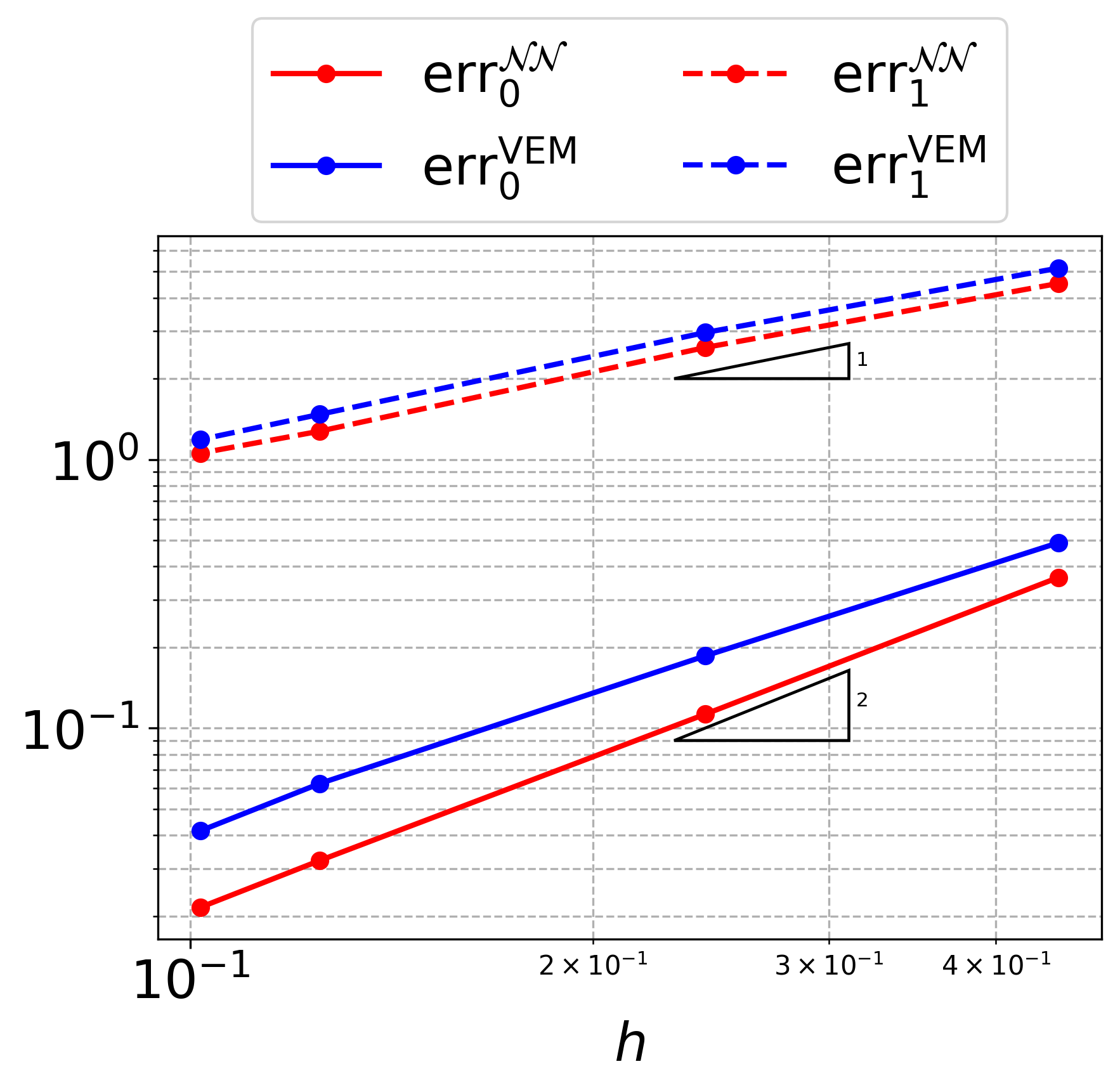}
    \caption{}
    \label{}
\end{subfigure}\quad
\begin{subfigure}{0.32\textwidth}
    \includegraphics[width=\textwidth]{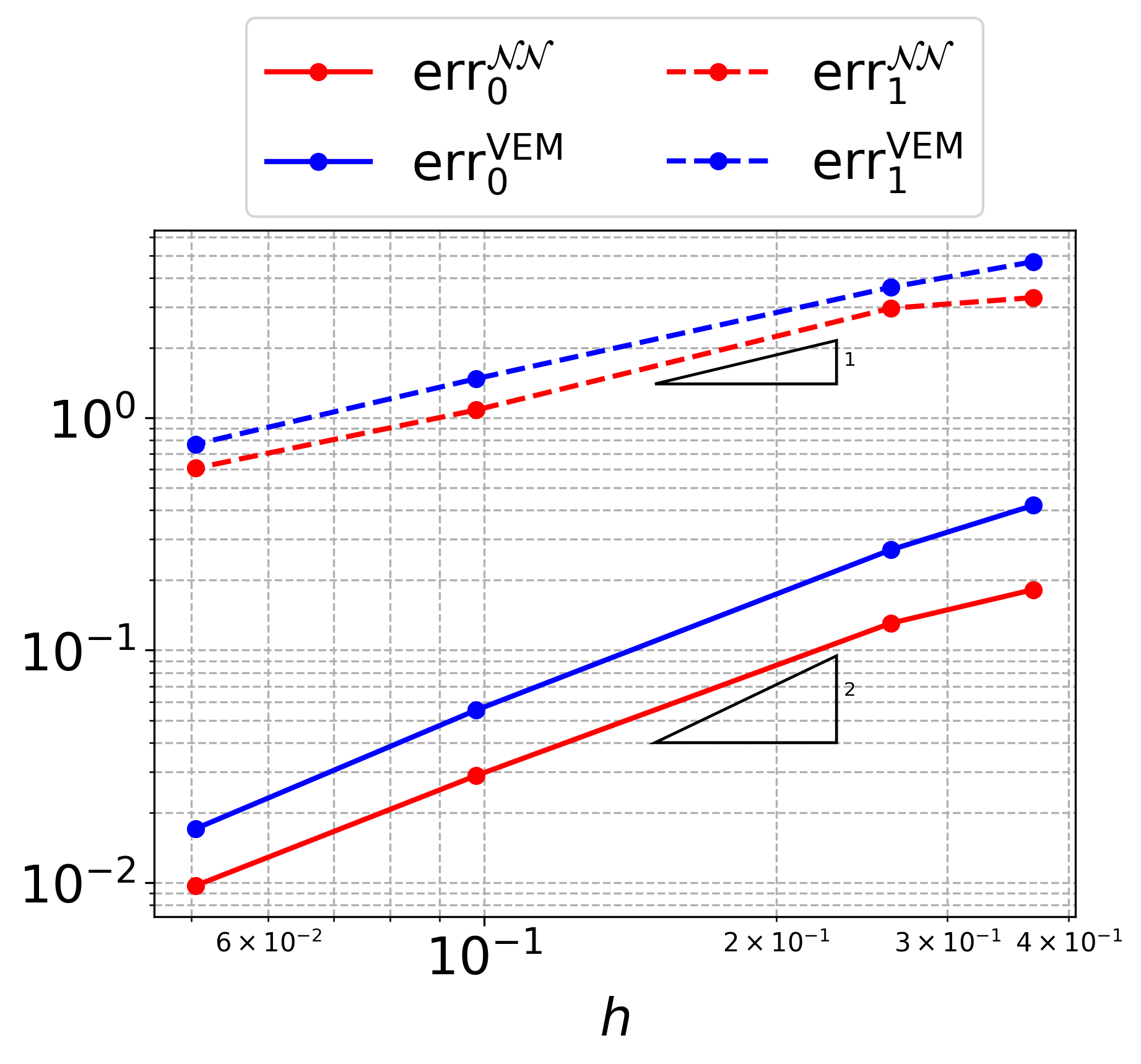}
    \caption{}
    \label{}
\end{subfigure}
\begin{subfigure}{0.32\textwidth}
    \includegraphics[width=\textwidth]{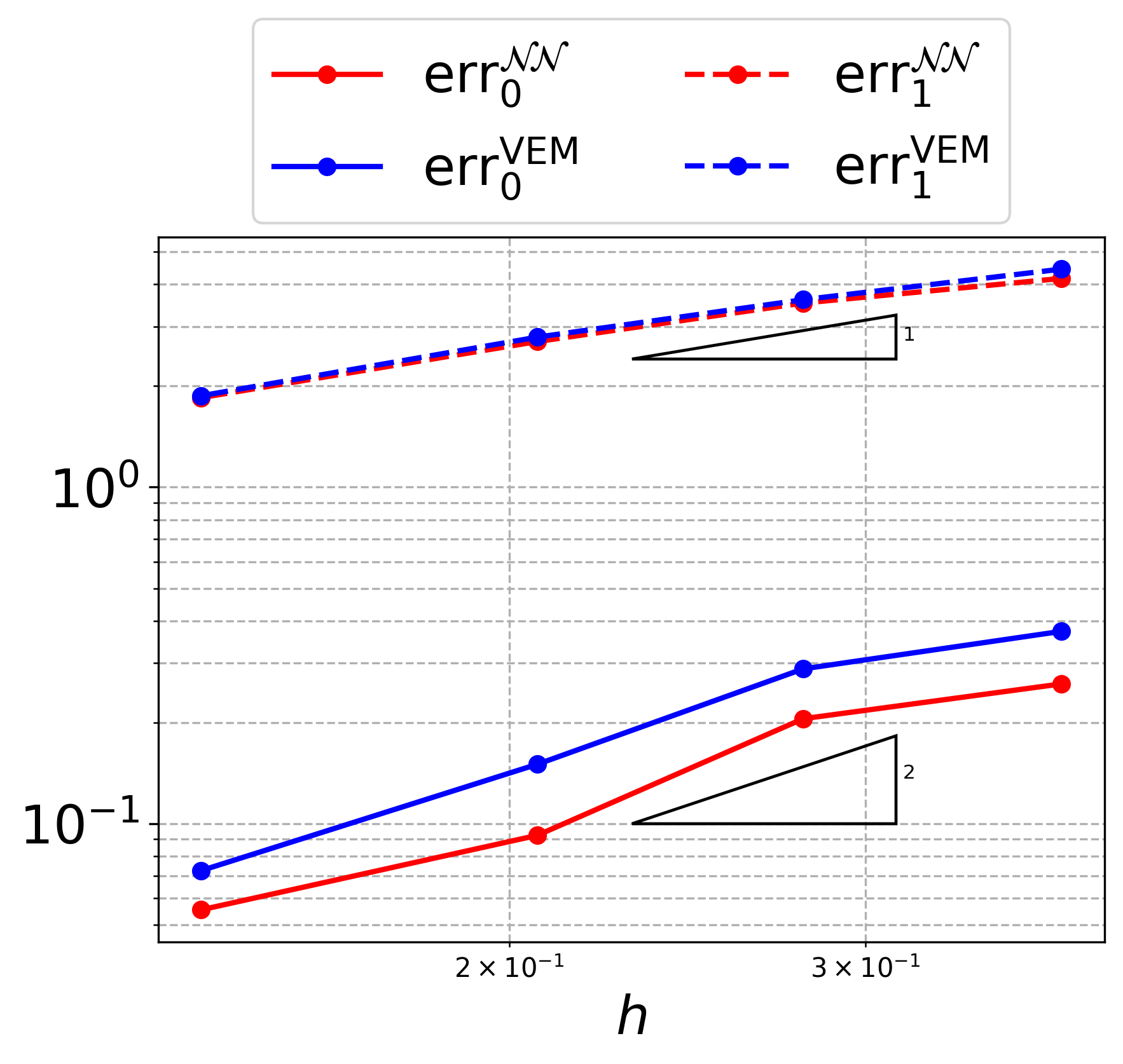}
    \caption{}
    \label{}
\end{subfigure}
\caption{Test 2: NAVEM and VEM errors w.r.t $h$. RDQM (Left), VM (Center) and HTM (Right).}
\label{fig:test2_errors}
\end{figure}

Let us now consider a boundary-value problem with a strongly anisotropic tensor $\DD$. In particular, we solve the problem \eqref{eq:diffreacadv} on $\Omega = (0,1)^2$ with
\begin{linenomath}
\begin{gather*}
    \DD = \mathbf{G}\widetilde\DD \mathbf{G}^T, \quad \widetilde \DD = \begin{bmatrix}
        1  & 0\\
        0 & 1.0e{-6}
    \end{bmatrix},\quad
    \mathbf{G} = \begin{bmatrix}
        \cos\left(\dfrac{\pi}6\right) & -\sin\left(\dfrac{\pi}6\right)\\[0.5cm]
        \sin\left(\dfrac{\pi}6\right) & \cos\left(\dfrac{\pi}6\right)
    \end{bmatrix}, \quad 
     \bbeta = \begin{bmatrix}
        0 \\ 0
    \end{bmatrix}, \quad \gamma = 0,
\end{gather*}
\end{linenomath}
where the matrix $\mathbf{G}$ is the Givens rotation matrix. As in the previous test cases, the forcing term and boundary data are chosen such that the exact solution is:
\begin{equation}
u(\xx) = \sin\left(3\cos\left(x-2y^2\right)^2+4\sin\left(y+2x\right)^2\right),
\label{eq:test2_exact_solution}
\end{equation}
which is shown in Figure \ref{fig:test2_solutions}.

The error convergence curves for both VEM and NAVEM are shown in Figure \ref{fig:test2_errors}. Again, we can observe that the knowledge of the virtual element basis functions is reflected in small error constants. We further highlight that, in this kind of problem, it is very difficult to design a proper stabilization term and a tuning strategy is advisable to choose a multiplicative stabilization coefficient \cite{Marcon2024}. The usage of a method like NAVEM removes this kind of issue.

\subsection{Test problem 3: Nonlinear problem}

\begin{figure}[t]
\centering
\includegraphics[width=0.5\textwidth]{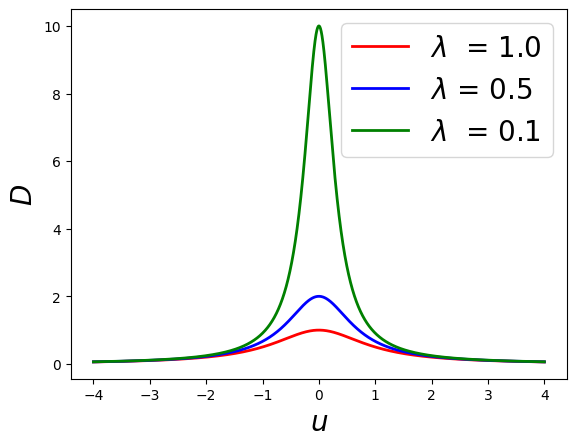}
\caption{Test 3: Diffusion coefficient $D(u, \lambda)$ as $u$ varies. The different curves are associated with different values of the parameter $\lambda$.}
\label{fig:test3_nl_coefficient}
\end{figure}

\begin{figure}[t]
\centering
\begin{subfigure}{0.32\textwidth}
    \includegraphics[width=\textwidth]{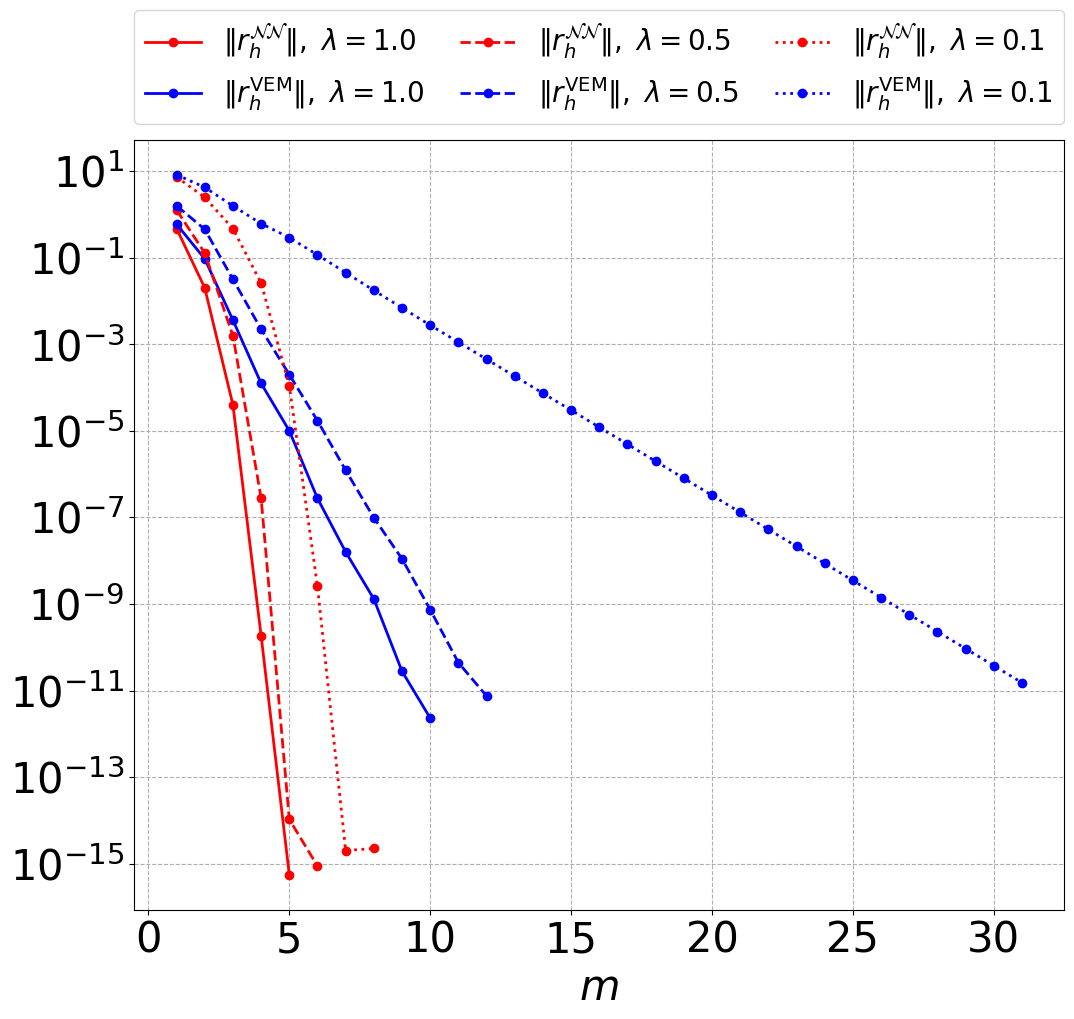}
    \caption{}
    \label{}
\end{subfigure}
\begin{subfigure}{0.32\textwidth}
    \includegraphics[width=\textwidth]{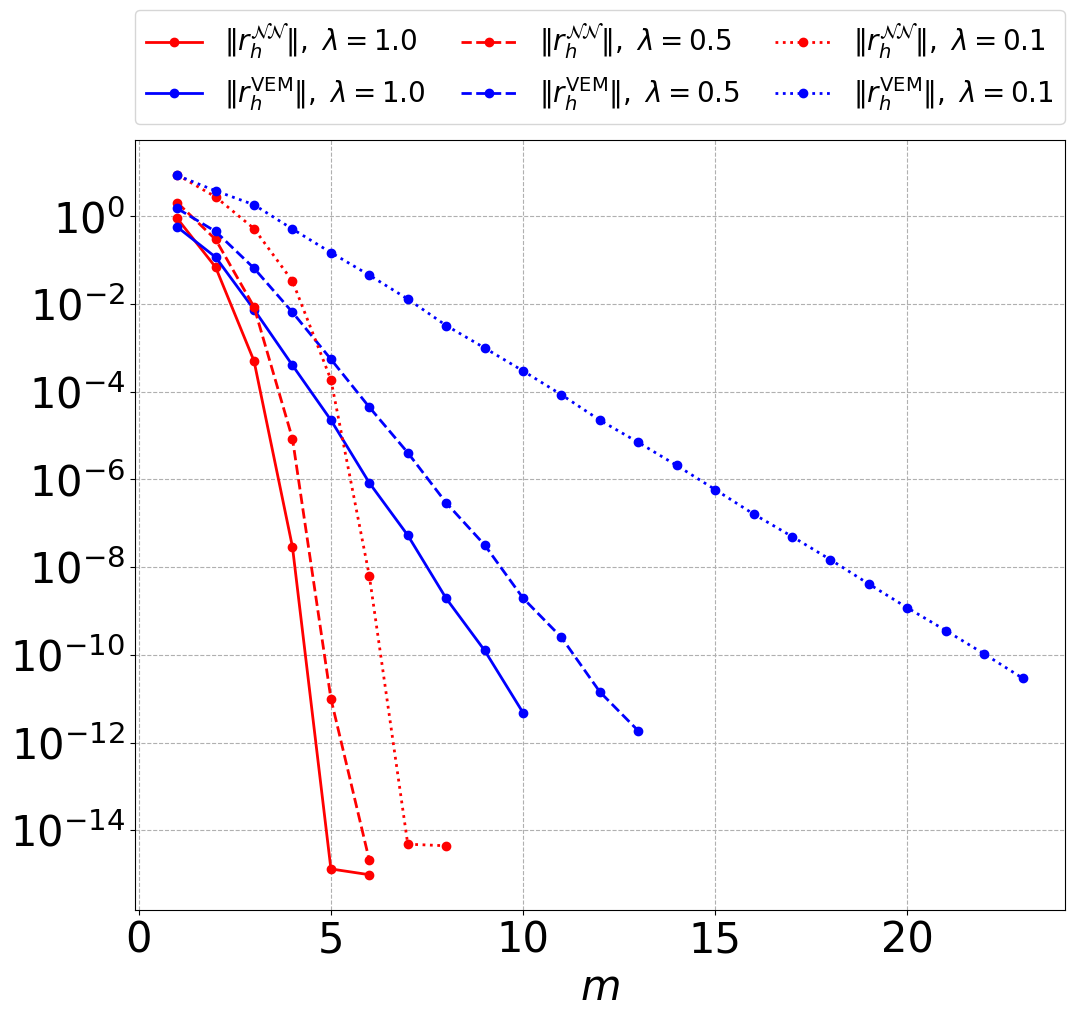}
    \caption{}
    \label{}
\end{subfigure}
\begin{subfigure}{0.32\textwidth}
    \includegraphics[width=\textwidth]{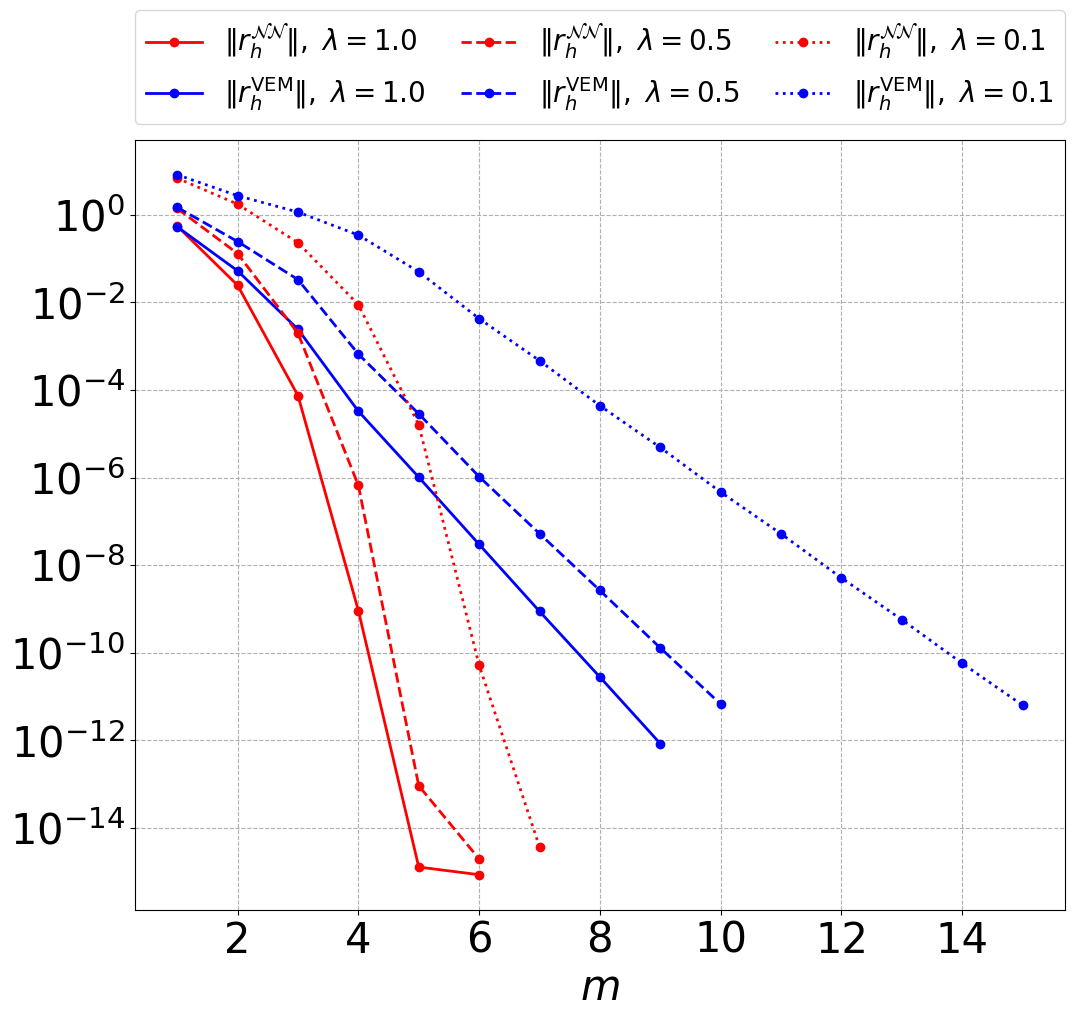}
    \caption{}
    \label{}
\end{subfigure}
\begin{subfigure}{0.32\textwidth}
    \includegraphics[width=\textwidth]{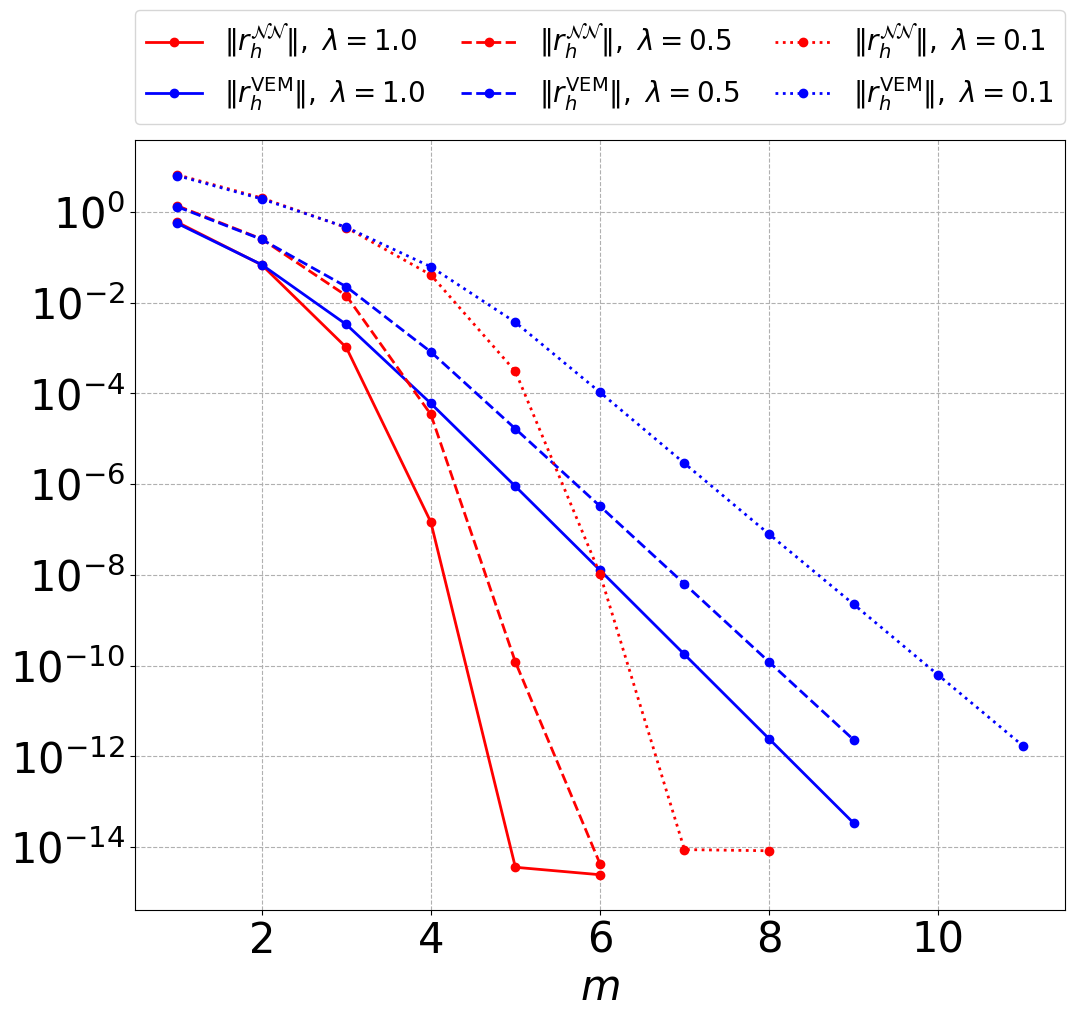}
    \caption{}
    \label{}
\end{subfigure}
\begin{subfigure}{0.32\textwidth}
    \includegraphics[width=\textwidth]{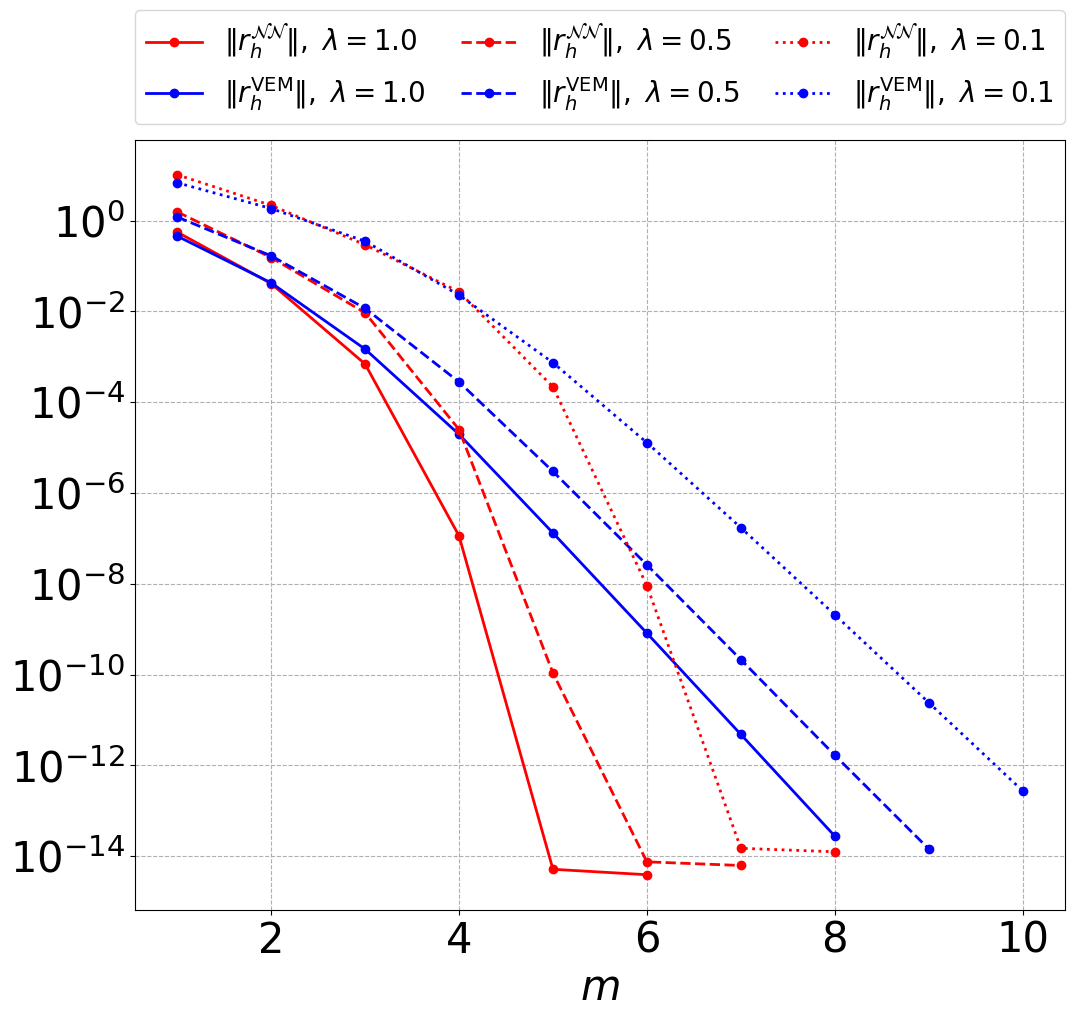}
    \caption{}
    \label{}
\end{subfigure}
\begin{subfigure}{0.32\textwidth}
    \includegraphics[width=\textwidth]{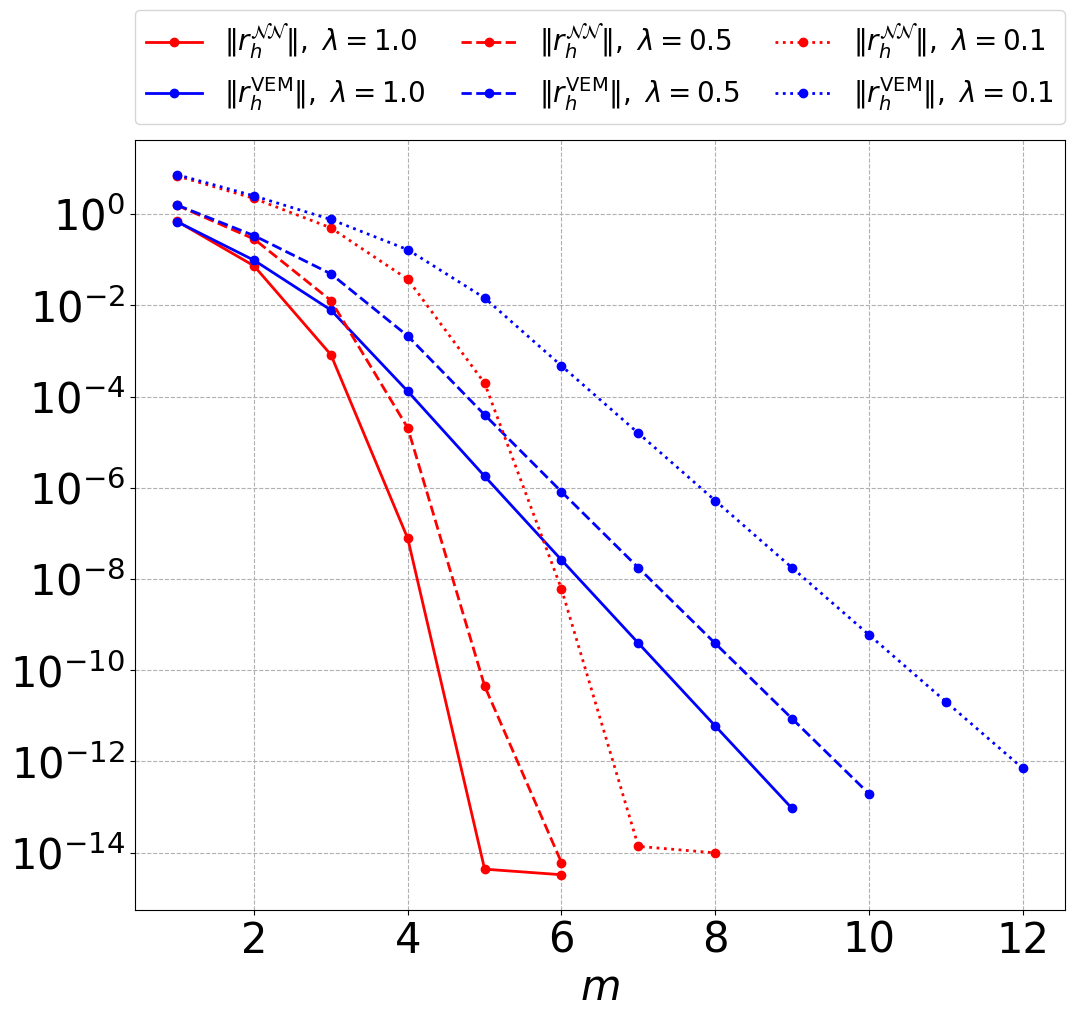}
    \caption{}
    \label{}
\end{subfigure}
\caption{Test 3: Residual norm as the number of non-linear iterations varies for both the VEM and the NAVEM. Different line styles represent different values of the parameter $\lambda$, while the two rows represent the first and the last mesh of each family, respectively. RDQM (Left), VM (Center) and HTM (Right).}
\label{fig:test3_nl_iterations}
\end{figure}

\begin{figure}[t]
\centering
\begin{subfigure}{0.32\textwidth}
    \includegraphics[width=\textwidth]{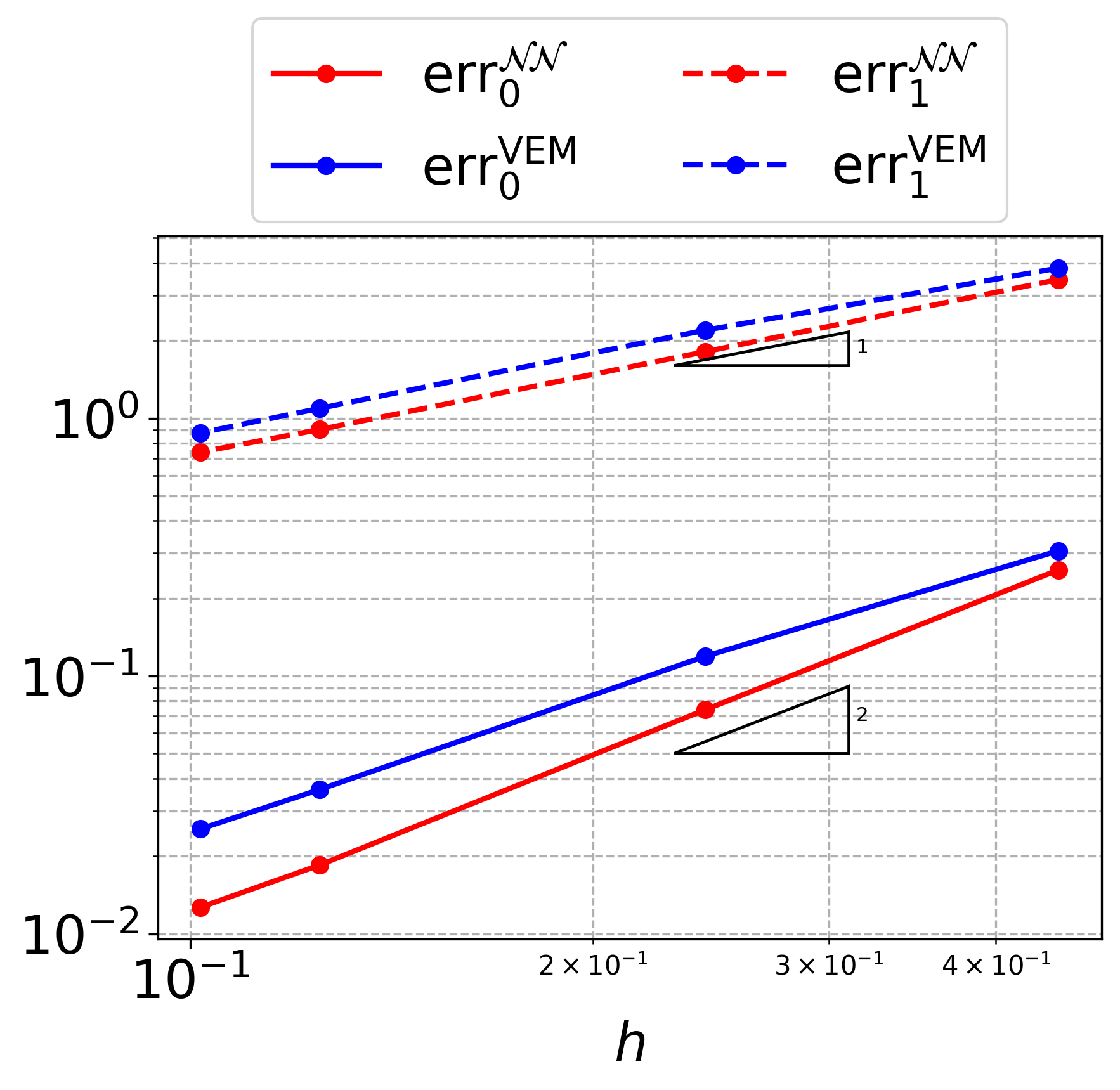}
    \caption{}
    \label{}
\end{subfigure}\quad
\begin{subfigure}{0.32\textwidth}
    \includegraphics[width=\textwidth]{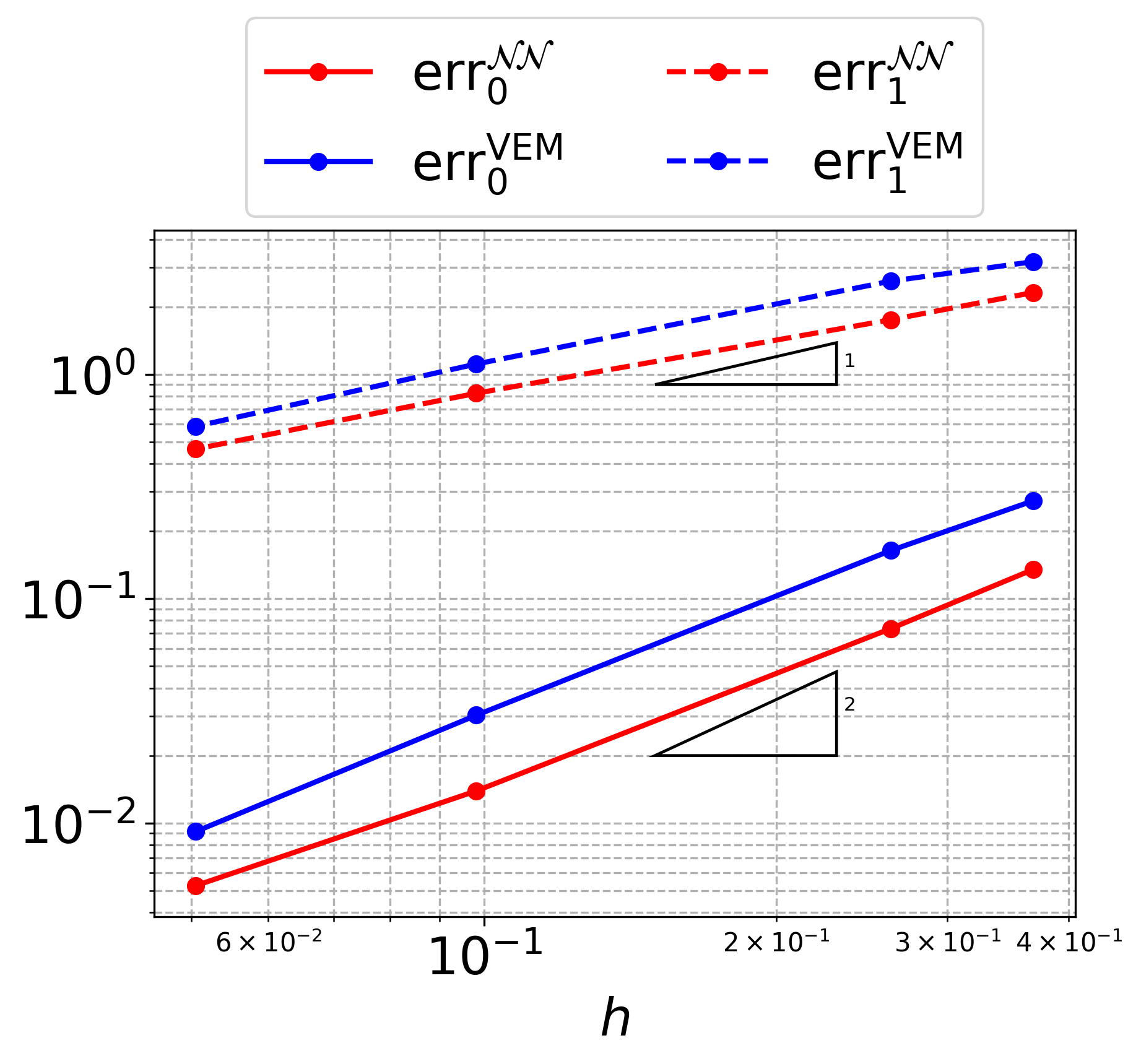}
    \caption{}
    \label{}
\end{subfigure}
\begin{subfigure}{0.32\textwidth}
    \includegraphics[width=\textwidth]{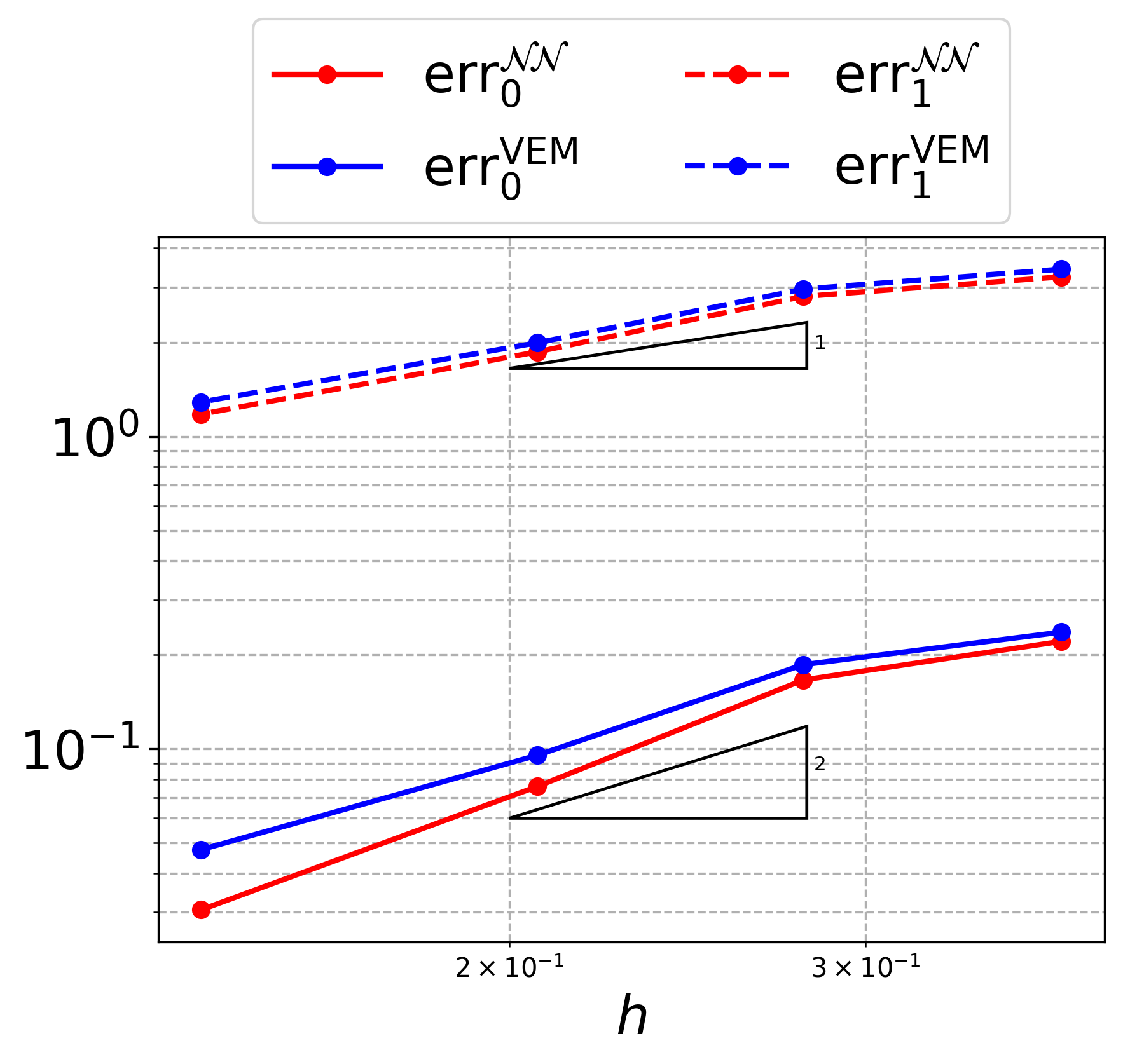}
    \caption{}
    \label{}
\end{subfigure}
\begin{subfigure}{0.32\textwidth}
    \includegraphics[width=\textwidth]{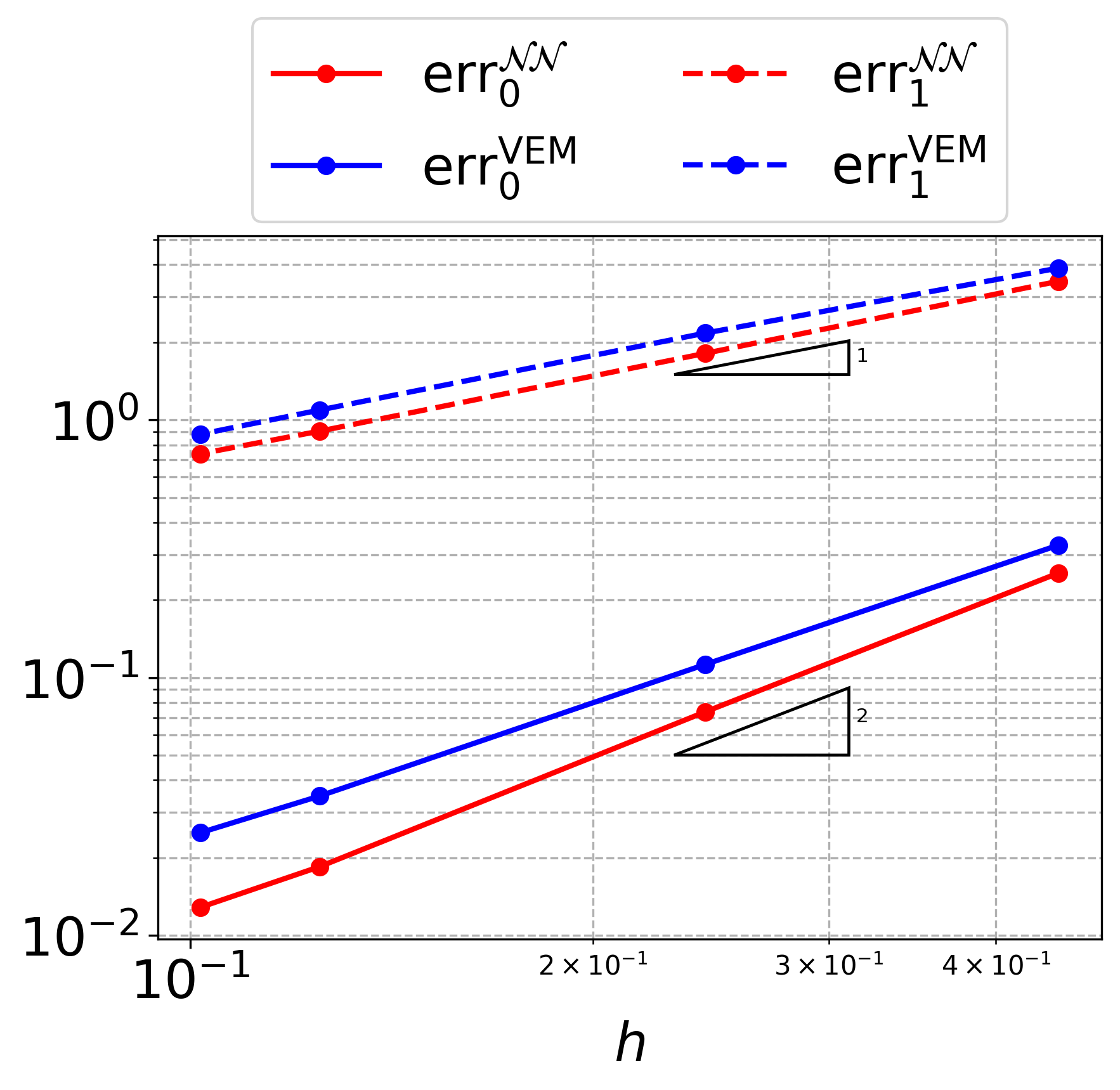}
    \caption{}
    \label{}
\end{subfigure}\quad
\begin{subfigure}{0.32\textwidth}
    \includegraphics[width=\textwidth]{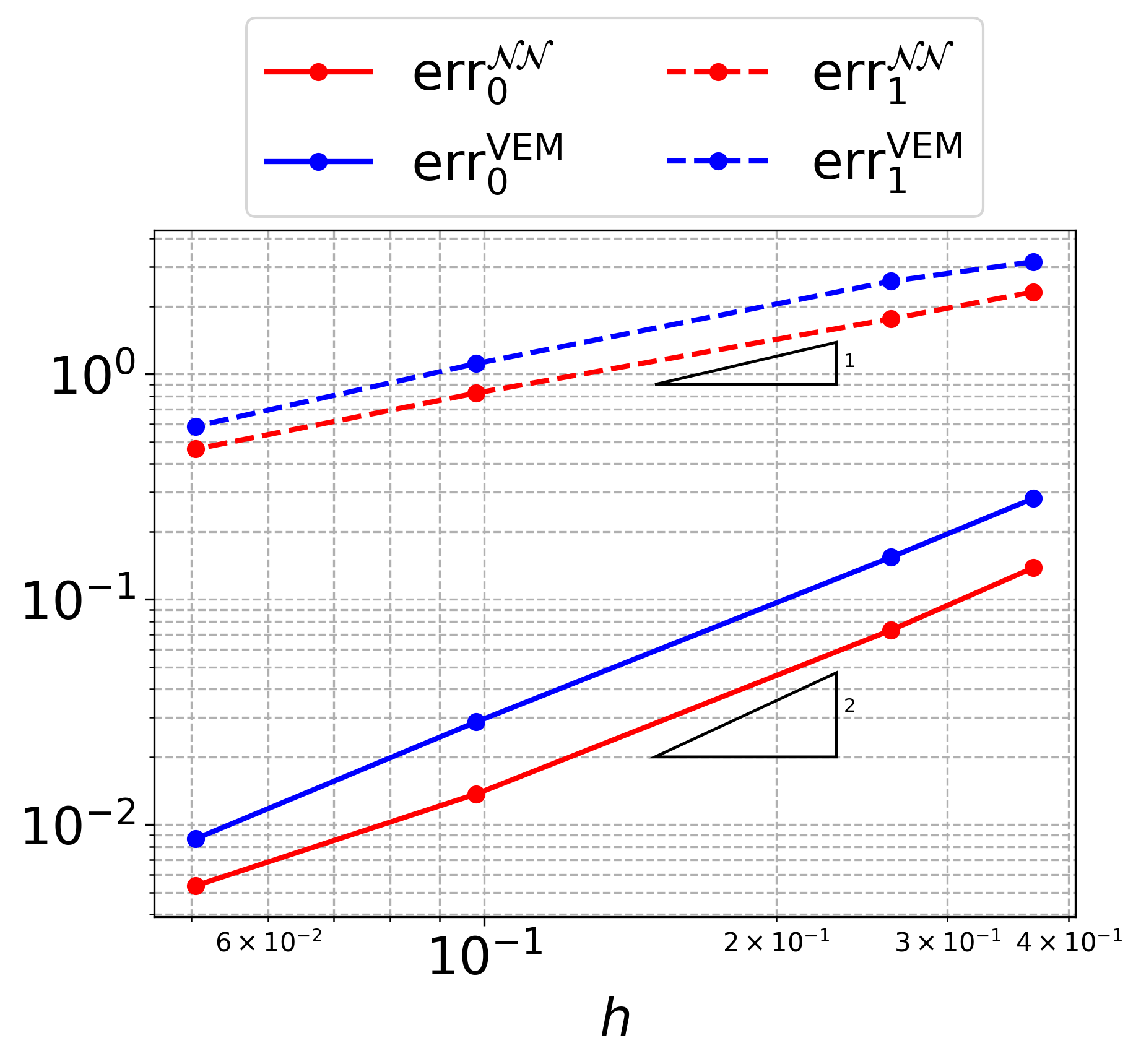}
    \caption{}
    \label{}
\end{subfigure}
\begin{subfigure}{0.32\textwidth}
    \includegraphics[width=\textwidth]{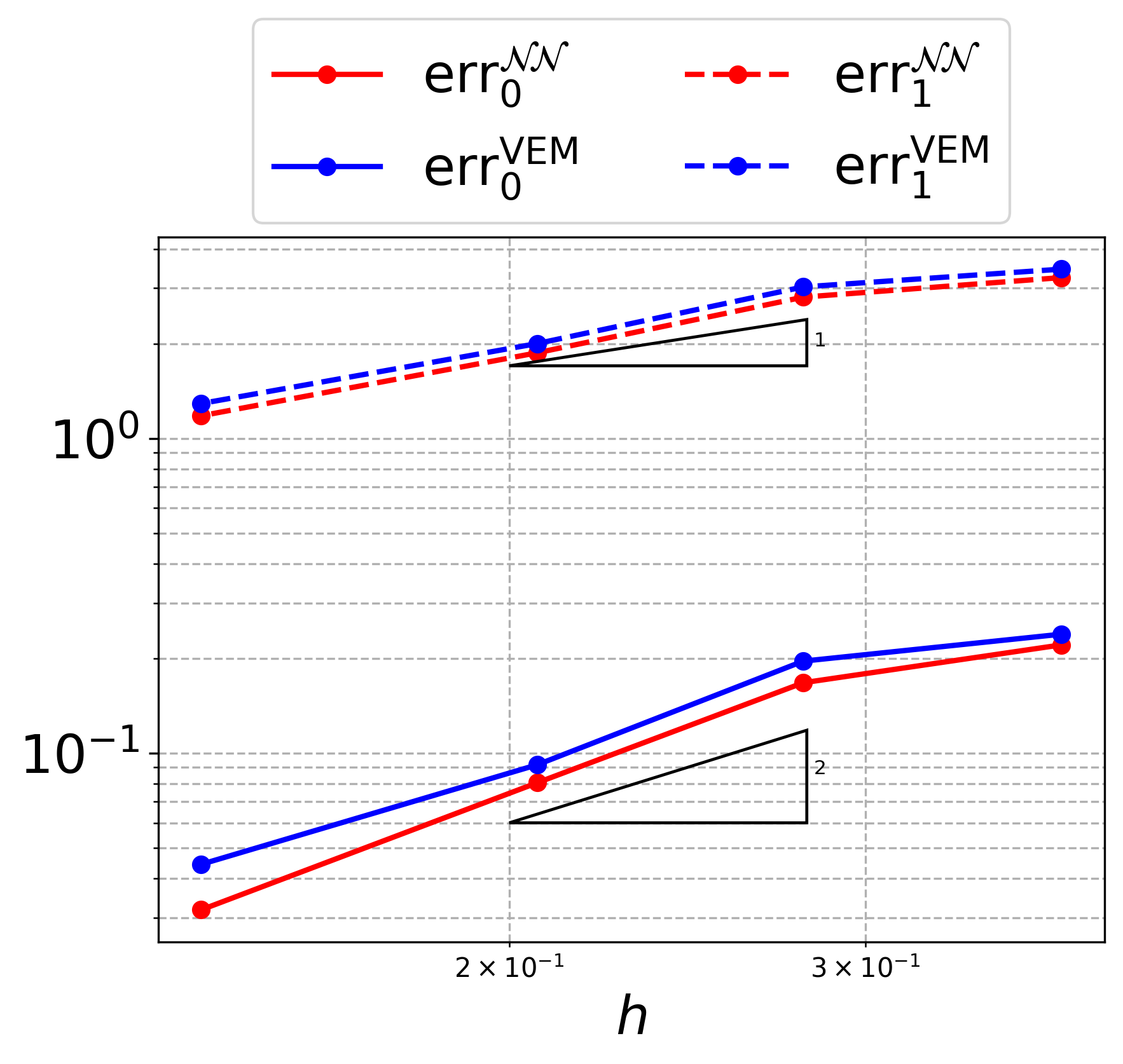}
    \caption{}
    \label{}
\end{subfigure}
\begin{subfigure}{0.32\textwidth}
    \includegraphics[width=\textwidth]{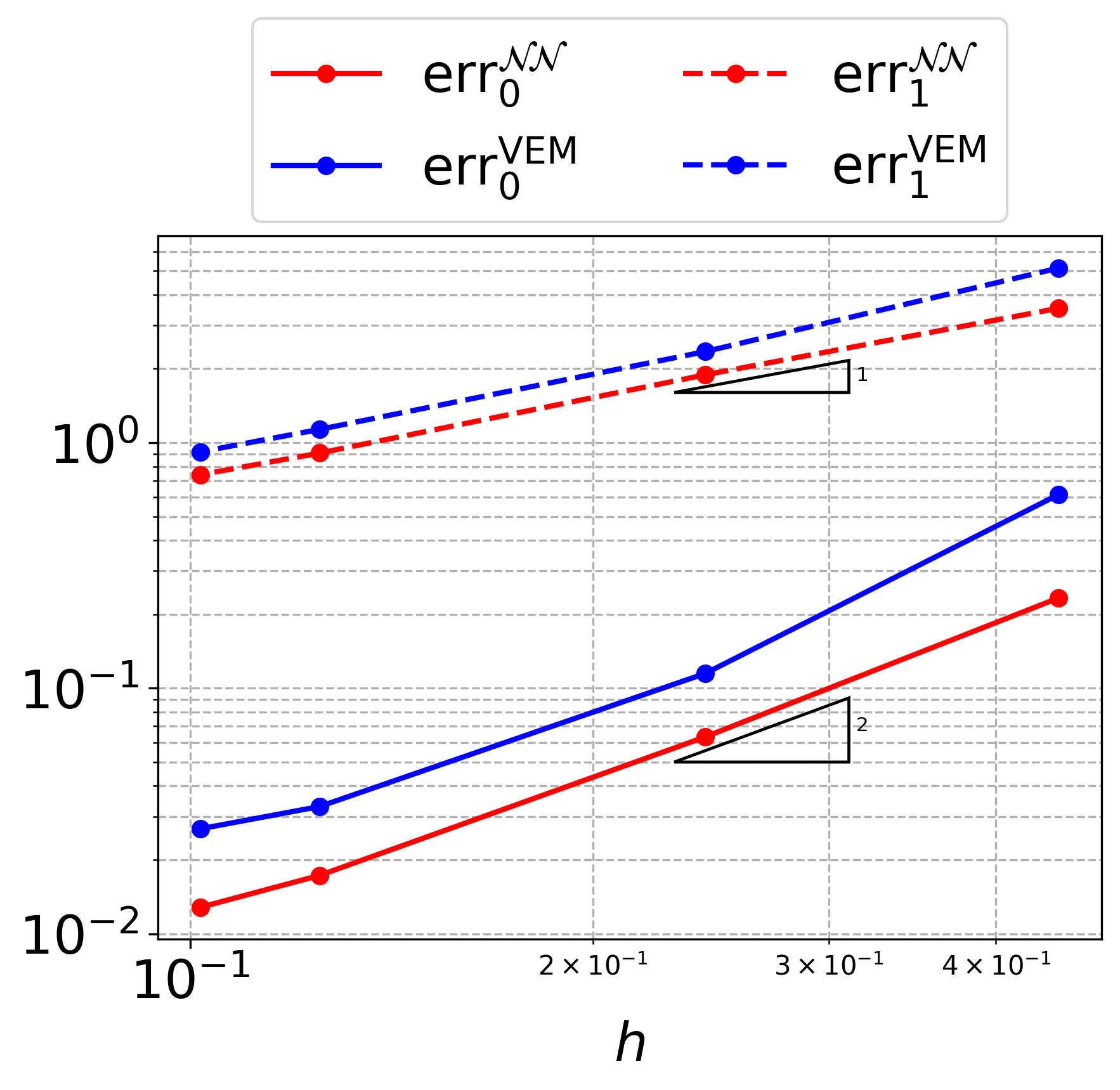}
    \caption{}
    \label{}
\end{subfigure}\quad
\begin{subfigure}{0.32\textwidth}
    \includegraphics[width=\textwidth]{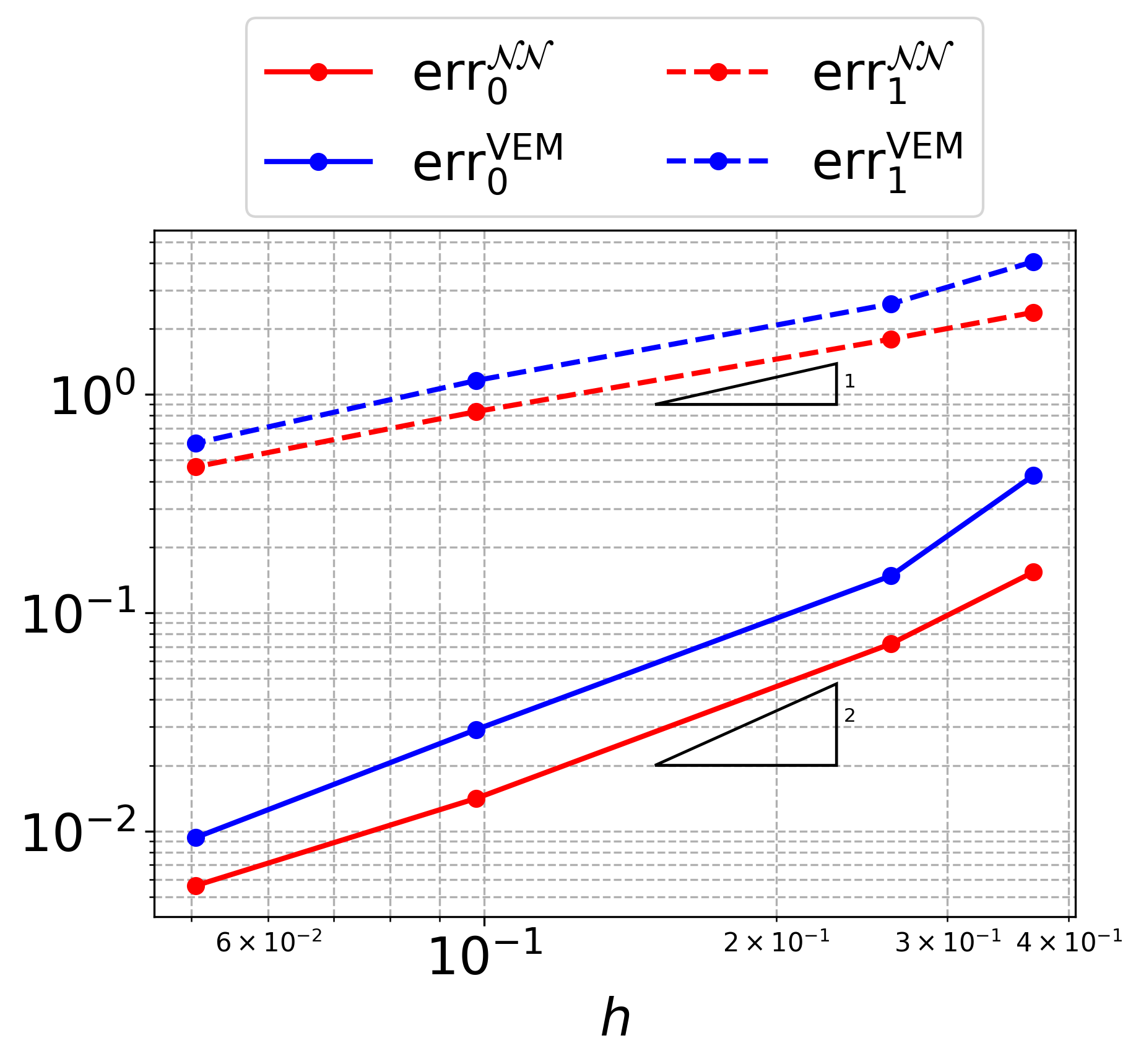} 
    \caption{}
    \label{}
\end{subfigure}
\begin{subfigure}{0.32\textwidth}
    \includegraphics[width=\textwidth]{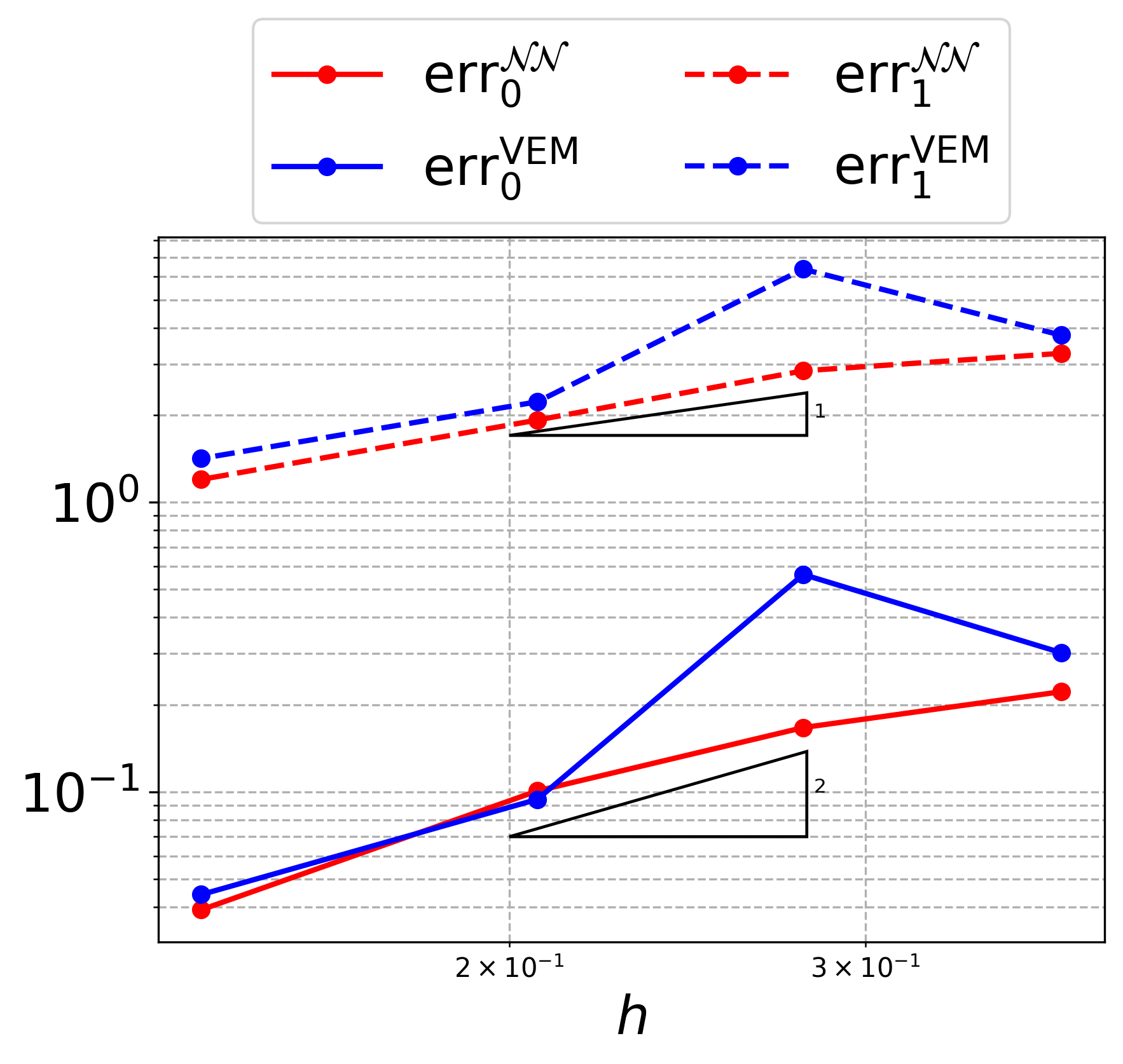}
    \caption{}
    \label{}
\end{subfigure}
\caption{Test 3: NAVEM and VEM errors w.r.t $h$. Each row corresponds to a different value of the parameter $\lambda = 1.0, 0.5, 0.1$ from top to bottom. RDQM (Left), VM (Center) and HTM (Right).}
\label{fig:test3_errors}
\end{figure}

Let us now consider the following nonlinear problem
\begin{equation}
    \begin{cases}
        - \nabla \cdot \left( D(u, \lambda) \nabla u\right)  = f & \text{in } \Omega,\\
        u = g_D & \text{on } \Gamma,
    \end{cases}
    \label{eq:nonlinear_prob}
\end{equation}
where the diffusion coefficient is given by
\begin{equation}
    D(u, \lambda) = \frac{1}{\lambda + u^2}, \text{ with } \lambda =1.0,\, 0.5,\, 0.1,
\end{equation}
whose graphical illustration is provided in Figure \ref{fig:test3_nl_coefficient} for $\lambda =1.0,\, 0.5,\, 0.1$. In order to compute the errors in \eqref{eq:nn_errors} and \eqref{eq:vem_errors}, we choose the Dirichlet boundary condition and the forcing term in such a way that, for any $\lambda$, the exact solution is
\begin{equation}
    u(\xx) = \frac{1}{8} \left(\sin(3 \pi ((x_1-0.5)^2 + (x_2 - 0.5)^2))\right)^3.
\label{eq:test3_exact_solution}
\end{equation}
The contour plot of the exact solution \eqref{eq:test3_exact_solution} is shown in Figure \ref{fig:test3_solutions}.

For comparison purposes, we briefly report here the Virtual Element formulation that we use to solve the problem \eqref{eq:nonlinear_prob}, which is introduced in \cite{Cangiani2019}. In the case of quasilinear elliptic problems, the local discrete virtual element bilinear form $\dbilinh[E]{}{}$ read as
\begin{linenomath}
\begin{align*}
    \dbilinh[E]{u_h}{v_h; z_h, \lambda} &= \int_E D(\proj{E,0}{1} z_h, \lambda) \proj{E,0}{0}\nabla u_h \cdot \proj{E,0}{0} \nabla v_h \\
    &\quad + \stab[E]{(I-\proj{E,\nabla}{1})u_h}{(I-\proj{E,\nabla}{1})v_h; z_h, \lambda},
\end{align*}
\end{linenomath}
where the VEM stabilizing form $\stab[E]{\cdot}{\cdot\ ; z_h, \lambda}$ is given by
\begin{linenomath}
\begin{equation*}
    \stab[E]{u_h}{v_h; z_h, \lambda} = D(\proj{E,0}{0} z_h, \lambda ) \sum_{i=1}^{\Ndof[E]} \dof_i^E(u_h) \dof_i^E(v_h).
\end{equation*}
\end{linenomath}
Now, we apply the Newton-Raphson method to deal with nonlinearities. Thus, given an initial iterates $u_h^0 \in \Vh{1}$, we define a sequence
\begin{linenomath}
\begin{equation*}
    u_h^{m + 1} = u_h^m + \delta_h^m\quad \forall m \geq 0,
\end{equation*}
\end{linenomath}
by solving at each nonlinear step $m$ the linearized problem: \textit{Find} $\delta_h^m \in \Vh{1}$ \textit{such that:}
\begin{linenomath}
\begin{align*}
    \sum_{E \in \Th} \left[\dbilinh[E]{\delta_h^m}{v_h; u_h^m, \lambda} + \jacbilinh[E]{\delta_h^m}{v_h; u_h^m, \lambda}\right] =  \sum_{E \in \Th}\left[ \rightlinh[E]{v_h} - \dbilinh[E]{u_h^m}{v_h; u_h^m, \lambda}\right]\quad \forall v_h \in \Vh{1},
\end{align*}
\end{linenomath}
where the extra term $\jacbilinh[E]{\cdot}{\cdot \ ; u_h^m, \lambda}$ stems from the linearization of both the consistency and the stabilization term and it is defined as
\begin{linenomath}
\begin{align*}
    \jacbilinh[E]{\delta_h^m}{v_h; u_h^m, \lambda} &= \int_E \frac{\partial D(\proj{E,0}{1} u_h^m, \lambda)}{\partial u} \proj{E,0}{1} \delta_h^m \proj{E,0}{0} \nabla u_h^m \cdot \nabla v_h \\
    &+ \frac{\partial D(\proj{E,0}{0} u_h^m, \lambda)}{\partial u} \proj{E,0}{0} \delta_h^m \sum_{i=1}^{\Ndof[E]} \dof_i^E((I-\proj{E,\nabla}{1})u_h) \dof_i^E((I-\proj{E,\nabla}{0})v_h).
\end{align*}
\end{linenomath}
Now, let us denote by $\bm{\delta}_h^m$ and $\bm{u}_h^m$ the vectors of coefficients of functions $\delta_h^m$ and $u_h^m$ with respect to the virtual element basis functions, we introduce the following matrices
\begin{linenomath}
\begin{align*}
    \mathbf{A}^{m}_h \in \R^{\Ndof \times \Ndof}: \quad (\mathbf{A}^{ m}_h)_{ji} = \sum_{E \in \Th} \left[\dbilinh[E]{\varphi_i}{\varphi_j; u_h^m, \lambda}\right],\\
    \bm{f}^{m}_h \in \R^{\Ndof}: \quad (\bm{f}^{m}_h)_{j} = \sum_{E \in \Th}\left[ \rightlinh[E]{\varphi_j} - \dbilinh[E]{u_h^m}{\varphi_j; u_h^m, \lambda}\right],
\end{align*}
\end{linenomath}
and define as stopping criteria
\begin{linenomath}
\begin{equation}
\begin{gathered}
    \norm[2]{\bm{r}_h^{VEM}} \leq \epsilon_{r,r} \norm[2]{\bm{f}^{0}_h} + \epsilon_{r,a} \text{ and }\\
    \norm[2]{\bm{\delta}_h^m} \leq \epsilon_{\delta,r} \norm[2]{\bm{u}_h^0} + \epsilon_{\delta,a},
    \end{gathered}
    \label{eq:stopping_criteria}
\end{equation}
\end{linenomath}
where
\begin{linenomath}
\begin{equation*}
    \bm{r}_h^{VEM} = \bm{f}^{m}_h - \mathbf{A}^{m}_h \bm{u}_h^m \quad \forall m \geq0.
\end{equation*}
\end{linenomath}
We observe that, in this kind of construction, the heavy usage of the polynomial projectors and stabilization term could become a great issue when the diffusion coefficient becomes highly non-linear, since the non-linearity may increase the distance between the virtual element solution and its projection. The usage of the NAVEM method helps to get rid of any stabilization term or projection operator, simplifying the discrete bilinear form used in the Newton-Raphson method. Given an initial guess $u_h^{\NN, 0} \in \nVh{1}$, at each nonlinear step $m=0,\dots,$ we solve the linearized problem: \textit{Find} $\delta_h^{\NN, m} \in \nVh{1}$ \textit{such that:}
\begin{linenomath}
\begin{align*}
    &\sum_{E \in \Th} \left[ \dbilinhn[E]{{\delta}^{\NN,m}_h}{v^{\NN}_h; \ u_h^{\NN, m}, \lambda}  + \jacbilinhn[E]{\delta^{\NN,m}_h}{v^{\NN}_h; \ u_h^{\NN,m}, \lambda} \right] = \\
    &\qquad \qquad \sum_{E \in \Th}\left[ \scal[E]{f}{v^{\NN}_h} - \dbilinhn[E]{u_h^{\NN,m}}{v^{\NN}_h; u_h^{\NN,m}, \lambda}\right]\quad \forall v^{\NN}_h \in \nVh{1},
\end{align*}
\end{linenomath}
where
\begin{linenomath}
\begin{equation*}
    u^{\NN, m+1}_h = \delta^{\NN,m}_h + u^{\NN,m}_h \quad m = 0,\dots.
\end{equation*}
\end{linenomath}
and the involved bilinear forms are reduced to
\begin{linenomath}
\begin{gather*}
    \dbilinhn[E]{u^{\NN}_h}{v^{\NN}_h; \ z_h^{\NN}, \lambda} = \int_E D(z_h^{\NN}, \lambda) \nabla u^{\NN}_h \cdot \nabla v^{\NN}_h,\\
    \jacbilinhn[E]{\delta^{\NN}_h}{v^{\NN}_h; \ u_h^{\NN}, \lambda} = \int_E \frac{\partial D(u_h^{\NN}, \lambda)}{\partial u}  \delta_h^{\NN} \ \nabla u_h^{\NN} \cdot \nabla v^{\NN}_h.
\end{gather*}
\end{linenomath}
In this numerical experiment, we set the initial guess as the all-zeros vector and we choose $\epsilon_{r,r} = \epsilon_{r,a} = 10^{-12}$ and $\epsilon_{\delta,r} = \epsilon_{\delta,a} = 10^{-10}$.

Figure \ref{fig:test3_nl_iterations} shows the behaviours of the residual as the number of non-linear iterations increases for both the methods, for the values of the parameter $\lambda = 1.0,\, 0.5,\, 0.1$ and for the first and the last mesh of each family. We observe that the coarser is the mesh and the smaller is $\lambda$, the larger is the number of iterations that the standard VEM employs to reach the desired tolerance. On the other hand, we observe that the number of iterations related to NAVEM is not strongly dependent on the parameter $h$ and its variability with respect to the parameter $\lambda$ is much weaker with respect to VEM. Moreover, we highlight that the plateau of the residual is due to the double stopping criteria imposed \eqref{eq:stopping_criteria}.

Figure \ref{fig:test3_errors} shows the convergence curves related to the errors \eqref{eq:nn_errors} and \eqref{eq:vem_errors} for NAVEM and VEM, respectively. Again, we observe a reduction in the error constants for each tested case with respect to VEM method, while highly reducing the number of iterations needed to achieve the desired tolerance. These results suggest that NAVEM can provide competitive accuracy while simplifying the formulation.

\section{Conclusions}\label{sec:conclusion} 
In this paper, we extend and describe the lowest-order Neural Approximated Virtual Element Method (NAVEM) on quite general polygonal elements. The NAVEM is a polygonal method used to solve partial differential equations which combines standard numerical techniques with neural networks, preserving the convergence rate of the standard numerical method, while exploiting the offline-online paradigm of neural networks to overcome the limitations of the standard procedure.
Indeed, it modifies the original VEM formulation in \cite{LBe13} by explicitly approximating the virtual element basis functions through suitable harmonic functions parameterized by a neural network and deleting issues related to the introduction of polynomial projections and stabilization operators as in standard VEM. 

Two different neural network architectures and related training strategies are described and theoretically justified. Few papers tackle theoretical discussions and results about the intersection between standard mesh-based solvers and neural network are available; some examples including \cite{Pintore2022, badia2024finite, berrone2022solving}.

Numerical results confirm that the presented method helps to avoid issues concerning the choice of the stabilization term and of accessing to the point-wise evaluation of basis functions without using polynomial projectors, showing good performances, especially in the case of highly non-linear problems. Furthermore, particular attention is devoted to the analysis of triangular meshes with hanging nodes given their relevance in the numerical field. 

We believe that this study could help exploring new advanced strategies for practical applications.